\documentclass{elsarticle}

\usepackage{bm}
\usepackage{fullpage}
\usepackage{amsfonts}
\usepackage{graphicx}
\usepackage{amsthm}
\usepackage{amssymb}
\usepackage{paralist}
\usepackage{mathtools}
\usepackage[section]{placeins}
\usepackage{algorithm}
\usepackage{algorithmic}
\usepackage{moreverb}
\usepackage[colorlinks, bookmarksopen, bookmarksnumbered, citecolor=red,
            urlcolor=red]{hyperref}
\usepackage{subcaption}
\usepackage[table]{xcolor}

\usepackage{tikz}
\usepackage{pgfplots}
\usepackage{pgfplotstable, booktabs}

\pgfplotsset{select coords between index/.style 2 args={
    x filter/.code={
        \ifnum\coordindex<#1\fi
        \ifnum\coordindex>#2\fi
    }
}}

\usepackage{bm}

\newcommand{\beq}{\begin{equation}}
\newcommand{\eeq}{\end{equation}}
\newcommand{\esplit}{\end{split}}
\newcommand{\beqalign}{\begin{array}{rl}}
\newcommand{\eeqalign}{\end{array}}

\newcommand{\norm}[1]{\ensuremath{|| #1 ||}}

\newcommand{\pder}[2]{\ensuremath{\frac{\partial #1}{\partial #2}}} 

\newcommand{\oder}[2]{\ensuremath{\frac{\mathrm{d} #1}{\mathrm{d} #2}}} 

\newcommand{\Fcal}{\ensuremath{\mathcal{F}}}
\newcommand{\Gcal}{\ensuremath{\mathcal{G}}}

\newcommand{\Jcal}{\ensuremath{\mathcal{J}}}

\newcommand{\Ocal}{\ensuremath{\mathcal{O}}}
\newcommand{\Pcal}{\ensuremath{\mathcal{P}}}

\newcommand{\Wcal}{\ensuremath{\mathcal{W}}}

\newcommand{\Lboldcal}{\ensuremath{\boldsymbol{\mathcal{L}}}}

\newcommand{\Mbb}{\ensuremath{\mathbb{M} }}

\newcommand{\Rbb}{\ensuremath{\mathbb{R} }}

\newcommand\Abm{{\ensuremath{\bm{A}}}}

\newcommand\Fbm{{\ensuremath{\bm{F}}}}
\newcommand\Gbm{{\ensuremath{\bm{G}}}}
\newcommand\Hbm{{\ensuremath{\bm{H}}}}
\newcommand\Ibm{{\ensuremath{\bm{I}}}}
\newcommand\Jbm{{\ensuremath{\bm{J}}}}

\newcommand\Qbm{{\ensuremath{\bm{Q}}}}
\newcommand\Rbm{{\ensuremath{\bm{R}}}}

\newcommand\Ubm{{\ensuremath{\bm{U}}}}
\newcommand\Vbm{{\ensuremath{\bm{V}}}}

\newcommand\Xbm{{\ensuremath{\bm{X}}}}

\newcommand\bbm{{\ensuremath{\bm{b}}}}
\newcommand\cbm{{\ensuremath{\bm{c}}}}

\newcommand\ebm{{\ensuremath{\bm{e}}}}
\newcommand\fbm{{\ensuremath{\bm{f}}}}

\newcommand\kbm{{\ensuremath{\bm{k}}}}

\newcommand\rbm{{\ensuremath{\bm{r}}}}

\newcommand\ubm{{\ensuremath{\bm{u}}}}
\newcommand\vbm{{\ensuremath{\bm{v}}}}
\newcommand\wbm{{\ensuremath{\bm{w}}}}
\newcommand\xbm{{\ensuremath{\bm{x}}}}
\newcommand\ybm{{\ensuremath{\bm{y}}}}

\newcommand\ebold{\ensuremath{\mathbf{e}}}

\newcommand\lambdabold{{\ensuremath{\boldsymbol{\lambda}}}}

\newcommand\mubold{{\ensuremath{\boldsymbol{\mu}}}}
\newcommand\kappabold{{\ensuremath{\boldsymbol{\kappa}}}}

\newcommand\taubold{{\ensuremath{\boldsymbol{\tau}}}}

\newcommand\nubold{{\ensuremath{\boldsymbol{\nu}}}}

\newcommand\Gammabold{{\ensuremath{\boldsymbol{\Gamma}}}}

\begin{document}

\title{A Fully Discrete Adjoint Method for Optimization of Flow Problems on
       Deforming Domains with Time-Periodicity Constraints}

\author[rvt1]{M.~J.~Zahr\fnref{fn1}\corref{cor1}}
\ead{mzahr@stanford.edu}

\author[rvt2]{P.-O.~Persson\fnref{fn2}}
\ead{persson@berkeley.edu}

\author[rvt2]{J.~Wilkening\fnref{fn2}}
\ead{wilken@math.berkeley.edu}
 
\address[rvt1]{Institute for Computational and Mathematical Engineering,
               Stanford University, Stanford, CA 94035.}
\address[rvt2]{Department of Mathematics and Lawrence Berkeley National Laboratory, University of California,
               Berkeley, CA 94720-3840.}
\cortext[cor1]{Corresponding author}
\fntext[fn1]{Graduate Student, Institute for Computational and Mathematical
             Engineering, Stanford University}
\fntext[fn2]{Associate Professor, Department of Mathematics, University of
             California, Berkeley.}

\begin{abstract}
A variety of shooting methods for computing fully discrete
time-periodic solutions of partial differential equations, including
Newton-Krylov and optimization-based methods, are discussed and used
to determine the periodic, compressible, viscous flow around a 2D
flapping airfoil. The Newton-Krylov method uses matrix-free GMRES to
solve the linear systems of equations that arise in the nonlinear
iterations, with matrix-vector products computed via the linearized
sensitivity evolution equations. The adjoint method is used to compute
gradients for the gradient-based optimization shooting methods. The
Newton-Krylov method is shown to exhibit superior convergence to the
optimal solution for these fluid problems, and fully leverages
quality starting data.

The central contribution of this work is the derivation of the
adjoint equations and the corresponding adjoint method for fully discrete,
time-periodically constrained partial differential equations. These adjoint
equations constitute a linear, two-point boundary value problem that is
provably solvable.  The periodic adjoint method is used to compute gradients
of quantities of interest along the manifold of time-periodic solutions of the
discrete partial differential equation, which is verified against a
second-order finite difference approximation. These gradients are then used in
a gradient-based optimization framework to determine the energetically optimal
flapping motion of a 2D airfoil in compressible, viscous flow over a single
cycle, such that the time-averaged thrust is identically zero. In less than
$20$ optimization iterations, the flapping energy was reduced nearly an order
of magnitude and the thrust constraint satisfied to $5$ digits of accuracy.
\end{abstract}

\maketitle

\section{Introduction}\label{sec:intro}
Cyclic steady-state motion of a system, i.e.~stable time-periodic
behavior, is of central importance in bio-locomotion and many branches
of engineering. Examples include flapping flight \cite{wang2000vortex,
wang2004unsteady, soueid2009optimization},
swimming at low or high Reynolds number
\cite{hosoi, loheac2013controllability, trouilloud2008soft,
      eloy2011optimisation},
helicopter aerodynamics \cite{dugundji1983some, hwang1997frequency,
                                    peters1994fast,
                                    verdult2004identification},
turbomachinery \cite{bucher2005detecting, allen2006floquet},
wind turbines \cite{bir1999operating, stol2002floquet} and vehicle tires
with treads \cite{oden:lin,tires1,tires2}, to name a few. A number of
sophisticated algorithms have recently been developed to compute
cyclic steady-states of systems governed by partial differential
equations. However, in applications, one often wishes to optimize a
quantity of interest over a cycle, such as minimizing energy subject to
lift and thrust constraints. The goal of the present paper is to
develop adjoint-based optimization techniques for such systems,
focusing on the challenges that arise due to time-periodicity
constraints.

A prerequisite to optimization is being able to accurately compute
time-periodic solutions. For low-dimensional systems (including
time-dependent PDEs with only one spatial dimension), orthogonal
collocation and (temporal) Fourier collocation algorithms such as
implemented in the software package AUTO \cite{doedel91a,doedel91b}
have proven to be robust and widely applicable
\cite{govaertz,keller:87,kevrekidis,tsai:jeng:94,okamura:10}.
In the aerodynamic optimization community, these methods are known as
harmonic balance \cite{hall2002computation}, time spectral
\cite{gopinath2005time}, or nonlinear frequency \cite{mcmullen2002application,
nadarajah2007boptimum} techniques. While these types of approaches can realize
spectral convergence in time, they quickly lead to extremely large-scale
computations as all time instances become coupled and the unknown state
vector includes all spatial degrees of freedom at every collocation point,
i.e., a tensor product between space and time.
At the other extreme, shooting methods \cite{keller:68,stoer:bulirsch}
treat only the initial conditions as unknowns and use a numerical
timestepping scheme to determine the state of the system at later
times. These methods are very effective for computing stable or
nearly stable time-periodic solutions in which nearby trajectories
do not diverge wildly from each other over the timescale of the
periodic solution and have a much smaller memory footprint than the
spectral collocation approaches. Examples of systems exhibiting this
non-chaotic behavior include mode-locked lasers
\cite{lasers}, water waves \cite{water1,water2,water3d},
viscoelastic fluid flows \cite{erica:thesis}, and rolling vehicle
tires \cite{tires1}. For chaotic systems, multiple-shooting
methods \cite{stoer:bulirsch} strike a balance between the
robustness of a temporal collocation method and the efficiency
of a shooting method for limiting the number of unknowns.

Shooting methods can further be classified by the method used to solve
for the unknown initial conditions.  If the periodic solution is
stable with fairly large decay rates relative to the period of the
solution, a simple and effective method is fixed point
iteration. However, if a high degree of accuracy is desired, the
slowest decaying modes often obstruct convergence in a reasonable
amount of simulation time. For example, Thomases and Shelley observed
``persistent oscillations'' in a Stokesian viscoelastic fluid for
which ``a simulation up to $t=10\,000$ reveals no decrease in their
amplitude,'' but could not conclude for certain that the limiting
oscillations were time-periodic.  Lust and Roose \cite{lust:roose:92}
devised a hybrid method in which some of the degrees of freedom are
solved for by a shooting method while others converge via fixed point
iteration. This works well if the number of unstable or mildly stable
modes is small. For neutrally stable problems such as computing
standing water waves, fixed point iteration cannot be used to solve
for any of the modes, and a genuinely large-scale nonlinear solver has
to be used. Mercer and Roberts developed a Newton-Raphson algorithm
for computing large-amplitude standing water waves \cite{mercer:92}.
Ambrose and Wilkening devised an adjoint-based minimization algorithm
for such problems based on the limited memory BFGS algorithm
\cite{bfgs,benj1,benj2,vtxs1}. This method was also used by Williams
et.~al.~for mode-locked lasers \cite{lasers}, and by Isaacson for computing
time-periodic solutions of the Oldroyd-B equations
\cite{erica:thesis}. Wilkening and Yu \cite{water1,water2} later
developed an overdetermined shooting method based on the
Levenberg-Marquardt method that takes advantage of consolidation of
work by computing multiple columns of the Jacobian matrix in
parallel. This approach was also used by Wilkening and Rycroft for
computing standing water waves in 3D. For studying transitions to
turbulence in Couette flow \cite{waleffe95, kawahara:kida,
  cvitanovic} or pipe flow, Viswanath developed a Newton-Krylov algorithm with a
locally constrained optimal hook step \cite{viswanath}. Similar work
was done by Schneider and collaborators in the context of turbulent
pipe flow \cite{schneider:07a, schneider:07b}.  Building on
these ideas, Govindjee, Potter and Wilkening developed a Newton-Krylov
approach for computing cyclic steady states in rolling vehicle tires
with treads \cite{tires1}.

Here we also adopt a Newton-Krylov approach, but we bring back adjoint
methods to optimize various quantities of interest over a cycle of the
(parameter-dependent) time-periodic solution.  Our motivation comes
from the goal of designing an energetically optimal flapping motion of
a 2D airfoil in a compressible, viscous fluid subject to a
time-averaged thrust constraint.  Early methods toward computing
energetically optimal flapping flight used derivative-free
optimization solvers \cite{pesavento2009flapping} and finite differences
\cite{culbreth2011high} or the sensitivity method to compute gradients for
first-order optimizers \cite{soueid2009optimization}.  The derivative-free
approach is limited to small parameter spaces and coarse flow discretizations
due to the large number of iterations required by such solvers. The
finite difference and sensitivity method also require small parameter spaces
as each entry of the gradient requires the solution of a nonlinear or
linearized, forward evolution equation, respectively.

A number of fully discrete, time-dependent, adjoint-based methods have
recently been introduced
\cite{nielsen2010discrete, van2013adjoint, jones2013adjoint, zahr2016dgopt}
that enable high-dimensional parameter spaces to be efficiently searched,
which leads to improved design and control.
Continuous adjoint methods have also been introduced for the
same purpose \cite{nadarajah2007optimum, economon2015unsteady}, although
they may result in inexact gradients and slow optimization convergence
if the spatial discretization employed is not adjoint consistent.
The method presented here improves
on these existing adjoint-based methods by incorporating time-periodicity
constraints that ensure a \emph{representative, in-flight} cycle is considered
during the optimization. The existing methods initialize the flow from the
steady-state or uniform flow and run several cycles to let non-physical
transients die out. These approaches have the disadvantage of requiring many
cycles to fully suppress the initial transients and will substantially
increase the cost of each flow simulation. Stanford and Beran
\cite{stanford2011cost} -- who considered structural optimization of a dry
flapping wing, i.e., without a surrounding fluid -- implicitly incorporated
time-periodicity constraints by using a spectral time discretization that only
supports time-periodic solutions. This approach, referred to as temporal
Fourier collocation above, has the disadvantage of coupling all timesteps,
which was mitigated in that work by considering the \emph{dry} structure and
using proper orthogonal decomposition-based model reduction methods.


The remainder of the paper is organized as follows. Section~\ref{sec:tpde}
discusses the numerical discretization of partial differential equations
and reviews Newton-Krylov and optimization-based shooting methods for computing
time-periodic solutions. A discussion on the stability of periodic orbits is
included. The fully discrete framework is emphasized throughout as this will
lead naturally to the corresponding fully discrete adjoint equations in
Section~\ref{sec:adj}. The derivation of the adjoint equations corresponding to
time-periodically constrained fully discrete partial differential equations is
provided in Section~\ref{subsec:adj-deriv}. These equations constitute
a linear, two-point boundary value problem; existence and uniqueness of
solutions is proved in \ref{sec:exist}. The corresponding adjoint
method for computing gradients of quantities of interest along the manifold of
solutions of the time-periodically constrained, fully discrete partial
differential equations is also provided in
Section~\ref{subsec:adj-deriv}. A matrix-free Krylov shooting method is
introduced in Section~\ref{subsec:adj-solver} for computing solutions of the
periodic adjoint equations. In Section~\ref{subsec:app-solvers}, the various
primal and dual shooting methods are compared side-by-side on a flapping
airfoil in compressible, viscous flow. Finally, in Section~\ref{subsec:app-opt},
the primal and dual shooting methods are used to compute optimization
functionals and gradients, respectively, to determine the energetically optimal
flapping motion in compressible, viscous flow, subject to a constraint on the
time-averaged thrust.




\section{Computing Time-Periodic Solutions of Partial Differential Equations}
\label{sec:tpde}
This section is devoted to the discretization and solution of
partial differential equations with time-periodicity constraints. This will
largely be a review of existing work on the topic
\cite{mercer:92,viswanath,benj1,vtxs1,water2,tires1}, although
emphasis will be placed on equations
that are \emph{parametrized}. This will lead to the main contribution of this
work, the fully discrete adjoint equations corresponding to time-periodic
solutions of partial differential equations and their use in computing
gradients of quantities of interest along the manifold of time-periodic
solutions.

Consider the general, nonlinear, time-periodically constrained system of
partial differential equations, parametrized by the vector
$\mubold \in \Rbb^{N_\mubold}$,
\begin{equation} \label{eqn:pde}
 \begin{aligned}
   \pder{\Ubm}{t} &= \Lboldcal(\Ubm, \mubold, t) \qquad \text{ in } \quad
   \Omega(\mubold, t) \times (0, T] \\
   \Ubm(\xbm, 0) &= \Ubm(\xbm, T),
 \end{aligned}
\end{equation}
where $\Lboldcal(\cdot, \mubold, t)$ is a spatial differential operator on the
parametrized, time-dependent domain $\Omega(\mubold, t) \subset \Rbb^{n_{sd}}$. The
boundary conditions have not been explicitly stated for brevity. This work
will only consider temporally first-order partial differential equations, or
those that have been recast as such. Without loss of generality, consider a
quantity of interest of the form
\begin{equation} \label{eqn:qoi-cont}
  \Fcal(\Ubm, \mubold) = \int_0^T \int_{\Gamma(\mubold, t)}
                         f(\Ubm, \mubold, t)\;dS\;dt,
\end{equation}
where $\Gamma(\mubold, t) \subseteq \partial \Omega(\mubold, t)$.
The generalization to other types of quantities of interest, such as volumetic
integrals and instantaneous or pointwise quantities of interest, is immediate
as the specific form of the quantity of interest will be abstracted away at the
fully discrete level. The form in (\ref{eqn:qoi-cont}) will be used in the
physical setup of the applications in Section~\ref{sec:app}.
In subsequent sections, this quantity of interest will
correspond to either the objective function or a constraint of an optimization
problem governed by a partial differential equation and subject to a
time-periodicity requirement. The remainder of this section will be concerned
with the numerical discretization and solution of (\ref{eqn:pde}).

\subsection{Numerical Discretization} \label{subsec:tpde-disc}
As the form of the spatial differential operator in the time-periodically
constrained system of partial differential equations in (\ref{eqn:pde}) was
not specified, an unspecified semi-discretization is applied to yield a system
of ordinary differential equations
\begin{equation} \label{eqn:semidisc}
  \begin{aligned}
    \Mbb\pder{\ubm}{t} &= \rbm(\ubm, \mubold, t) \\
    \ubm(0) &= \ubm(T),
  \end{aligned}
\end{equation}
where $\ubm(t) \in \Rbb^{N_\ubm}$ is the state vector of the
semi-discretization, $\rbm$ is the spatial discretization of $\Lboldcal$,
and $\Mbb$ is the mass matrix arising from the discretization of
$\displaystyle{\pder{\Ubm}{t}}$.

In Section~\ref{sec:app}, a high-order discontinuous Galerkin method is used
to semi-discretize the compressible Navier-Stokes equations.
In general, a partial differential equation with a parametrized,
time-dependent spatial domain will lead to ordinary differential equations
with a parametrized, time-dependent mass matrix. However,
Section~\ref{sec:app} will discuss an Arbitrary-Lagrangian-Eulerian formulation
of general systems of conservation laws that solves a transformed set of
equations on a \emph{fixed} domain, leading to a constant mass matrix.
Therefore, attention will be restricted to the case of a fixed mass matrix.

As the motivating applications for this work are viscous fluid dynamics
problems, an implicit temporal discretization is used since timestep
sizes tend to be stability-limited. Specifically, Diagonally Implicit
Runge-Kutta (DIRK) schemes -- Runge-Kutta schemes with a lower triangular
Butcher tableau -- are used since stable, high-order discretizations are
possible without incurring the large cost of coupling all stages.
An $s$-stage DIRK discretization of (\ref{eqn:semidisc}) leads to the
fully discrete, nonlinear evolution equations
\begin{equation} \label{eqn:dirk}
  \begin{aligned}
    \ubm^{(n)} &= \ubm^{(n-1)} + \sum_{i = 1}^s b_i\kbm^{(n)}_i \\
    \Mbb\kbm^{(n)}_i &= \Delta t_n\rbm\left(\ubm_i^{(n)},~\mubold,~
                                            t_{n-1} + c_i\Delta t_n\right),
  \end{aligned}
\end{equation}
where
\begin{equation}\label{eqn:dirk-stage}
  \ubm_i^{(n)}
               = \ubm^{(n-1)} + \sum_{j = 1}^i a_{ij}\kbm^{(n)}_j.
\end{equation}
From  (\ref{eqn:dirk}), it is clear that each timestep requires
the solution of $s$ nonlinear systems of equations of size $N_\ubm$.
Time-periodicity may then be expressed as the constraint
\begin{equation}\label{eq:time:periodicity}
  \ubm^{(0)} = \ubm^{(N_t)},
\end{equation}
where $N_t$ is the time index of the cycle period.

While the form of the fully discrete time-periodically constrained partial
differential equations are specific to a DIRK temporal discretization, this is
not a fundamental restriction of this work. Extension of the analysis and
derivations in Sections~\ref{subsec:tpde-soln}~and~\ref{sec:adj} to other
classes of temporal discretizations, whether they are implicit or explicit,
is possible.

With the full discretization of the partial differential equation addressed,
the quantity of interest, $\Fcal$, must also be discretized to be computable.
While any numerical quadrature formula can be used to perform the
discretization of the space-time integral, this may lead to a truncation error
of different orders for the governing equations and the quantity of
interest. This implies there is wasted effort since the largest order will
dominate. A solver-consistent discretization of the quantities of interest
\cite{zahr2016dgopt} -- where the spatial and temporal discretization of the
partial differential equation are also used for the quantity of interest --
can be used to circumvent this wasted effort.  In this work, the discontinuous
Galerkin shape functions will be used to discretize the spatial surface
integral and the DIRK scheme to discretize the temporal integral. Regardless
of the discretization chosen for the quantity of interest, the fully discrete
version takes the form
\begin{equation} \label{eqn:qoi-disc}
  F(\ubm^{(0)}, \dots, \ubm^{(N_t)}, \kbm_1^{(1)}, \dots, \kbm_s^{(N_t)}).
\end{equation}
Note that the aforementioned solver-consistent discretization leads to the
dependence of the fully discrete quantity of interest on the Runge-Kutta
stages, $\kbm_i^{(n)}$, combined in such a way that the expected truncation
order is attained, see \cite{zahr2016dgopt} for details. A standard numerical
quadrature, such as midpoint rule or Simpson's rule, would \emph{not} use the
stages, which are only low-order solutions of (\ref{eqn:semidisc})
\cite{zahr2016aiaa}.

With the full numerical discretization of the system of partial differential
equations complete, the next section will discuss methods for solving the
fully discrete, time-periodically constrained partial differential equations.
The periodicity constraint, i.e. $\displaystyle{\ubm^{(0)} = \ubm^{(N_t)}}$,
turns the problem into a nonlinear two-point boundary value problem, which
eliminates the possibility of using traditional evolution methods (since the
initial conditions are unknown).

\subsection{Numerical Solvers: Shooting Methods} \label{subsec:tpde-soln}
This section provides a brief, non-exhaustive review of methods which have
been introduced for solving time-periodic partial differential equations.
A distinguishing feature of this work is that we directly consider the
fully discrete form of the governing equations, whereas previous work has
focused on the continuous \cite{vandewalle1992efficient} or
semi-discrete \cite{srzednicki1994periodic} levels. The section
will conclude with a discussion of a Newton-Krylov shooting method using a
purely matrix-free Krylov solver to solve the linear systems of equations that
arise, which extends the work in \cite{tires1}.

Define $\ubm^{(N_t)}(\ubm_0; \mubold)$ as the solution of the following
initial-value problem
\begin{equation} \label{eqn:dirk2}
  \begin{aligned}
    \ubm^{(0)} &= \ubm_0 \\
    \ubm^{(n)} &= \ubm^{(n-1)} + \sum_{i = 1}^s b_i\kbm^{(n)}_i \\
    \Mbb\kbm^{(n)}_i &= \Delta t_n\rbm\left(\ubm_i^{(n)},~\mubold,~
                                            t_{n-1} + c_i\Delta t_n\right),
  \end{aligned}
\end{equation}
which can be solved using a traditional evolution algorithm that advances
the solution from timestep $n$ to $n+1$. Notice that this overloads the
notation introduced in Section~\ref{subsec:tpde-disc}, which defines
$\ubm^{(N_t)}$ as the discrete approximation of the time-periodic solution
of the system of partial differential equations at the final time. Here, it is
a nonlinear function that maps a state $\ubm_0 \in \Rbb^{N_\ubm}$ to the
state $\ubm^{(N_t)}(\ubm_0; \mubold)$. From (\ref{eqn:dirk}) and
(\ref{eqn:dirk2}), it is clear that $\ubm_0$ is the time-periodic initial
condition of the fully discrete partial differential equation, $\ubm^{(0)}$,
if it is a fixed point of $\ubm^{(N_t)}(\cdot; \mubold)$, namely
\begin{equation} \label{eqn:fixedpt}
  \ubm^{(N_t)}(\ubm_0; \mubold) = \ubm_0.
\end{equation}
Then, provided the mapping $\ubm_0 \rightarrow \ubm^{(N_t)}(\ubm_0; \mubold)$
is a contraction mapping, the Banach Fixed Point Theorem implies the existence
of the fixed point and provides a convergent algorithm for finding it, see
Algorithm~\ref{alg:fixedpt}. This is a convenient algorithm as it only relies
on solution of the nonlinear evolution equation (\ref{eqn:dirk2}), but
is known to suffer from poor convergence rates and lack of
convergence if the mapping under consideration is not a contraction.
\begin{algorithm}[!htbp]
  \caption{Fixed Point Iteration Time-Periodic Solutions of PDE}
  \label{alg:fixedpt}
  \begin{algorithmic}[1]
    \REQUIRE Initial guess for periodic initial condition, $\ubm_0$;
             parameter configuration, $\mubold$
    \ENSURE Periodic initial condition, $\ubm^{(0)}$
    \WHILE{$\norm{\ubm^{(N_t)}(\ubm_0; \mubold) - \ubm_0}_2 > \epsilon$}
      \STATE Update
      \begin{equation*}
        \ubm_0 \leftarrow \ubm^{(N_t)}(\ubm_0; \mubold)
      \end{equation*}
    \ENDWHILE
    \STATE Define periodic initial condition
    \begin{equation*}
      \ubm^{(0)} = \ubm_0
    \end{equation*}
  \end{algorithmic}
\end{algorithm}

Another class of solvers for time-periodically constrained partial differential
equations rely on unconstrained, gradient-based optimization techniques.
Define the function
\begin{equation} \label{eqn:opt-obj}
 j(\ubm_0) = \frac{1}{2}\norm{\ubm^{(N_t)}(\ubm_0; \mubold) - \ubm_0}_2^2
\end{equation}
and consider the unconstrained optimization problem
\begin{equation}
  \underset{\ubm_0 \in \Rbb^{N_\ubm}}{\text{minimize}}~~j(\ubm_0),
\end{equation}
which can be solved using gradient-based optimization techniques such
as steepest descent, the Broyden-Fletcher-Goldfarb-Shanno (BFGS)
algorithm, or its limited-memory version, L-BFGS
\cite{gill1981practical, zhu1997algorithm, nocedal2006numerical}.
The gradient of (\ref{eqn:opt-obj}), $\displaystyle{\oder{j}{\ubm_0}}$,
is usually computed using the adjoint method since the large number of
optimization variables, $N_\ubm$, renders the finite differences method or the
linearized forward method impractical \citep{gunzburger2003perspectives}.
Throughout this work, the notation
$\displaystyle{\oder{(\cdot)}{\mubold}}$ will be used to denote the total
derivative of a quantity of interest with respect to parameters --
including the explicit dependence as well as the implicit dependence through
the solution of the governing equation -- and the partial derivative notation
$\displaystyle{\pder{(\cdot)}{\mubold}}$ will be used elsewhere.
The adjoint equations for the fully discrete evolution equations
in (\ref{eqn:dirk2}) corresponding to the quantity of interest,
$j(\ubm_0)$, with parameter $\ubm_0$ are
\begin{equation} \label{eqn:uns-disc-adj-dirk}
  \begin{aligned}
    \lambdabold^{(N_t)} &= \ubm^{(N_t)}(\ubm_0; \mubold) - \ubm_0 \\
    \lambdabold^{(n-1)} &= \lambdabold^{(n)} +
      \sum_{i=1}^s \Delta t_n\pder{\rbm}{\ubm}\left(\ubm_i^{(n)},~\mubold,~
                             t_{n-1}+c_i\Delta t_n\right)^T\kappabold_i^{(n)} \\
    \Mbb^T\kappabold_i^{(n)} &= b_i\lambdabold^{(n)} +
                         \sum_{j=i}^s a_{ji}\Delta t_n\pder{\rbm}{\ubm}
                         \left(\ubm_j^{(n)},
                         ~\mubold,~t_{n-1}+c_j\Delta t_n\right)^T
                         \kappabold_j^{(n)}
  \end{aligned}
\end{equation}
for $n = 1, \dots, N_t$ and $i = 1, \dots, s$. The gradient of $j(\ubm_0)$ is
reconstructed from the dual variables as
\begin{equation} \label{eqn:funcl-grad-nosens-dirk}
  \oder{F}{\mubold} = {\lambdabold^{(0)}}^T +
                      \ubm_0 - \ubm^{(N_t)}(\ubm_0; \mubold).
\end{equation}
See \cite{zahr2016dgopt} for the derivation. These methods have been used
with considerable success to solve a variety of time-periodic partial
differential equations, including the Benjamin-Ono equation
\citep{benj1}, a wave-guide array mode-locked laser
system \cite{lasers}, and the vortex sheet with surface tension
\cite{vtxs1}.
Unfortunately, the underlying optimization algorithms suffer from relatively
slow convergence, requiring many line-searches before becoming superlinear,
and never achieve quadratic convergence.


An attractive alternative is to recast the fixed point iteration as a
nonlinear system of equations and use either the Newton-Raphson method or the
Levenberg-Marquardt method to reap the
benefits of quadratic convergence. To this end, define the nonlinear
system of equations
\begin{equation} \label{eqn:nlsys}
  \Rbm(\ubm_0) = \ubm^{(N_t)}(\ubm_0; \mubold) - \ubm_0 = 0
\end{equation}
with Jacobian matrix
\begin{equation} \label{eqn:jac-prim}
  \Jbm(\ubm_0) = \pder{\Rbm}{\ubm_0}(\ubm_0) =
                 \pder{\ubm^{(N_t)}}{\ubm_0}(\ubm_0; \mubold) - \Ibm
\end{equation}
where $\Ibm$ is the $N_\ubm \times N_\ubm$ identity matrix. The crucial
component of the Newton-Raphson method is the solution of a linear system of
equations with the Jacobian (\ref{eqn:jac-prim}), i.e.~the solution of
$\displaystyle{\Jbm(\ubm_0)\xbm = \bbm}$, given $\ubm_0 \in \Rbb^{N_\ubm}$ and
$\bbm \in \Rbb^{N_\ubm}$.  A \emph{linear} evolution equation defining
$\displaystyle{\pder{\ubm^{(N_t)}}{\ubm_0}}$, i.e. the sensitivity of the
final state with respect to perturbations in the initial state, is introduced
by linearizing the fully discrete evolution equation in (\ref{eqn:dirk2})
about the primal state $\ubm^{(n)}$, $\kbm_i^{(n)}$ with respect to the
initial state $\ubm_0$. Direct differentiation of (\ref{eqn:dirk2}) with
respect to $\ubm_0$ leads to the forward sensitivity equations
\begin{equation} \label{eqn:dirk-sens}
  \begin{aligned}
    \pder{\ubm^{(0)}}{\ubm_0} &= \Ibm \\
    \pder{\ubm^{(n)}}{\ubm_0} &= \pder{\ubm^{(n-1)}}{\ubm_0} +
                                \sum_{i = 1}^s b_i\pder{\kbm^{(n)}_i}{\ubm_0} \\
    \Mbb\pder{\kbm^{(n)}_i}{\ubm_0} &=
    \Delta t_n\pder{\rbm}{\ubm}\left(\ubm_i^{(n)},~\mubold,~
                                     t_{n-1} + c_i\Delta t_n\right)\left[
                        \pder{\ubm^{(n-1)}}{\ubm_0} +
                        \sum_{j = 1}^i a_{ij}\pder{\kbm^{(n)}_j}{\ubm_0}\right].
  \end{aligned}
\end{equation}
In general, $\displaystyle{\pder{\ubm^{(N_t)}}{\ubm_0}}$ is a large
($N_\ubm \times N_\ubm$), dense matrix that requires the solution of
$N_\ubm$ linear evolution equations to form. While it is true that the columns
of the matrix can be solved in parallel, formation and storage of this
matrix may be impractical, particularly for the large-scale computational fluid
dynamics problems that motivate this work. For non-dissipative problems
such as standing waves in the free-surface Euler equations
\citep{water2, water3d}, this is worth
the expense since all perturbation directions have to be explored (as opposed
to letting the evolution over a cycle damp out high frequency transients).
But for viscous problems such as that explored in Section~\ref{sec:app}
below, solving the Newton-Raphson equations by Krylov subspace methods
requires many fewer iterations than there are columns of the Jacobian.

Formation and storage of $\displaystyle{\pder{\ubm^{(N_t)}}{\ubm_0}}$ can be
completely avoided if a matrix-free Krylov method \cite{keyes} is used to solve the
linear systems arising in the Newton-Raphson method, i.e.
$\Jbm(\ubm_0)\xbm = \bbm$. In this case, only matrix-vector products of the
form
\begin{equation} \label{eqn:jac-prim-matvec}
  \Jbm(\ubm_0)\vbm = \pder{\Rbm}{\ubm_0}(\ubm_0) \vbm =
                 \pder{\ubm^{(N_t)}}{\ubm_0}(\ubm_0; \mubold)\vbm - \vbm
\end{equation}
for any $\vbm \in \Rbb^{N_\ubm}$, are required. For efficiency, these must be
computed without explicitly forming the matrix
$\displaystyle{\pder{\ubm^{(N_t)}}{\ubm_0}}$. This is accomplished by
considering the forward sensitivity equations in (\ref{eqn:dirk-sens}) in the
direction defined by $\vbm$. Multiplying (\ref{eqn:dirk-sens}) by the vector
$\vbm$ leads to the system of linear evolution equations
\begin{equation} \label{eqn:dirk-sens-mvp}
  \begin{aligned}
    \pder{\ubm^{(0)}}{\ubm_0}\vbm &= \vbm \\
    \pder{\ubm^{(n)}}{\ubm_0}\vbm &= \pder{\ubm^{(n-1)}}{\ubm_0}\vbm +
                            \sum_{i = 1}^s b_i\pder{\kbm^{(n)}_i}{\ubm_0}\vbm \\
    \Mbb\pder{\kbm^{(n)}_i}{\ubm_0}\vbm &=
    \Delta t_n\pder{\rbm}{\ubm}\left(\ubm_i^{(n)},~\mubold,~
                                     t_{n-1} + c_i\Delta t_n\right)\left[
                    \pder{\ubm^{(n-1)}}{\ubm_0}\vbm +
                    \sum_{j = 1}^i a_{ij}\pder{\kbm^{(n)}_j}{\ubm_0}\vbm\right].
  \end{aligned}
\end{equation}
These can be solved for $\displaystyle{\pder{\ubm^{(n)}}{\ubm_0} \cdot \vbm}$
and $\displaystyle{\pder{\kbm_i^{(n)}}{\ubm_0} \cdot \vbm}$ directly, only
requiring \emph{one} linear evolution for each $\vbm$. Since the equations in
(\ref{eqn:dirk-sens-mvp}) are linear, the underlying linear solver must be
converged to high accuracy if accurate sensitivities are to be obtained. This
mitigates the speedup with respect to the nonlinear, primal solves whose
linear systems are usually solved to low precision. For the problems
considered in Section~\ref{sec:app}, the primal equations were, on average,
$2$ times more expensive than the sensitivity equations, even though
$5$ nonlinear iterations were required for convergence.
This implies the cost of evaluating $\Rbm(\ubm_0)$ is approximately $2$
times as expensive as a Jacobian-vector product $\Jbm(\ubm_0)\vbm$. The
Newton-Krylov method, with Jacobian-vector products computed as the solution of
(\ref{eqn:dirk-sens-mvp}), is summarized in Algorithm~\ref{alg:primal}.  If
the linear system of equations arising at each iteration is solved to
sufficient accuracy, this algorithm will converge quadratically.
The starting guess can be obtained by fixed point iteration
      (Algorithm~\ref{alg:fixedpt}), or, if solutions come in families
      parametrized by an amplitude, by numerical continuation
      \cite{govaertz,keller:87,doedel91b,benj1,vtxs1,water1,water3d}.
  The latter approach is particularly useful
  if the system is not dissipative and externally driven, as fixed
  point iteration relies on transient modes being damped by the evolution
  equations, i.e.~on the periodic solution being stable and attracting.

Given this exposition on methods for computing time-periodic solutions of
partial differential equations, we turn our attention to determining the
\emph{stability} of the corresponding periodic orbit.

\begin{algorithm}[!htbp]
  \caption{Newton-Krylov Shooting Method for Time-Periodic
           Solutions of PDE} \label{alg:primal}
  \begin{algorithmic}[1]
    \REQUIRE Initial guess for periodic initial condition, $\ubm_0$;
             parameter configuration, $\mubold$
    \ENSURE Periodic initial condition, $\ubm^{(0)}$
    \WHILE{$\norm{\ubm^{(N_t)}(\ubm_0; \mubold) - \ubm_0}_2 > \epsilon$}
      \STATE Solve \emph{unsymmetric} linear system of equations
      \begin{equation*}
        \pder{\ubm^{(N_t)}}{\ubm_0}(\ubm_0;\mubold) \cdot \Delta\ubm =
                                          \ubm^{(N_t)}(\ubm_0; \mubold) - \ubm_0
      \end{equation*}
      using a \emph{matrix-free} Krylov method with matrix-vector products
      \begin{equation*}
        \pder{\ubm^{(N_t)}}{\ubm_0}(\ubm_0;\mubold) \cdot \vbm
      \end{equation*}
      computed as the solution of the directional sensitivity equations
      (\ref{eqn:dirk-sens-mvp})
      \STATE Update solution
      \begin{equation*}
        \ubm_0 \leftarrow \ubm_0 - \Delta\ubm
      \end{equation*}
    \ENDWHILE
    \STATE Define periodic initial condition
    \begin{equation*}
      \ubm^{(0)} = \ubm_0
    \end{equation*}
  \end{algorithmic}
\end{algorithm}

\subsection{Stability of Periodic Orbits of Fully Discrete Partial
            Differential Equations} \label{subsec:tpde-stab} In this
section, the concept of stability of a periodic orbit of fully
discrete partial differential equations is introduced and a method for
determining the stability of a periodic solution presented.
Fixed-point iteration will fail if the underlying solution is
  not stable and attracting. Moreover, the performance of the Newton-Krylov
  approach depends critically on whether ``most'' perturbations are
  damped out rapidly by the evolution itself.  For some problems
  without dissipation or forcing, e.g.~for standing water waves in an
  ideal fluid \cite{water1,water2,water3d}, the system is neutrally
  stable, with all Floquet multipliers lying on the unit circle. In
  that case, it pays to parallelize a full Jacobian calculation in a
  Newton-Raphson or Levenberg-Marquardt method since all perturbation
  directions have to be explored before convergence is achieved in
  practice using an iterative approach. However, for the flapping
  airfoil studied in Section~\ref{sec:app}, the Newton-Krylov approach
  converges in many fewer iterations than there are degrees of
  freedom. The purpose of this section and the results of
  Figure~\ref{fig:stab} below is to explore the stability properties
  that lead to this behavior for the flapping airfoil.

Recall the interpretation
of $\ubm^{(N_t)}$ as a function that propagates an initial condition $\ubm_0$
to its final state $\ubm^{(N_t)}(\ubm_0;~\mubold)$. Let $\ubm_0^*(\mubold)$
be the time-periodic solution of the fully discrete partial differential
equation in (\ref{eqn:dirk}) at parameter configuration $\mubold$,
i.e., $\ubm_0^*(\mubold) = \ubm^{(N_t)}(\ubm_0^*;~\mubold)$.
A periodic orbit is stable if there is a $\delta > 0$ such that
\begin{equation}\label{eqn:stab-def}
 \lim_{n \rightarrow \infty}
 \norm{\ubm^{(n \cdot N_t)}(\ubm_0^*(\mubold) + \Delta\ubm;~\mubold)
                                   - \ubm_0^*(\mubold)} = 0
\end{equation}
if $\norm{\Delta \ubm} < \delta$, where
\begin{equation}
 \ubm^{(n \cdot N_t)}(\ubm_0;~\mubold) =
         \ubm^{(N_t)}(\cdot;~\mubold) \circ \dots \circ \ubm^{(N_t)}(\ubm_0;~\mubold).
\end{equation}
A Taylor expansion of $\ubm^{(N_t)}$ about the periodic solution leads to
\begin{equation} \label{eqn:pert1}
 \ubm^{(N_t)}(\ubm_0^*(\mubold)+\Delta\ubm;~\mubold) = \ubm_0^*(\mubold) +
      \pder{\ubm^{(N_t)}}{\ubm_0}(\ubm_0^*(\mubold);~\mubold)\cdot \Delta \ubm +
      \Ocal(\norm{\Delta \ubm}^2)
\end{equation}
where time-periodicity of $\ubm_0^*(\mubold)$ was used. Repeated application of
(\ref{eqn:pert1}) leads to
\begin{equation} \label{eqn:pertn}
   \ubm^{(n \cdot N_t)}(\ubm_0^*(\mubold)+\Delta\ubm;~\mubold) =
   \ubm_0^*(\mubold) +
    \left[\pder{\ubm^{(N_t)}}{\ubm_0}(\ubm_0^*(\mubold);~\mubold)
          \right]^n\Delta\ubm + \Ocal(\norm{\Delta\ubm}^{n+1}).
\end{equation}
Taking $\delta < 1$, the stability criteria in (\ref{eqn:stab-def}) is
satisfied if all eigenvalues of
$\displaystyle{\pder{\ubm^{(N_t)}}{\ubm_0}(\ubm_0^*(\mubold);~\mubold)}$
have modulus strictly less than $1$. In Section~\ref{sec:app}, the stability
of the periodic flow around a flapping airfoil is verified using this method.

Given this exposition on solvers for time-periodically constrained partial
differential equations, we turn our attention to deriving the corresponding
fully discrete adjoint equations.

\section{Fully Discrete Time-Periodic Adjoint Method} \label{sec:adj}
In this section, the adjoint equations corresponding to the fully discrete
time-periodically constrained partial differential equations (\ref{eqn:dirk})
and quantity of interest
$\displaystyle{F(\ubm^{(0)}, \dots, \ubm^{(N_t)},
                 \kbm_1^{(1)}, \dots, \kbm_s^{(N_t)}, \mubold)}$,
will be derived. For the remainder of this section,
$\displaystyle{\ubm^{(0)}, \dots, \ubm^{(N_t)},
               \kbm_1^{(1)}, \dots, \kbm_s^{(N_t)}}$
will be taken as the time-periodic solution of the fully discrete
partial differential equations (\ref{eqn:dirk}) at parameter $\mubold$.
The adjoint equations will be derived by linearizing the fully discrete
equations about this periodic solution. This highlights the importance of an
efficient periodic solver -- the subject of Section~\ref{subsec:tpde-soln} --
as it is a prerequisite for the adjoint method.

Before proceeding to the derivation of the adjoint equations, the following
definitions are introduced for the fully discrete time-periodic constraint and
Runge-Kutta stage equations and state updates
\begin{equation} \label{eqn:abs-govern-disc}
  \begin{aligned}
    \tilde\rbm^{(0)}(\ubm^{(0)},~\ubm^{(N_t)})
    =& ~\ubm^{(0)} - \ubm^{(N_t)} = 0\\
    \tilde\rbm^{(n)}(\ubm^{(n-1)},~\ubm^{(n)},~
                        \kbm_1^{(n)}, \dots, \kbm_s^{(n)},~
                        \mubold)
    =& ~\ubm^{(n)}-\ubm^{(n-1)}-\sum_{i=1}^sb_i\kbm_i^{(i)} = 0\\
    \Rbm_i^{(n)}(\ubm^{(n-1)},
                  \kbm_1^{(n)}, \dots, \kbm_i^{(n)},
                  \mubold)
    =& ~\Mbb\kbm_i^{(n)} - \Delta t_n\rbm\left(\ubm_i^{(n)},~\mubold,~
       t_{n-1}+c_i\Delta t_n\right) = 0
  \end{aligned}
\end{equation}
for $n = 1, \dots, n$ and $i = 1, \dots, s$.

\subsection{Derivation} \label{subsec:adj-deriv}
The derivation of the fully discrete adjoint equations corresponding to the
output functional, $F$, begins with the introduction of test variables
\begin{equation}
  \lambdabold^{(0)},~\lambdabold^{(n)},~\kappabold_i^{(n)} \in \Rbb^{N_\ubm}
\end{equation}
for $n = 1, \dots, N_t$ and $i = 1, \dots, s$. Since
$\displaystyle{\ubm^{(0)}, \dots, \ubm^{(N_t)},
               \kbm_1^{(1)}, \dots, \kbm_s^{(N_t)}}$
are taken as the solution of the fully discrete time-periodic problem
in (\ref{eqn:abs-govern-disc}), the following identity holds,
\emph{for any} $\mubold \in \Rbb^{N_\mubold}$,
\begin{equation} \label{eqn:uns-disc-adj-deriv-0}
 F = F + 0 = F - {\lambdabold^{(0)}}^T\tilde\rbm^{(0)}
       - \sum_{n=1}^{N_t} {\lambdabold^{(n)}}^T\tilde\rbm^{(n)}
       - \sum_{n=1}^{N_t}\sum_{i=1}^s {\kappabold_i^{(n)}}^T\Rbm_i^{(n)}
\end{equation}
for \emph{any} value of the test functions $\lambdabold^{(n)}$ and
$\kappabold_i^{(n)}$. In (\ref{eqn:uns-disc-adj-deriv-0}), arguments have been
dropped for brevity; it is understood that all terms are evaluated at the
periodic solution of (\ref{eqn:dirk}) at parameter $\mubold$.
Since (\ref{eqn:abs-govern-disc})
holds for any $\mubold \in \Rbb^{N_\mubold}$, provided
$\displaystyle{\ubm^{(0)}, \dots, \ubm^{(N_t)},
               \kbm_1^{(1)}, \dots, \kbm_s^{(N_t)}}$
is the corresponding periodic solution,
differentiation with respect to $\mubold$ leads to
\begin{equation} \label{eqn:uns-disc-adj-derive-1}
  \begin{aligned}
    \oder{F}{\mubold} = \pder{F}{\mubold} &+
                        \sum_{n=0}^{N_t} \pder{F}{\ubm^{(n)}}
                                         \pder{\ubm^{(n)}}{\mubold} + 
                        \sum_{n=1}^{N_t} \sum_{i=1}^s\pder{F}{\kbm_i^{(n)}}
                                         \pder{\kbm_i^{(n)}}{\mubold} -
                 {\lambdabold^{(0)}}^T\left[\pder{\tilde\rbm^{(0)}}{\mubold}
      + \pder{\tilde\rbm^{(0)}}{\ubm^{(0)}}\pder{\ubm^{(0)}}{\mubold}
      + \pder{\tilde\rbm^{(0)}}{\ubm^{(N_t)}}\pder{\ubm^{(N_t)}}{\mubold}\right]
          \\ &- \sum_{n=1}^{N_t} {\lambdabold^{(n)}}^T\left[
            \pder{\tilde\rbm^{(n)}}{\mubold} +
            \pder{\tilde\rbm^{(n)}}{\ubm^{(n)}}\pder{\ubm^{(n)}}{\mubold} +
            \pder{\tilde\rbm^{(n)}}{\ubm^{(n-1)}}\pder{\ubm^{(n-1)}}{\mubold} +
            \sum_{p=1}^s \pder{\tilde\rbm^{(n)}}{\kbm_p^{(n)}}
                                             \pder{\kbm_p^{(n)}}{\mubold}\right]
          \\ &- \sum_{n=1}^{N_t} \sum_{i=1}^s {\kappabold_i^{(n)}}^T
        \left[\pder{\Rbm_i^{(n)}}{\mubold} +
             \pder{\Rbm_i^{(n)}}{\ubm^{(n-1)}}\pder{\ubm^{(n-1)}}{\mubold} +
             \sum_{j=1}^i\pder{\Rbm_i^{(n)}}{\kbm_j^{(n)}}
                         \pder{\kbm_j^{(n)}}{\mubold}\right].
   \end{aligned}
\end{equation}
Re-arrangement of terms in (\ref{eqn:uns-disc-adj-derive-1}) such that the
state variable sensitivities are isolated leads to the
following expression for $\displaystyle{\oder{F}{\mubold}}$
\begin{equation} \label{eqn:uns-disc-adj-derive-2}
  \begin{aligned}
    \oder{F}{\mubold} = \pder{F}{\mubold} &+
   \left[\pder{F}{\ubm^{(N_t)}} - {\lambdabold^{(N_t)}}^T
    \pder{\tilde\rbm^{(N_t)}}{\ubm^{(N_t)}} - {\lambdabold^{(0)}}^T
    \pder{\tilde\rbm^{(0)}}{\ubm^{(N_t)}}\right]\pder{\ubm^{(N_t)}}{\mubold} -
        \sum_{n=0}^{N_t}{\lambdabold^{(n)}}^T\pder{\tilde\rbm^{(n)}}{\mubold} - 
        \sum_{n=1}^{N_t}\sum_{p=1}^s {\kappabold_p^{(n)}}^T
                                                  \pder{\Rbm_p^{(n)}}{\mubold}\\
         &+ \sum_{n=1}^{N_t}\left[\pder{F}{\ubm^{(n-1)}} -
                 {\lambdabold^{(n-1)}}^T\pder{\tilde\rbm^{(n-1)}}{\ubm^{(n-1)}}-
                 {\lambdabold^{(n)}}^T\pder{\tilde\rbm^{(n)}}{\ubm^{(n-1)}}-
                        \sum_{i=1}^s{\kappabold_i^{(n)}}^T
                                     \pder{\Rbm_i^{(n)}}{\ubm^{(n-1)}}\right]
                        \pder{\ubm^{(n-1)}}{\mubold} \\
         &+ \sum_{n=1}^{N_t} \sum_{p=1}^s \left[\pder{F}{\kbm_p^{(n)}} -
               {\lambdabold^{(n)}}^T\pder{\tilde\rbm^{(n)}}{\kbm_p^{(n)}}-
                \sum_{i=p}^s {\kappabold_i^{(n)}}^T
                                   \pder{\Rbm_i^{(n)}}{\kbm_p^{(n)}}
                                   \right]\pder{\kbm_p^{(n)}}{\mubold}.
   \end{aligned}
\end{equation}
The dual variables, $\lambdabold^{(n)}$ and $\kappabold_i^{(n)}$, which have
remained arbitrary to this point, are chosen such that the bracketed terms in
(\ref{eqn:uns-disc-adj-derive-2}) vanish
\begin{equation} \label{eqn:uns-disc-adj}
  \begin{aligned}
    \pder{\tilde\rbm^{(0)}}{\ubm^{(N_t)}}^T\lambdabold^{(0)} +
    \pder{\tilde\rbm^{(N_t)}}{\ubm^{(N_t)}}^T\lambdabold^{(N_t)} &=
                                                       \pder{F}{\ubm^{(N_t)}} \\
     \pder{\tilde\rbm^{(n)}}{\ubm^{(n-1)}}^T\lambdabold^{(n)}
   + \pder{\tilde\rbm^{(n-1)}}{\ubm^{(n-1)}}^T\lambdabold^{(n-1)} &=
     \pder{F}{\ubm^{(n-1)}}^T - \sum_{i=1}^s \pder{\Rbm_i^{(n)}}{\ubm^{(n-1)}}^T
                                                            \kappabold_i^{(n)}\\
     \sum_{j=i}^s \pder{\Rbm_j^{(n)}}{\kbm_i^{(n)}}^T\kappabold_j^{(n)} &=
      \pder{F}{\kbm_i^{(n)}} - \pder{\tilde\rbm^{(n)}}{\kbm_i^{(n)}}^T
                                                               \lambdabold^{(n)}
  \end{aligned}
\end{equation}
for $n = 1, \dots, N_t$ and $i = 1, \dots, s$.  These are the
\emph{fully discrete adjoint equations} corresponding to the time-periodic
primal evolution equations in (\ref{eqn:abs-govern-disc}), discrete quantity
of interest $F$, and parameter $\mubold$. Defining
the dual variables as the solution of the adjoint equations in
(\ref{eqn:uns-disc-adj}), the expression for
$\displaystyle{\oder{F}{\mubold}}$ in (\ref{eqn:uns-disc-adj-derive-2})
reduces to
\begin{equation} \label{eqn:funcl-grad-nosens}
  \oder{F}{\mubold} = \pder{F}{\mubold} -
        \sum_{n=0}^{N_t}{\lambdabold^{(n)}}^T\pder{\tilde\rbm^{(n)}}{\mubold} - 
        \sum_{n=1}^{N_t}\sum_{p=1}^s {\kappabold_p^{(n)}}^T
                                                  \pder{\Rbm_p^{(n)}}{\mubold}.
\end{equation}
This provides a means of computing the total derivative
$\displaystyle{\oder{F}{\mubold}}$ without explicitly computing the large,
dense state sensitivities since the expression in (\ref{eqn:funcl-grad-nosens})
is independent of them. Direct differentiation of
$\tilde\rbm^{(n)}$ and $\Rbm_i^{(n)}$ from their definitions in
(\ref{eqn:abs-govern-disc}) leads to the final form of the adjoint equations
of the fully discrete, time-periodically constrained partial differential
equations in (\ref{eqn:dirk})
\begin{equation} \label{eqn:uns-disc-adj-dirk}
  \begin{aligned}
    \lambdabold^{(N_t)} &= \lambdabold^{(0)} + \pder{F}{\ubm^{(N_t)}}^T \\
    \lambdabold^{(n-1)} &= \lambdabold^{(n)} + \pder{F}{\ubm^{(n-1)}}^T +
      \sum_{i=1}^s \Delta t_n\pder{\rbm}{\ubm}\left(\ubm_i^{(n)},~\mubold,~
                             t_{n-1}+c_i\Delta t_n\right)^T\kappabold_i^{(n)} \\
    \Mbb^T\kappabold_i^{(n)} &= \pder{F}{\kbm_i^{(n)}}^T + b_i\lambdabold^{(n)}+
                         \sum_{j=i}^s a_{ji}\Delta t_n\pder{\rbm}{\ubm}
                         \left(\ubm_j^{(n)},
                         ~\mubold,~t_{n-1}+c_j\Delta t_n\right)^T
                         \kappabold_j^{(n)}
  \end{aligned}
\end{equation}
for $n = 1, \dots, N_t$ and $i = 1, \dots, s$.
Similarly, the total
derivative of $F$, independent of state sensitivities, takes the form
\begin{equation} \label{eqn:funcl-grad-nosens-dirk}
  \oder{F}{\mubold} = \pder{F}{\mubold} +
                     \sum_{n=1}^{N_t} \Delta t_n \sum_{i=1}^s
                      {\kappabold_i^{(n)}}^T\pder{\rbm}{\mubold}(\ubm_i^{(n)},~
                                                \mubold,~t_{n-1}+c_i\Delta t_n).
\end{equation}
From (\ref{eqn:uns-disc-adj-dirk}), it can be seen that the fully discrete
adjoint equations take the form of a \emph{linear, two-point boundary-value
problem} and cannot be solved directly as an evolution equation.
\ref{sec:exist} proves existence and uniqueness of solutions to
(\ref{eqn:uns-disc-adj-dirk}). The next section will discuss solvers for the
discrete time-periodic adjoint equations in (\ref{eqn:uns-disc-adj-dirk}).

\subsection{Numerical Solver: Matrix-Free Krylov Method}
\label{subsec:adj-solver}
As the adjoint equations corresponding to the fully discrete time-periodic
partial differential equation are linear, this section will consider
matrix-free Krylov methods to solve them. Alternatively, any of the 
methods discussed in Section~\ref{subsec:tpde-soln} could be used.

Define $\lambdabold^{(0)}(\lambdabold_{N_t}; \mubold)$ as the solution of
the linear, backward evolution equations
\begin{equation} \label{eqn:uns-disc-adj-dirk2}
  \begin{aligned}
    \lambdabold^{(N_t)} &= \lambdabold_{N_t} \\
    \lambdabold^{(n-1)} &= \lambdabold^{(n)} + \pder{F}{\ubm^{(n-1)}}^T +
      \sum_{i=1}^s \Delta t_n\pder{\rbm}{\ubm}\left(\ubm_i^{(n)},~\mubold,~
                             t_{n-1}+c_i\Delta t_n\right)^T\kappabold_i^{(n)} \\
    \Mbb^T\kappabold_i^{(n)} &= \pder{F}{\kbm_i^{(n)}} + b_i\lambdabold^{(n)} +
                         \sum_{j=i}^s a_{ji}\Delta t_n\pder{\rbm}{\ubm}
                         \left(\ubm_j^{(n)},
                         ~\mubold,~t_{n-1}+c_j\Delta t_n\right)^T
                         \kappabold_j^{(n)},
  \end{aligned}
\end{equation}
which can be directly evolved, backward-in-time. Similar to
Section~\ref{subsec:tpde-soln} this constitutes a notation overload since
$\lambdabold^{(0)} \in \Rbb^{N_\ubm}$ is the initial solution of the
adjoint equations corresponding to the fully discrete periodic
partial differential equations, as well as the linear function that takes
a state $\lambdabold_{N_t}$ to $\lambdabold^{(0)}(\lambdabold_{N_t}; \mubold)$.
Then, $\lambdabold^{(0)}(\lambdabold_{N_t}; \mubold)$ is the initial solution
of (\ref{eqn:uns-disc-adj-dirk}) if the following \emph{linear} equation is
satisfied
\begin{equation} \label{eqn:fixedpt-dual}
  \lambdabold^{(0)}(\lambdabold_{N_t}; \mubold, t) = \lambdabold_{N_t} -
                                                       \pder{F}{\ubm^{(N_t)}}^T.
\end{equation}
This is a linear system of equations of the form, $\Abm \xbm = \bbm$
where
\begin{equation} \label{eqn:jac-dual}
 \Abm = \pder{\lambdabold^{(0)}}{\lambdabold_{N_t}} - \Ibm.
\end{equation}
The columns of the linear operator $\Abm$ can be formed
by considering perturbations of (\ref{eqn:uns-disc-adj-dirk2}) with respect
to the final state $\lambdabold_{N_t}$. Differentiation of
(\ref{eqn:uns-disc-adj-dirk2}) with respect to $\lambdabold_{N_t}$ leads to
the adjoint sensitivity equations
\begin{equation} \label{eqn:uns-disc-adj-dirk-sens}
  \begin{aligned}
    \pder{\lambdabold^{(N_t)}}{\lambdabold_{N_t}} &= \Ibm \\
    \pder{\lambdabold^{(n-1)}}{\lambdabold_{N_t}} &=
    \pder{\lambdabold^{(n)}}{\lambdabold_{N_t}} +
      \sum_{i=1}^s \Delta t_n\pder{\rbm}{\ubm}\left(\ubm_i^{(n)},~\mubold,~
                             t_{n-1}+c_i\Delta t_n\right)^T
                                 \pder{\kappabold_i^{(n)}}{\lambdabold_{N_t}} \\
    \Mbb^T\pder{\kappabold_i^{(n)}}{\lambdabold_{N_t}} &=
                         b_i\pder{\lambdabold^{(n)}}{\lambdabold_{N_t}} +
                         \sum_{j=i}^s a_{ji}\Delta t_n\pder{\rbm}{\ubm}
                         \left(\ubm_j^{(n)},
                         ~\mubold,~t_{n-1}+c_j\Delta t_n\right)^T
                         \pder{\kappabold_j^{(n)}}{\lambdabold_{N_t}}.
  \end{aligned}
\end{equation}
Similar to the situation for the primal problem, the matrix
$\displaystyle{\pder{\lambdabold^{(0)}}{\lambdabold_{N_t}}}$
is an $N_\ubm \times N_\ubm$ dense matrix that requires $N_\ubm$ linear
evolution equations to form. As this is impractical for large problems,
a matrix-free Krylov method is used to solve (\ref{eqn:fixedpt-dual}),
which only requires matrix-vector products of the form
\begin{equation}
 \Abm\vbm = \pder{\lambdabold^{(0)}}{\lambdabold_{N_t}}\vbm - \vbm.
\end{equation}
The first term in this matrix-vector product can be computed directly by
considering the adjoint sensitivity equations in a given direction $\vbm$
\begin{equation} \label{eqn:uns-disc-adj-dirk-sens-mvp}
  \begin{aligned}
    \pder{\lambdabold^{(N_t)}}{\lambdabold_{N_t}}\vbm &= \vbm \\
    \pder{\lambdabold^{(n-1)}}{\lambdabold_{N_t}}\vbm &=
    \pder{\lambdabold^{(n)}}{\lambdabold_{N_t}}\vbm +
      \sum_{i=1}^s \Delta t_n\pder{\rbm}{\ubm}\left(\ubm_i^{(n)},~\mubold,~
                             t_{n-1}+c_i\Delta t_n\right)^T
                             \pder{\kappabold_i^{(n)}}{\lambdabold_{N_t}}\vbm \\
    \Mbb^T\pder{\kappabold_i^{(n)}}{\lambdabold_{N_t}}\vbm &=
                         b_i\pder{\lambdabold^{(n)}}{\lambdabold_{N_t}}\vbm +
                         \sum_{j=i}^s a_{ji}\Delta t_n\pder{\rbm}{\ubm}
                         \left(\ubm_j^{(n)},
                         ~\mubold,~t_{n-1}+c_j\Delta t_n\right)^T
                         \pder{\kappabold_j^{(n)}}{\lambdabold_{N_t}}\vbm.
  \end{aligned}
\end{equation}
The equations in (\ref{eqn:uns-disc-adj-dirk-sens-mvp}) can be solved for
$\displaystyle{\pder{\lambdabold^{(0)}}{\lambdabold_{N_t}}}\cdot \vbm$ at the
cost of one linear evolution solution for each $\vbm$. The adjoint sensitivity
equations in (\ref{eqn:uns-disc-adj-dirk-sens-mvp}) are \emph{independent} of
the quantity of interest, $F$. If there are multiple quantities of interest,
fast multiple right-hand side solvers \citep{simoncini1995iterative,
  chan1997analysis, gutknecht2006block} could be used to solve $\Abm\xbm =
\bbm$ as the matrix $\Abm$ will be fixed and only the right-hand side
varied. Furthermore, the adjoint sensitivity equations in
(\ref{eqn:uns-disc-adj-dirk-sens-mvp}) and the adjoint equations in
(\ref{eqn:uns-disc-adj-dirk2}) are identical, with the exception of the terms
$\displaystyle{\pder{F}{\ubm^{(n-1)}}}$ and
$\displaystyle{\pder{F}{\kbm_i^{(n)}}}$. Therefore, the adjoint sensitivities
are less expensive to compute than the adjoint states and the savings become
substantial if 
$\displaystyle{\pder{F}{\ubm^{(n-1)}}}$ and
$\displaystyle{\pder{F}{\kbm_i^{(n)}}}$ are expensive to compute.
Algorithm~\ref{alg:dual} below details the use of a matrix-free GMRES method
to solve (\ref{eqn:fixedpt-dual}) with matrix-vector products defined by
(\ref{eqn:uns-disc-adj-dirk-sens-mvp}).

\begin{algorithm}[htbp]
  \caption{GMRES for Solution of Fully Discrete, Time-Periodic Adjoint PDE}
  \label{alg:dual}
  \begin{algorithmic}[1]
    \REQUIRE Initial guess for periodic adjoint final condition,
             $\lambdabold_{N_t, 0}$; parameter configuration, $\mubold$;
             periodic primal solution,
             $\ubm^{(0)}, \dots, \ubm^{(N_t)},
              \kbm_1^{(1)}, \dots, \kbm_s^{(N_t)}$
    \ENSURE Periodic adjoint final condition, $\lambdabold^{(N_t)}$
    \STATE Compute
           \begin{equation*}
             \rbm_0 = \lambdabold^{(0)}(\lambdabold_{N_t, 0}, \mubold) +
                      \pder{F}{\ubm^{(N_t)}}^T - \lambdabold_{N_t, 0}
           \end{equation*}
    \STATE Set $\beta = \norm{\rbm_0}_2$, $\vbm_1 = \rbm_0/\beta$, and
               $\lambdabold_{N_t} = \lambdabold_{N_t, 0}$
    \WHILE{$\displaystyle{\norm{\lambdabold^{(0)}(\lambdabold_{N_t}, \mubold) +
                   \pder{F}{\ubm^{(N_t)}}^T - \lambdabold_{N_t}}_2 > \epsilon}$}
      \FOR{j = 1, 2, \dots, m}
        \STATE Compute
               \begin{equation*}
                 \wbm_j = \pder{\lambdabold^{(0)}}{\lambdabold_{N_t}}\vbm_j -
                          \vbm_j
               \end{equation*}
               as the solution of the adjoint sensitivity equations
               (\ref{eqn:uns-disc-adj-dirk-sens-mvp})
        \FOR{i = 1, \dots, j}
          \STATE $h_{ij} = (\wbm_j, \vbm_i)$
          \STATE $\wbm_j = \wbm_j - h_{ij}\vbm_i$
        \ENDFOR
        \STATE $h_{j+1, j} = \norm{\wbm_j}_2$
        \STATE $\vbm_{j+1} = \wbm_j/h_{j+1, j}$
      \ENDFOR
      \STATE Compute
             \begin{equation*}
               \ybm_m = \arg \min \norm{\beta\ebm_1 - \Hbm_m\ybm}_2,
             \end{equation*}
             where $\ebm_1$ is the first canonical until vector in
             $\Rbb^{N_\ubm}$ and
             $\Hbm = \{h_{ij}\}_{1 \leq i \leq m+1, 1 \leq j \leq m}$
      \STATE Update solution
             \begin{equation*}
               \lambdabold_{N_t} = \lambdabold_{N_t, 0} + \Vbm_m\ybm_m
             \end{equation*}
             where
             \begin{equation*}
               \Vbm_m = \begin{bmatrix} \vbm_1 & \cdots & \vbm_m\end{bmatrix}
             \end{equation*}
    \ENDWHILE
    \STATE Define adjoint periodic final condition
    \begin{equation*}
      \lambdabold^{(N_t)} = \lambdabold_{N_t}
    \end{equation*}
  \end{algorithmic}
\end{algorithm}

With the solution of the fully discrete primal and dual time-periodic
problems fully specified, from numerical discretization to solution
algorithms, we close this section with an algorithm that uses the
fully discrete adjoint method to compute the gradient of the
quantity of interest \emph{on the manifold of periodic solutions}.
First, the fully discrete time-periodic solution (\ref{eqn:dirk}) must
be computed, i.e. using a matrix-free Newton-Krylov method, to yield
$\displaystyle{\ubm^{(0)}, \dots, \ubm^{(N_t)},
               \kbm_1^{(1)}, \dots, \kbm_s^{(N_t)}}$.
Next, the corresponding fully discrete adjoint equations are defined about
this periodic solution and solved, i.e. using a matrix-free Krylov method,
for
$\displaystyle{\lambdabold^{(0)}, \dots, \lambdabold^{(N_t)},
               \kappabold_1^{(1)}, \dots, \kappabold_s^{(N_t)}}$.
Finally, (\ref{eqn:funcl-grad-nosens-dirk}) is used to reconstruct the desired
gradient $\displaystyle{\oder{F}{\mubold}}$.
This procedure is summarized in Algorithm~\ref{alg:grad}.

\begin{algorithm}[htbp]
  \caption{Gradients on Manifold of Time-Periodic Solutions of PDEs}
  \label{alg:grad}
  \begin{algorithmic}[1]
    \REQUIRE Parameter configuration, $\mubold$, and fully discrete quantity
             of interest,
             $F(\ubm^{(0)}, \dots, \ubm^{(N_t)},
              \kbm_1^{(1)}, \dots, \kbm_s^{(N_t)})$
    \ENSURE Gradient, $\displaystyle{\oder{F}{\mubold}}$, on manifold of
            time-periodic solutions
    \STATE For parameter $\mubold$, compute time-periodic solution of
           fully discrete PDE in (\ref{eqn:dirk}), i.e. using the Newton-Krylov
           shooting method in Algorithm~\ref{alg:primal}
           \begin{equation*}
            \ubm^{(0)}, \dots, \ubm^{(N_t)}, \kbm_1^{(1)}, \dots, \kbm_s^{(N_t)}
           \end{equation*}
    \STATE For fully discrete functional
           $F(\ubm^{(0)}, \dots, \ubm^{(N_t)},
              \kbm_1^{(1)}, \dots, \kbm_s^{(N_t)})$, compute adjoint solution
           of fully discrete time-periodic PDE in
           (\ref{eqn:uns-disc-adj-dirk}), i.e. using
           GMRES shooting method in Algorithm~\ref{alg:dual} with matrix-vector
           products computed from the backward evolution of the adjoint
           sensitivity equations in (\ref{eqn:uns-disc-adj-dirk-sens-mvp})
           \begin{equation*}
            \lambdabold^{(0)}, \dots, \lambdabold^{(N_t)},
            \kappabold_1^{(1)}, \dots, \kappabold_s^{(N_t)}
           \end{equation*}
    \STATE Reconstruct $\displaystyle{\oder{F}{\mubold}}$ using dual variables
           according to (\ref{eqn:funcl-grad-nosens-dirk})
  \end{algorithmic}
\end{algorithm}

\subsection{Generalized Reduced-Gradient Method for PDE Optimization
            with Time-Periodicity Constraints} \label{subsec:adj-grg}
Consider the fully discrete time-dependent PDE-constrained optimization problem
\begin{equation} \label{opt:disc0}
  \begin{aligned}
    & \underset{\substack{\ubm^{(0)},~\dots,~\ubm^{(N_t)} \in \Rbb^{N_\ubm},\\
                        \kbm_1^{(1)},~\dots,~\kbm_s^{(N_t)} \in \Rbb^{N_\ubm},\\
                        \mubold \in \Rbb^{N_\mubold}}}
                        {\text{minimize}}
    & & F(\ubm^{(0)},~\dots,~\ubm^{(N_t)},~\kbm_1^{(1)},~\dots,~\kbm_s^{(N_t)},~
          \mubold) \\
    & \text{subject to}
    & & \cbm(\ubm^{(0)},~\dots,~\ubm^{(N_t)},~\kbm_1^{(1)},~\dots,
            ~\kbm_s^{(N_t)},~\mubold) \geq 0 \\
    & & &   \ubm^{(0)} = \ubm_0 \\
    & & & \ubm^{(n)} = \ubm^{(n-1)} + \sum_{i = 1}^s b_i\kbm^{(n)}_i \\
    & & & \Mbb\kbm^{(n)}_i = \Delta t_n\rbm\left(\ubm_i^{(n)},~\mubold,~
                                                 t_{n-1} + c_i\Delta t_n\right)
  \end{aligned}
\end{equation}
where $F$ is a fully discrete output functional of the partial differential
equation and $\cbm$ is a vector of such output functionals. The nested or
Generalized Reduced-Gradient (GRG) approach to solve (\ref{opt:disc0})
explicitly enforces the PDE constraint at each optimization iteration. The
implicit function theorem states that the solution of the discretized PDE,
can be considered an implicit function of the parameter $\mubold$, i.e.,
$\ubm^{(n)} = \ubm^{(n)}(\mubold)$ and $\kbm_i^{(n)} = \kbm_i^{(n)}(\mubold)$.
Strict enforcement of the discretized partial differential equation allows the
PDE variables and equations to be removed from the optimization problem
\begin{equation} \label{opt:disc1}
  \begin{aligned}
    & \underset{\mubold \in \Rbb^{N_\mubold}}{\text{minimize}}
    & & F(\ubm^{(0)}(\mubold),~\dots,~\ubm^{(N_t)}(\mubold),
         ~\kbm_1^{(1)}(\mubold),~\dots,~\kbm_s^{(N_t)}(\mubold),~
            \mubold) \\
    & \text{subject to}
    & & \cbm(\ubm^{(0)}(\mubold),~\dots,~\ubm^{(N_t)}(\mubold),
            ~\kbm_1^{(1)}(\mubold),~\dots,~\kbm_s^{(N_t)}(\mubold),~\mubold) \geq 0.
  \end{aligned}
\end{equation}
To solve this optimization problem using gradient-based techniques, the terms
$\displaystyle{\oder{F}{\mubold}}$ and $\displaystyle{\oder{\cbm}{\mubold}}$ --
gradients of quantities of interest along the manifold of solutions of the
PDE -- are required. Depending on the relative number of variables in
$\mubold$ to the number of constraints in $\cbm$, the direct or adjoint method
can be efficiently used to compute these gradients \emph{without relying on
finite differences}.

Now consider the optimization problem in (\ref{opt:disc0}) with the
time-periodicity constraint added
\begin{equation} \label{opt:disc2}
  \begin{aligned}
    & \underset{\substack{\ubm^{(0)},~\dots,~\ubm^{(N_t)} \in \Rbb^{N_\ubm},\\
                        \kbm_1^{(1)},~\dots,~\kbm_s^{(N_t)} \in \Rbb^{N_\ubm},\\
                        \mubold \in \Rbb^{N_\mubold}}}
                        {\text{minimize}}
    & & F(\ubm^{(0)},~\dots,~\ubm^{(N_t)},~\kbm_1^{(1)},~\dots,~\kbm_s^{(N_t)},~
          \mubold) \\
    & \text{subject to}
    & & \cbm(\ubm^{(0)},~\dots,~\ubm^{(N_t)},~\kbm_1^{(1)},~\dots,
            ~\kbm_s^{(N_t)},~\mubold) \geq 0 \\
    & & &   \ubm^{(0)} = \ubm^{(N_t)} \\
    & & & \ubm^{(n)} = \ubm^{(n-1)} + \sum_{i = 1}^s b_i\kbm^{(n)}_i \\
    & & & \Mbb\kbm^{(n)}_i = \Delta t_n\rbm\left(\ubm_i^{(n)},~\mubold,~
                                                 t_{n-1} + c_i\Delta t_n\right).
  \end{aligned}
\end{equation}
Strict enforcement of the time-periodic partial differential equations leads to
an application of the implicit function theorem, similar to that above, i.e.,
$\ubm^{(n)} = \ubm^{(n)}(\mubold)$ and $\kbm_i^{(n)} = \kbm_i^{(n)}(\mubold)$,
where $\ubm^{(n)}$ and $\kbm_i^{(n)}$ are the time-periodic solution of the
discrete partial differential equations. This results in an optimization
problem identical to that in (\ref{opt:disc1}) with this new definition of
$\ubm^{(n)}(\mubold)$ and $\kbm_i^{(n)}(\mubold)$. The novel periodic adjoint
method, derived in Section~\ref{subsec:adj-deriv}, can be used to compute
gradients along the manifold of time-periodic solutions of the fully discrete
PDE, i.e.
$\displaystyle{\oder{F}{\mubold}}$ and $\displaystyle{\oder{\cbm}{\mubold}}$,
for the use in gradient-based optimization.

\section{Applications} \label{sec:app}
This section will present two numerical experiments that study time-periodic
solutions of the compressible Navier-Stokes equations, and the behavior of the
methods discussion. The first experiment will solely consider solutions of
the fully discrete primal and dual periodic two-point boundary-value problems,
and study convergence of the methods introduced. The second will apply the
novel, fully discrete, periodic adjoint method to solve an optimal control
problem governed by a time-periodically constrained partial differential
equation.

The compressible Navier-Stokes equations take the form
\begin{align}
\frac{\partial \rho}{\partial t}  + \frac{\partial}{\partial x_i}
(\rho u_i) &= 0, \label{ns1}\\
\frac{\partial}{\partial t} (\rho u_i) +
\frac{\partial}{\partial x_i} (\rho u_i u_j+ p)  &=
\phantom{+}\frac{\partial \tau_{ij}}{\partial x_j}
\quad\text{for }i=1,2,3, \label{ns2} \\
\frac{\partial}{\partial t} (\rho E) +
\frac{\partial}{\partial x_i} \left(u_j(\rho E+p)\right) &=
-\frac{\partial q_j}{\partial x_j}
+\frac{\partial}{\partial x_j}(u_j\tau_{ij}), \label{ns3}
\end{align}
in $\Omega(\mubold, t)$ where $\rho$ is the fluid density, $u_1,u_2,u_3$ are the
velocity components, and $E$ is the total energy. The viscous stress tensor and
heat flux are given by
\begin{align}
\tau_{ij} = \mu
\left( \frac{\partial u_i}{\partial x_j} +
\frac{\partial u_j}{\partial x_i} -\frac23
\frac{\partial u_k}{\partial x_k} \delta_{ij} \right)
\qquad \text{ and } \qquad
q_j = -\frac{\mu}{\mathrm{Pr}} \frac{\partial}{\partial x_j}
\left( E+\frac{p}{\rho} -\frac12 u_k u_k \right).
\end{align}
Here, $\mu$ is the viscosity coefficient and $\mathrm{Pr = 0.72}$ is
the Prandtl number which we assume to be constant. For an ideal gas,
the pressure $p$ has the form
\begin{align}
p=(\gamma-1)\rho \left( E - \frac12 u_k u_k\right), \label{ns5}
\end{align}
where $\gamma = 1.4$ is the adiabatic gas constant.  Furthermore, entropy of the
system is assumed constant, i.e.~isentropic, which is equivalent to the flow
being adiabatic and reversible.  For a perfect gas, the entropy is defined as
\begin{equation}\label{eqn:entropy}
  s = p/\rho^\gamma.
\end{equation}
From the isentropic assumption, (\ref{eqn:entropy}) explicitly relates the
pressure and density of the flow, rendering the energy equation redundant.
This effectively reduces the square system of PDEs of size $n_{sd}+2$ to one
of size $n_{sd}+1$, where $n_{sd}$ is the number of spatial dimensions.
It can be shown, under suitable assumptions, that the solution of the
isentropic approximation of the Navier-Stokes equations converges to
the solution of the incompressible Navier-Stokes equations as the Mach
number approaches zero \cite{lin1995incompressible, desjardins1999incompressible,
froehle2013high}.

\begin{figure}[h]
  \begin{subfigure}{0.49\textwidth}
    \centering
    \includegraphics[width=2.5in]{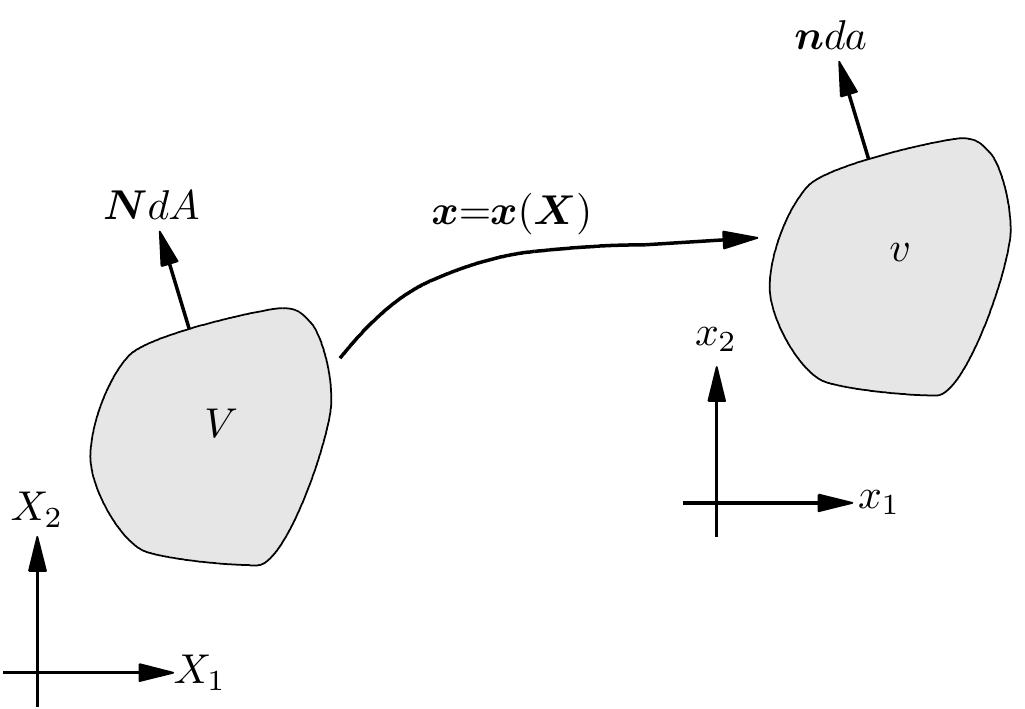}
  \end{subfigure} \hfill
  \begin{subfigure}{0.49\textwidth}
    \centering
    \includegraphics[width=2.25in]{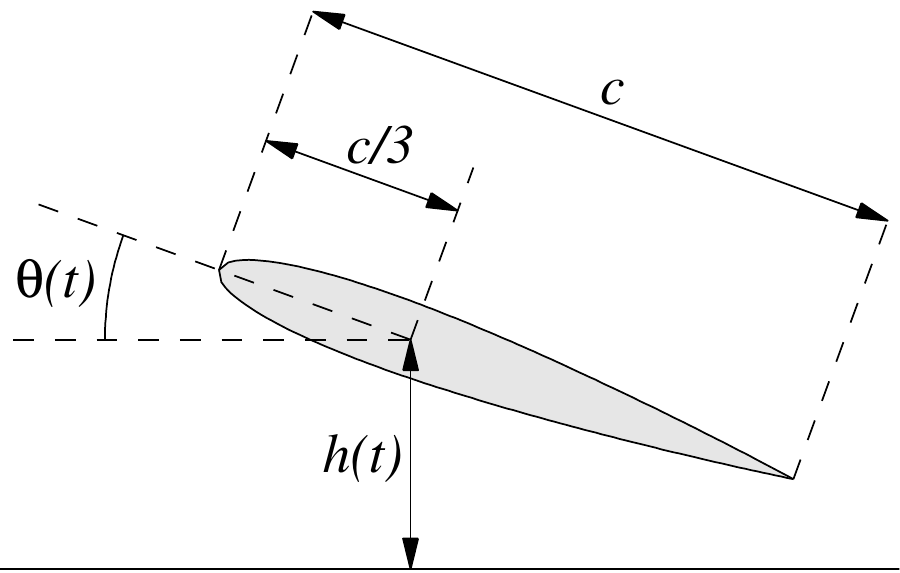}
  \end{subfigure}
  \caption{Time-dependent mapping between reference and physical domain (left)
           and kinematic description of body under consideration, NACA0012
           airfoil (right).} \label{fig:mapfoil}
\end{figure}

The compressible Navier-Stokes equations can be written in conservation
form as
\begin{equation}\label{eqn:claw-phys}
  \pder{\Ubm}{t} + \nabla\cdot\Fbm(\Ubm,~\nabla\Ubm) = 0
                                                 \quad \text{in}~~\Omega(\mubold, t),
\end{equation}
where $\Omega(\mubold, t)$ is, in general, a time-dependent, parametrized domain.
The Arbitrary-Lagrangian-Eulerian (ALE) form of the conservation law are
obtained -- following the work in \citep{persson2009dgdeform} -- by
introducing a mapping $\Gcal$ from the deforming domain $\Omega(\mubold, t)$ to a
fixed, reference domain $V$ (Figure~\ref{fig:mapfoil}). This mapping is used
to transform the conservation law on $\Omega(\mubold, t)$ to one on $V$.
The transformed system of conservation laws, defined on the reference domain
$V$ takes the form
\begin{equation} \label{eqn:ns-cons-ref}
  \left.\pder{\Ubm_\Xbm}{t}\right|_\Xbm +
  \nabla_\Xbm\cdot\Fbm_\Xbm(\Ubm_\Xbm,~\nabla_\Xbm\Ubm_\Xbm) = 0
\end{equation}
where $\nabla_\Xbm$ denotes spatial derivatives with respect to the reference
variables, $\Xbm$.  The transformed state vector, $\Ubm_\Xbm$, and its
corresponding spatial gradient with respect to the reference configuration
take the form
\begin{equation} \label{eqn:transform1-nogcl}
\Ubm_\Xbm = g\Ubm, \qquad
\nabla_\Xbm \Ubm_\Xbm = g^{-1}\Ubm_\Xbm\pder{g}{\Xbm} +
                        g\nabla\Ubm \cdot \Gbm,
\end{equation}
where $\Gbm = \nabla_\Xbm \Gcal$, $g = \text{det}(\Gbm)$,
$\displaystyle{\vbm_\Gbm = \pder{\xbm}{t} = \pder{\Gcal}{t}}$,
and the arguments have been dropped, for brevity. The transformed fluxes are
\begin{equation} \label{eqn:transform2-nogcl}
  \begin{aligned}
    \Fbm_\Xbm(\Ubm_\Xbm, \nabla_\Xbm\Ubm_\Xbm) &= \Fbm_\Xbm^{inv}(\Ubm_\Xbm) +
                     \Fbm_\Xbm^{vis}(\Ubm_\Xbm, \nabla_\Xbm\Ubm_\Xbm), \qquad \\
    \Fbm_\Xbm^{inv}(\Ubm_\Xbm) &= g\Fbm^{inv}(g^{-1}\Ubm_\Xbm)\Gbm^{-T} -
                                 \Ubm_\Xbm\otimes\Gbm^{-1}\vbm_\Gbm, \qquad \\
    \Fbm_\Xbm^{vis}(\Ubm_\Xbm, \nabla_\Xbm\Ubm_\Xbm) &=
     g\Fbm^{vis}\left(g^{-1}\Ubm_\Xbm, g^{-1}\left[\nabla_\Xbm\Ubm_\Xbm
                - g^{-1}\Ubm_\Xbm\pder{g}{\Xbm}\right]\Gbm^{-1}\right)\Gbm^{-T}.
  \end{aligned}
\end{equation}
For details regarding the derivation of the transformed equations and
modifications that ensure the Geometric Conservation Law
\cite{thomas1979geometric} is satisfied, the reader is referred to
\cite{persson2009dgdeform, zahr2016dgopt}. The system of equations
\eqref{eqn:ns-cons-ref} is discretized using a nodal discontinuous Galerkin
(DG) method on unstructured meshes of triangles, with polynomial degrees 3
within each element. The viscous fluxes are chosen according to the compact DG
method \cite{peraire2008compact} method, and the our implementation is fully
implicit with exact Jacobian matrices and a range of parallel iterative
solvers \cite{persson08newtongmres}.  The resulting semi-discrete system has
the form of our general system of ODEs \eqref{eqn:semidisc}.
All partial derivatives of the semi-discrete governing equations
and corresponding quantities of interest, namely
$\displaystyle{\pder{\rbm}{\ubm},~\pder{\rbm}{\mubold},~\pder{f_h}{\ubm},~
               \pder{f_h}{\mubold}}$
are computed via automatic symbolic differentiation at the element-level with
the MAPLE software \cite{maple1994waterloo} and subsequent assembly. The
semi-discrete quantity of interest $f_h$ is defined as the approximation of
$\displaystyle{\int_{\Gammabold(\mubold,\,t)} f(\Ubm,\,\mubold,\,t)\,dS}$ in
(\ref{eqn:qoi-cont}) using the DG shape functions and required, along with the
temporal discretization scheme, to compute the discrete output functional $F$
in (\ref{eqn:qoi-disc}). Additional details regarding computation of the
partial derivatives with respect to $\mubold$ in the case of a parametrized,
deforming domain are provided in \cite{zahr2016dgopt}.

The remainder of this document will consider the time-periodic solution and
optimization of a flapping NACA0012 airfoil, shown in Figure~\ref{fig:mapfoil}.
Two quantities of interest that will be considered are the total work exerted
by the fluid on the airfoil, $\Wcal$, and the impulse in the $x$-direction
imparted on the airfoil by the fluid, $\Jcal_x$, which take the form
\begin{equation} \label{eqn:qoi}
\Wcal(\Ubm, \mubold) =
            \int_0^T\int_\Gammabold \fbm(\Ubm, \mubold, t) \cdot \dot\xbm~dS~dt
\qquad \text{ and } \qquad
\Jcal_x(\Ubm, \mubold)=
            \int_0^T\int_\Gammabold \fbm(\Ubm, \mubold, t) \cdot \ebold_1~dS~dt
\end{equation}
In this case, $\Gammabold$ is the surface of the airfoil,
$\ebm_1 \in \Rbb^{n_{sd}}$ is the $1$st canonical unit vector,
$\fbm(\Ubm, \mubold, t) \in \Rbb^{n_{sd}}$ is the instantaneous force that the
fluid exerts on the airfoil, and $\dot\xbm$ is the pointwise velocity of
airfoil. The solver-consistent discretization, discussed in
Section~\ref{subsec:tpde-disc} and \cite{zahr2016dgopt}, of these quantities
results in the fully discrete approximations
$W(\ubm^{(0)}, \dots, \ubm^{(N_t)}, \kbm_1^{(1)}, \dots, \kbm_s^{(N_t)},
   \mubold)$
and
$J_x(\ubm^{(0)}, \dots, \ubm^{(N_t)}, \kbm_1^{(1)}, \dots, \kbm_s^{(N_t)},
     \mubold)$.
The \emph{instantaneous} quantities of interest corresponding to those in
(\ref{eqn:qoi}) are the power and $x$-directed force the fluid exerts on the
airfoil, which take the form
\begin{equation*}
\Pcal(\Ubm, \mubold, t) =
                       \int_\Gammabold \fbm(\Ubm, \mubold, t) \cdot \dot\xbm~dS
\qquad \text{ and } \qquad
\Fcal_x(\Ubm, \mubold, t)=
                       \int_\Gammabold \fbm(\Ubm, \mubold, t) \cdot \ebold_1~dS.
\end{equation*}
Define $\Pcal^h(\ubm, \mubold, t)$ and $\Fcal_x^h(\ubm, \mubold, t)$ as the
solver-consistent semi-discretization of these instantaneous quantities of
interest.

\subsection{Time-Periodic Solutions of the Compressible Navier-Stokes Equations}
\label{subsec:app-solvers}
In this section, the various solvers discussed in this document for
determining primal and dual time-periodic solutions of partial differential
equations are compared for a flapping airfoil in an isentropic, viscous
flow. The stability of the periodic orbit is verified by performing an
eigenvalue analysis of $\displaystyle{\pder{\ubm^{(N_t)}}{\ubm_0}}$. The
section closes with validation of the adjoint method, introduced for
efficient gradient computation of quantities of interest, against a
second-order finite difference approximation.

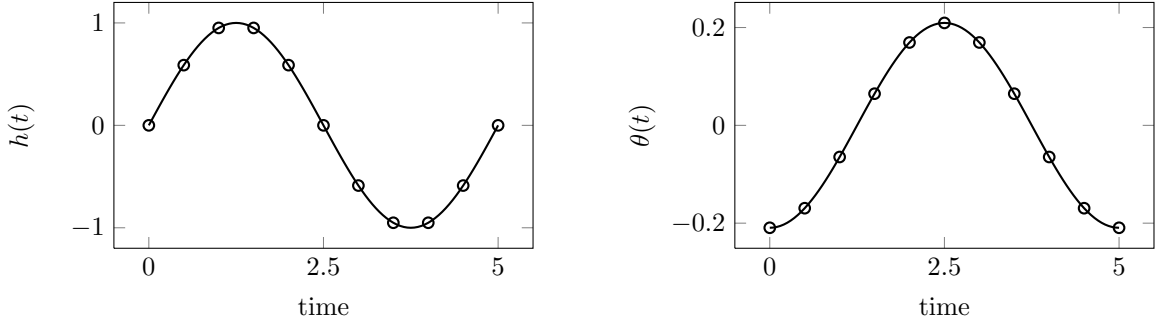
\begin{figure}
\hfill
\begin{subfigure}{0.45\textwidth}
 \begin{tikzpicture}

\begin{axis}[
  scale only axis,
  width=0.75\textwidth,
  height=0.15\textheight,
  xtick={0, 2.5, 5},
  xlabel=time,
  ylabel=$h(t)$]
\addplot [black, solid, thick, mark=o, mark repeat=10]  table[x expr=0.05*\coordindex, y index=2] {foil2d/cfd/flap00/periodic/prim/app0_mesh-how1-0-3_cfd-1p4-0p2-1000-0p72-0-1-_1-0__tdisc2-0p2-100-_1-1_-3.xyth};
\end{axis}

\end{tikzpicture}
\end{subfigure} \hfill
\begin{subfigure}{0.45\textwidth}
 \begin{tikzpicture}

\begin{axis}[
  scale only axis,
  width=0.75\textwidth,
  height=0.15\textheight,
  xtick={0, 2.5, 5},
  xlabel=time,
  ylabel=$\theta(t)$]
\addplot [black, solid, thick, mark=o, mark repeat=10]  table[x expr=0.05*\coordindex, y index=3] {foil2d/cfd/flap00/periodic/prim/app0_mesh-how1-0-3_cfd-1p4-0p2-1000-0p72-0-1-_1-0__tdisc2-0p2-100-_1-1_-3.xyth};
\end{axis}

\end{tikzpicture}
\end{subfigure}
\hfill
\caption{Trajectories of $h(t)$ and $\theta(t)$ that define the motion of the
         airfoil in Figure~\ref{fig:mapfoil} and will be used to study primal
         and dual time-periodic solvers.}
\label{fig:traj}
\end{figure}

\begin{figure}
 \centering
 \begin{subfigure}{0.3\textwidth}
  \centering
  \includegraphics[width=\textwidth]{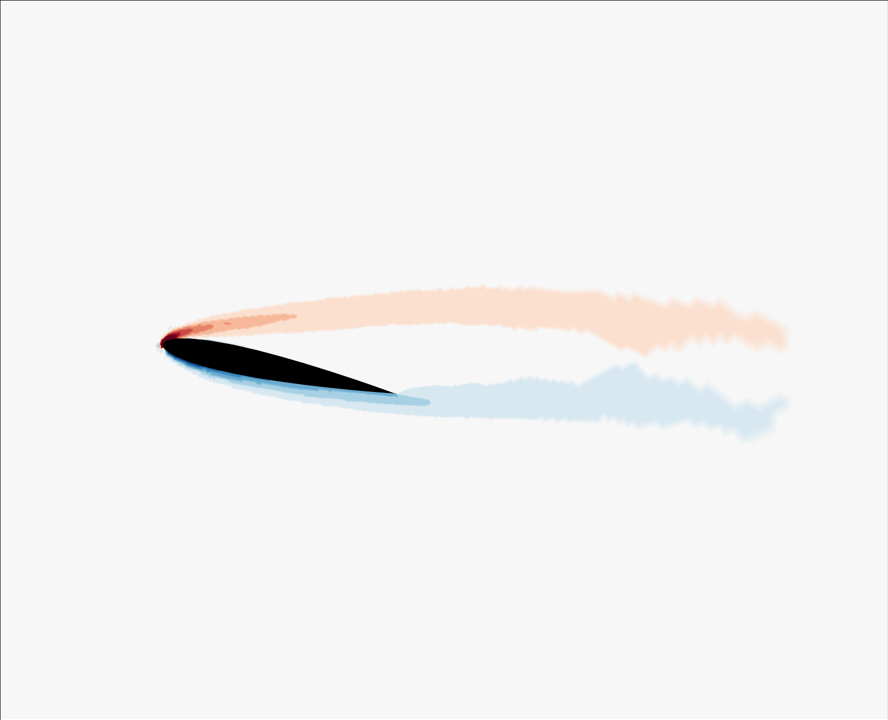}
 \end{subfigure} \hfill
 \begin{subfigure}{0.3\textwidth}
  \centering
  \includegraphics[width=\textwidth]{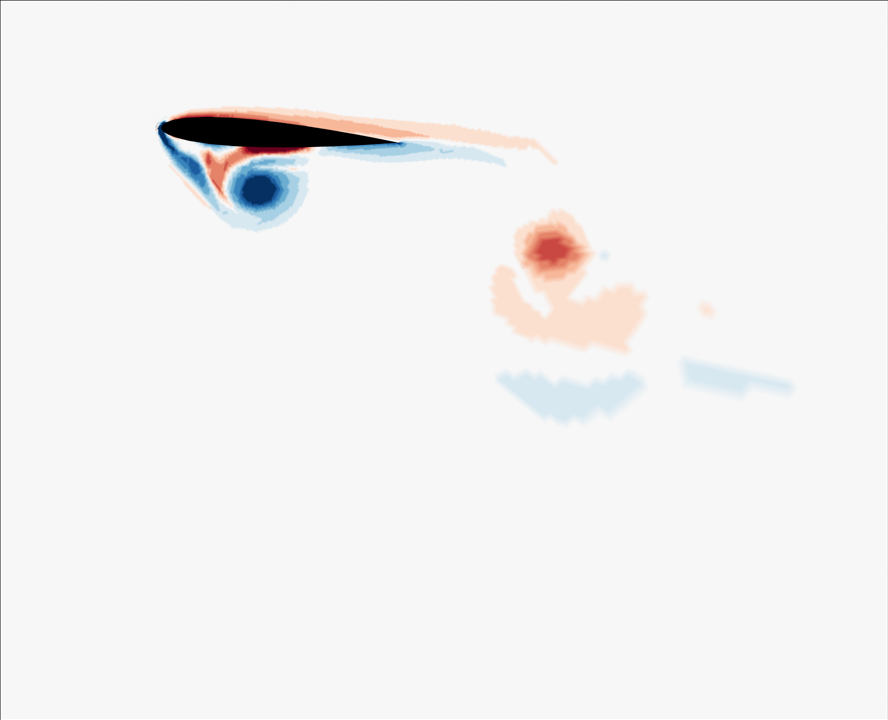}
 \end{subfigure} \hfill
 \begin{subfigure}{0.3\textwidth}
  \centering
  \includegraphics[width=\textwidth]{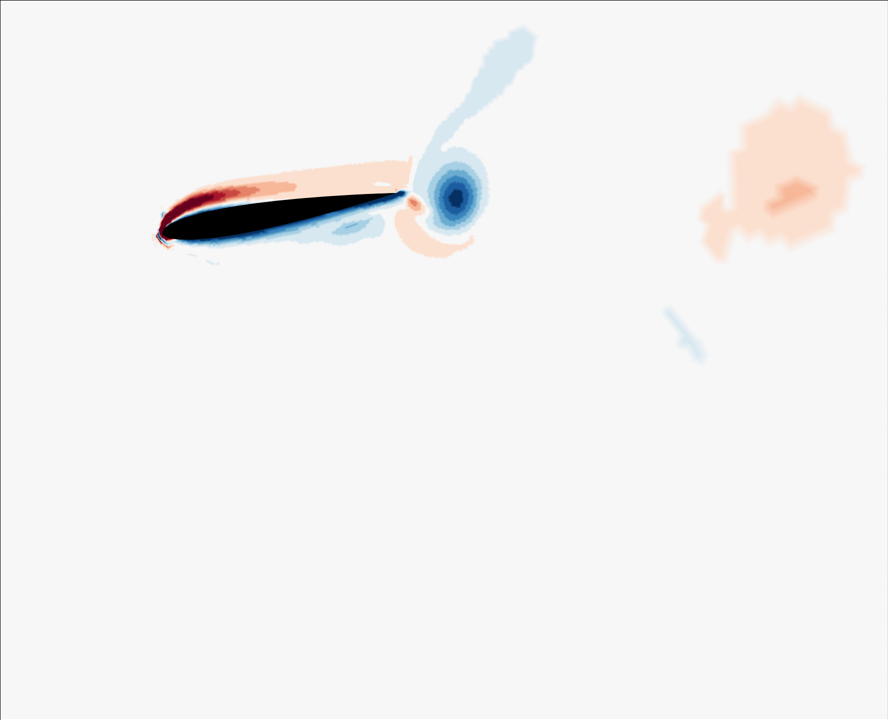}
 \end{subfigure} \\
 \begin{subfigure}{0.3\textwidth}
  \centering
  \includegraphics[width=\textwidth]{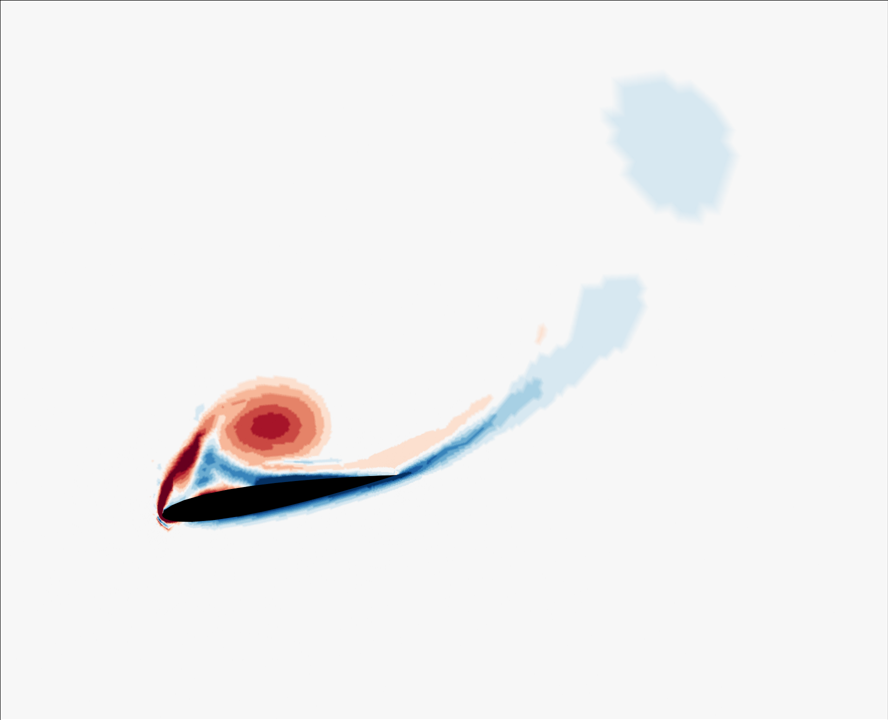}
 \end{subfigure} \hfill
 \begin{subfigure}{0.3\textwidth}
  \centering
  \includegraphics[width=\textwidth]{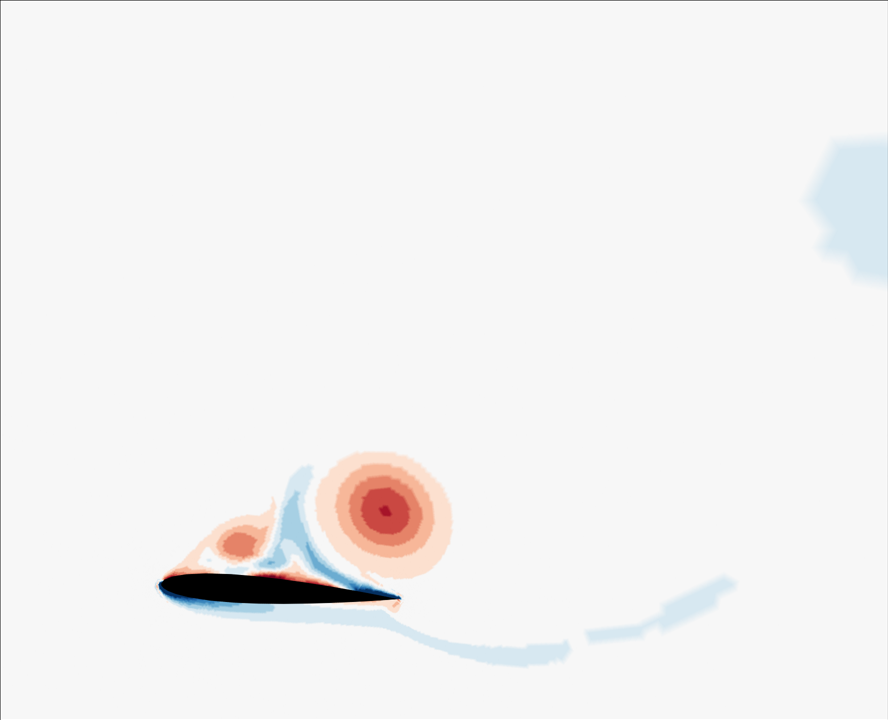}
 \end{subfigure} \hfill
 \begin{subfigure}{0.3\textwidth}
  \centering
  \includegraphics[width=\textwidth]{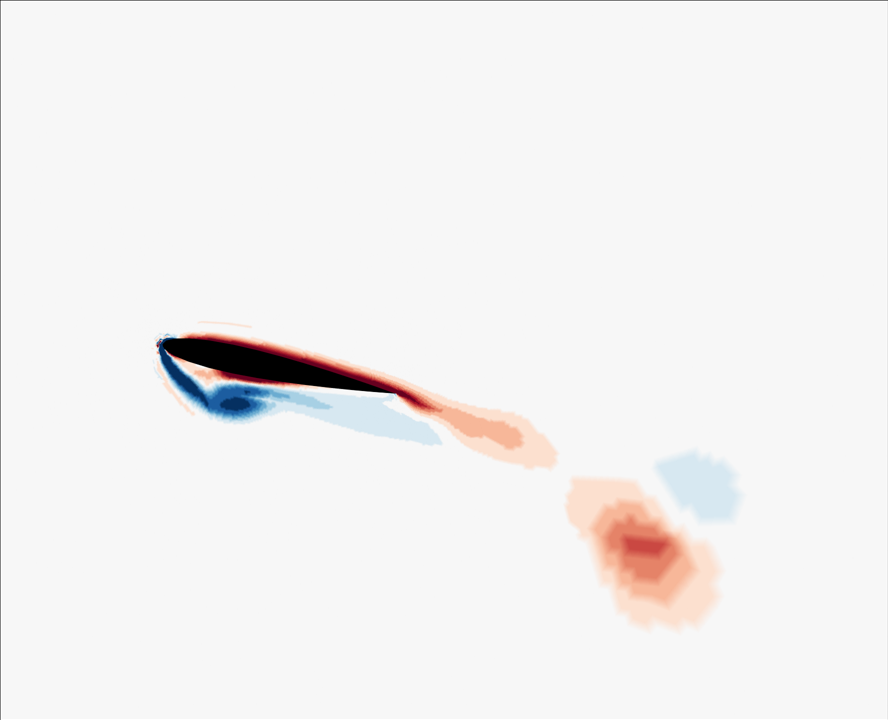}
 \end{subfigure}
 \caption{Flow vorticity around heaving/pitching airfoil for simulation
          initialized from steady state flow. Non-physical transients are
          introduced at the beginning of the time interval that result in
          non-trivial errors in integrated quantities of interests.
          Snapshots taken at times $t = 0.0,~1.0,~2.0,~3.0,~4.0,~5.0$.}
 \label{fig:from-ss}
\end{figure}

\begin{figure}
 \centering
 \begin{subfigure}{0.3\textwidth}
  \centering
  \includegraphics[width=\textwidth]{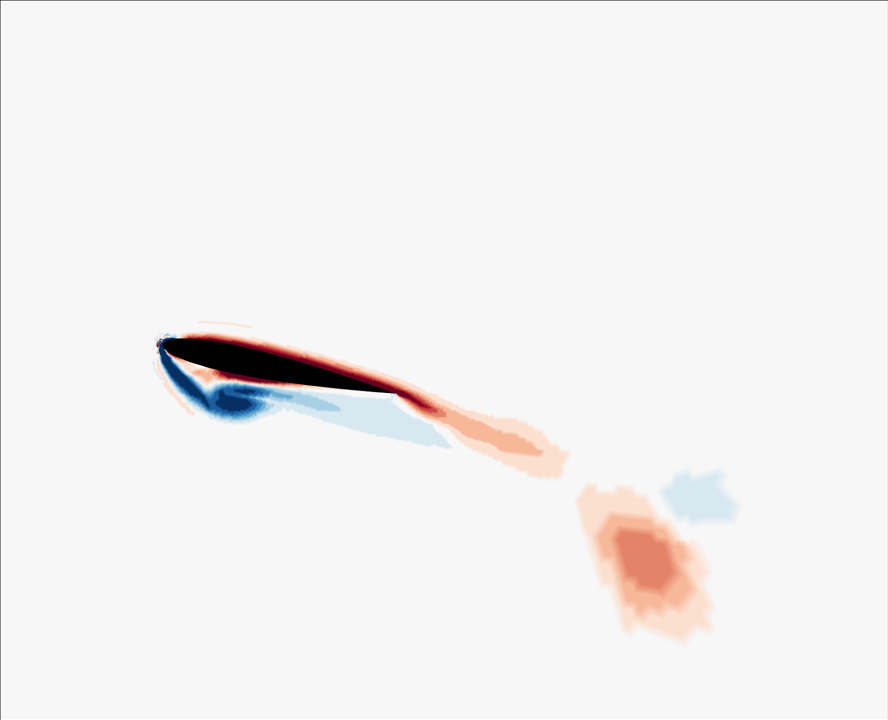}
 \end{subfigure} \hfill
 \begin{subfigure}{0.3\textwidth}
  \centering
  \includegraphics[width=\textwidth]{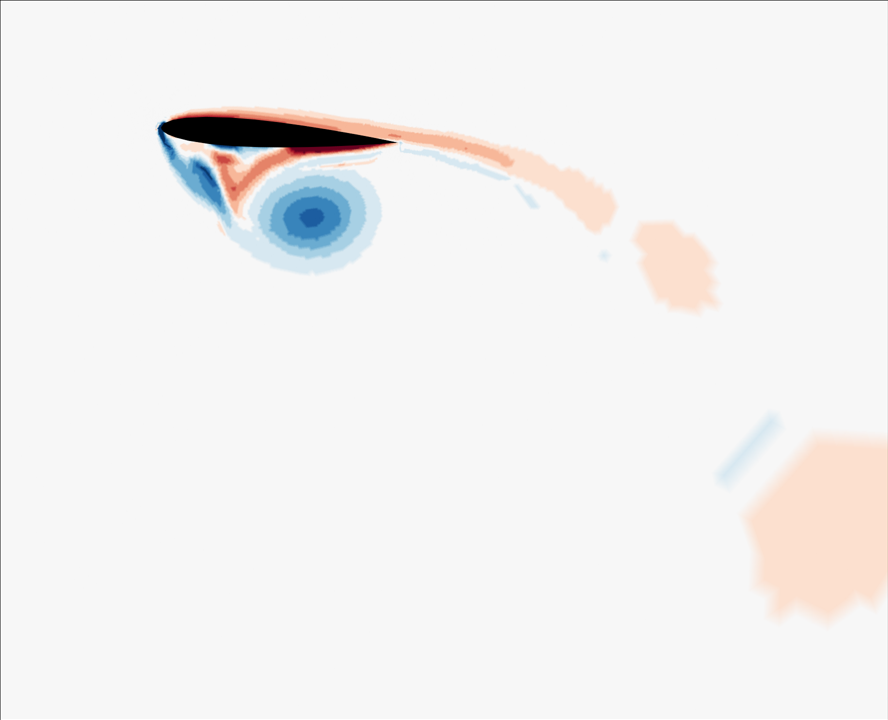}
 \end{subfigure} \hfill
 \begin{subfigure}{0.3\textwidth}
  \centering
  \includegraphics[width=\textwidth]{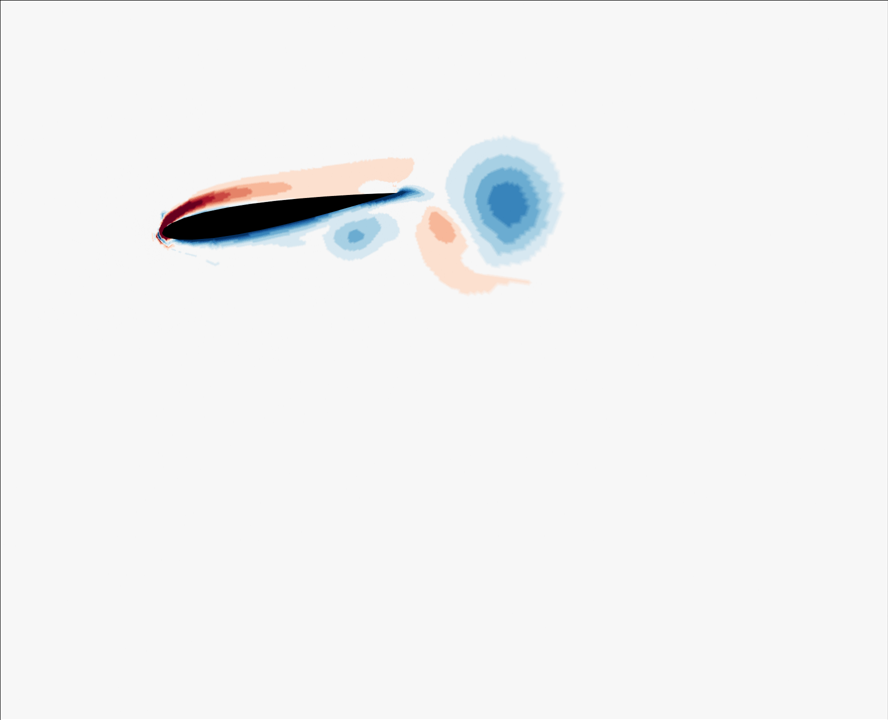}
 \end{subfigure} \\
 \begin{subfigure}{0.3\textwidth}
  \centering
  \includegraphics[width=\textwidth]{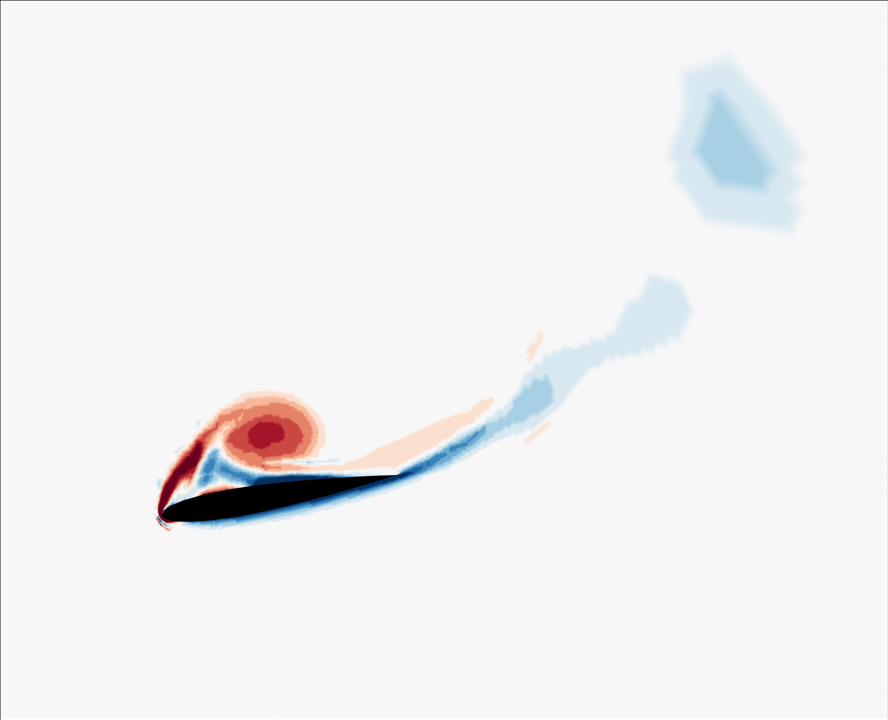}
 \end{subfigure} \hfill
 \begin{subfigure}{0.3\textwidth}
  \centering
  \includegraphics[width=\textwidth]{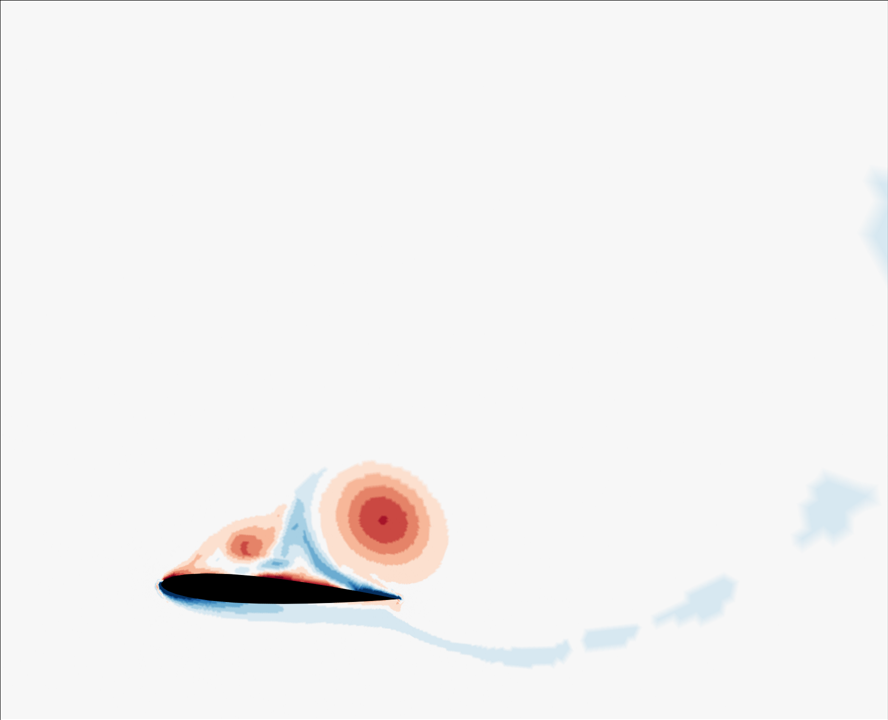}
 \end{subfigure} \hfill
 \begin{subfigure}{0.3\textwidth}
  \centering
  \includegraphics[width=\textwidth]{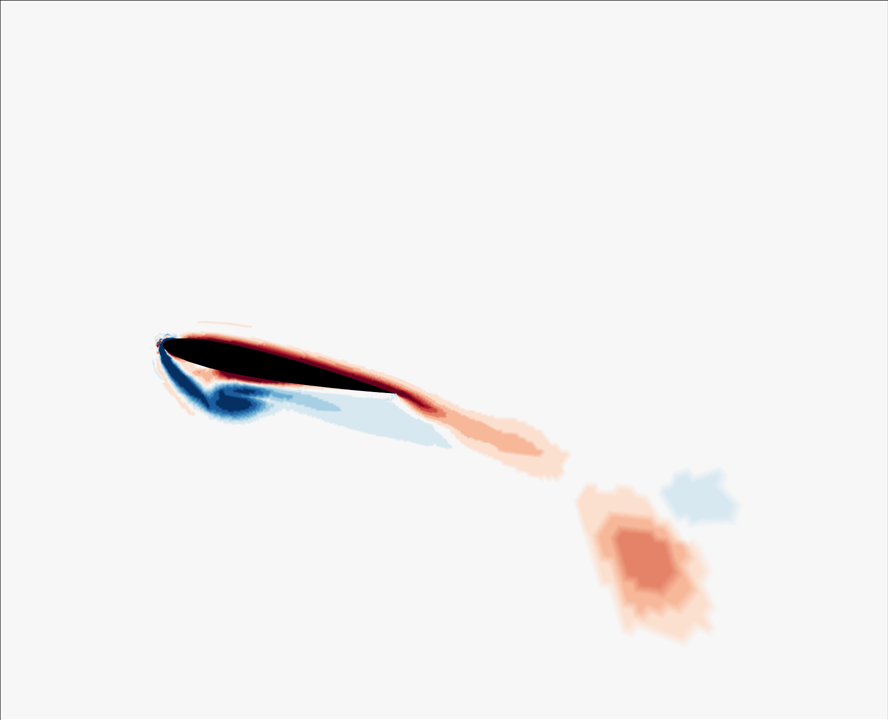}
 \end{subfigure}
 \caption{Time-periodic flow vorticity around heaving/pitching airfoil,
          i.e., initialized from periodic initial condition. The time-periodic
          initial condition ensures transients are not introduced at the
          beginning of the simulation; the result is a seamless transition
          between periods, as would be experienced in-flight, and trusted
          integrated quantities of interest. Snapshots taken at times
          $t = 0.0,~1.0,~2.0,~3.0,~4.0,~5.0$.}
 \label{fig:from-periodic}
\end{figure}

Consider the NACA0012 airfoil in Figure~\ref{fig:mapfoil} immersed in an
isentropic, viscous flow with Reynolds number set to $1000$. The Mach number
is set to $0.2$, in order to get a good compromise between approximate
incompressibility and well-conditioned equations. The kinematic motion of the
foil is parametrized with a single Fourier mode, i.e.,
\begin{equation}
 \begin{aligned}
  h(\mubold,\,t) &= A_h\sin(\omega_h t + \phi_h) + c_h \\
  \theta(\mubold,\,t) &= A_\theta\sin(\omega_\theta t + \phi_\theta) + c_\theta.
 \end{aligned}
\end{equation}
The vector of parameters is fixed for the remainder of this section
\begin{equation} \label{eqn:nominal}
\mubold = \begin{bmatrix} A_h & \omega_h & \phi_h & c_h &
                  A_\theta & \omega_\theta & \phi_\theta & c_\theta\end{bmatrix}
        = \begin{bmatrix} 1.0 & 0.4\pi & 0.0 & 0.0 &
                  \frac{\pi}{15} & 0.4\pi & \frac{\pi}{2} & 0.0\end{bmatrix},
\end{equation}
and corresponds to the motion in Figure~\ref{fig:traj} with period $T = 5$.
The mapping $\Gcal(\Xbm,\,t)$ from the fixed reference domain $V$ to
the physical domain $\Omega(\mubold,\,t)$ takes the form of a parametrized
rigid body motion
\begin{equation} \label{eqn:defmap}
 \Gcal(\Xbm,\,t) = \vbm(\mubold,\,t)+\Qbm(\mubold,\,t)(\Xbm - \xbm_0)+\xbm_0,
\end{equation}
where $\xbm_0$ is the location of pitching axis in the reference configuration
(the $1/3$ chord) and
\begin{equation*}
 \Qbm(\mubold,\,t) = 
 \begin{bmatrix}
  \cos \theta(\mubold,\,t) & \sin \theta(\mubold,\,t) \\
  -\sin \theta(\mubold,\,t) & \cos \theta(\mubold,\,t) \\
 \end{bmatrix} \qquad \qquad
 \vbm(\mubold,\,t) = 
 \begin{bmatrix}
  0 \\ h(\mubold,\,t)
 \end{bmatrix}.
\end{equation*}
The isentropic Navier-Stokes equations are discretized with the
discontinuous Galerkin scheme of Section~\ref{sec:app} using $978$ triangular
$p = 3$ elements. No-slip boundary conditions are imposed on the airfoil wall
and characteristic free-stream boundary conditions at the far-field. The
temporal discretization uses a third-order diagonally implicit Runge-Kutta
solver with $100$ equally spaced steps to discretize a single period of the
motion.
The airfoil and surrounding fluid vorticity field are shown in
Figures~\ref{fig:from-ss}~and~\ref{fig:from-periodic} with the flow field
initialized from steady-state flow and the time-periodic initial condition,
respectively. It is clear that the flow in Figure~\ref{fig:from-periodic} will
seamlessly transition between periods. The initialization from the steady-state
solution in Figure~\ref{fig:from-ss} will introduce non-physcial transients
into the flow as discussed in the next section.

First, the solvers introduced in Section~\ref{subsec:tpde-soln} are compared
for different initial guesses for the time-periodic initial condition. In the
absence of any a-priori information regarding the time-periodic solution, a
reasonable initial guess is the steady-state flow. Since the problem under
consideration is being forced by an input -- the periodic motion of the
foil -- a mechanism for improving the initial guess is to simulate the flow
field for $m$ periods of the foil motion and use the final state of the final
period as the initial guess. This corresponds to using $m$ iterations of
fixed point iteration (Algorithm~\ref{alg:fixedpt}) as a nonlinear
preconditioner for the nonlinear system of equations (\ref{eqn:fixedpt}) that
enforces time-periodicity of the flow.

\begin{figure}
 \centering
  \begin{subfigure}{0.28\textwidth}
    \centering
    \begin{tikzpicture}

\begin{axis}[
name=plt1,
scale only axis,
width=0.9\textwidth,
height=0.2\textheight,
xmax=100,
ymin=1e-10,
ymax=1e3,
xmode=log,
ymode=log,
xlabel=iterations (primal solves),
ylabel=$\norm{\ubm^{(N_t)} - \ubm_0}_2$]
\addplot [black, solid, thick, mark=*, mark repeat=10]  table[x expr=\coordindex+1, y index=1] {foil2d/cfd/flap00/periodic/prim/fixedpt/testcase00/ascii/app1_mesh-how1-0-3_cfd-1p4-0p2-1000-0p72-0-1-_1-0__tdisc2-0p2-100-_1-1_-3.fixedpt.conv};
\addplot [red, solid, thick, mark=o]  table[x expr=\coordindex+1, y index=1] {foil2d/cfd/flap00/periodic/prim/newtgmres/testcase04/neq0/ascii/app1_mesh-how1-0-3_cfd-1p4-0p2-1000-0p72-0-1-_1-0__tdisc2-0p2-100-_1-1_-3.newtraph.conv};
\addplot [blue, solid, thick, mark=x]  table[x expr=\coordindex+1, y index=1] {foil2d/cfd/flap00/periodic/prim/newtgmres/testcase05/neq0/ascii/app1_mesh-how1-0-3_cfd-1p4-0p2-1000-0p72-0-1-_1-0__tdisc2-0p2-100-_1-1_-3.newtraph.conv};
\addplot [magenta, solid, thick, mark=+]  table[x expr=\coordindex+1, y index=1] {foil2d/cfd/flap00/periodic/prim/newtgmres/testcase06/neq0/ascii/app1_mesh-how1-0-3_cfd-1p4-0p2-1000-0p72-0-1-_1-0__tdisc2-0p2-100-_1-1_-3.newtraph.conv};
\end{axis}
\end{tikzpicture}
  \end{subfigure} \hfill \hspace{1.5cm}
  \begin{subfigure}{0.28\textwidth}
    \centering
    \begin{tikzpicture}

\begin{axis}[
name=plt1,
scale only axis,
width=0.9\textwidth,
height=0.2\textheight,
xmax=100,
ymin=1e-10,
ymax=1e3,
xmode=log,
ymode=log,
yticklabels={,,},
xlabel=iterations (primal solves)]
\addplot [black, solid, thick, mark=*, mark repeat=10]  table[skip first n=1, x expr=\coordindex+1, y index=1] {foil2d/cfd/flap00/periodic/prim/fixedpt/testcase00/ascii/app1_mesh-how1-0-3_cfd-1p4-0p2-1000-0p72-0-1-_1-0__tdisc2-0p2-100-_1-1_-3.fixedpt.conv};
\addplot [cyan, solid, thick, mark=square*, mark repeat=10]  table[x expr=\coordindex+1, y index=1] {foil2d/cfd/flap00/periodic/prim/steepest/testcase00/neq1/ascii/app0_mesh-how1-0-3_cfd-1p4-0p2-1000-0p72-0-1-_1-0__tdisc2-0p2-100-_1-1_-3.steepest.conv};
\addplot [yellow!90!black, solid, thick, mark=square*, mark repeat=10]  table[x expr=\coordindex+1, y index=1] {foil2d/cfd/flap00/periodic/prim/lbfgs/testcase00/neq1/ascii/app0_mesh-how1-0-3_cfd-1p4-0p2-1000-0p72-0-1-_1-0__tdisc2-0p2-100-_1-1_-3.lbfgs.conv};
\addplot [red, solid, thick, mark=o]  table[x expr=\coordindex+1, y index=1] {foil2d/cfd/flap00/periodic/prim/newtgmres/testcase04/neq1/ascii/app1_mesh-how1-0-3_cfd-1p4-0p2-1000-0p72-0-1-_1-0__tdisc2-0p2-100-_1-1_-3.newtraph.conv};
\addplot [blue, solid, thick, mark=x]  table[x expr=\coordindex+1, y index=1] {foil2d/cfd/flap00/periodic/prim/newtgmres/testcase05/neq1/ascii/app1_mesh-how1-0-3_cfd-1p4-0p2-1000-0p72-0-1-_1-0__tdisc2-0p2-100-_1-1_-3.newtraph.conv};
\addplot [magenta, solid, thick, mark=+]  table[x expr=\coordindex+1, y index=1] {foil2d/cfd/flap00/periodic/prim/newtgmres/testcase06/neq1/ascii/app1_mesh-how1-0-3_cfd-1p4-0p2-1000-0p72-0-1-_1-0__tdisc2-0p2-100-_1-1_-3.newtraph.conv};
\end{axis}
\end{tikzpicture}
  \end{subfigure} \hfill
  \begin{subfigure}{0.28\textwidth}
    \centering
    \begin{tikzpicture}

\begin{axis}[
name=plt1,
scale only axis,
width=0.9\textwidth,
height=0.2\textheight,
xmax=100,
ymin=1e-10,
ymax=1e3,
xmode=log,
ymode=log,
yticklabels={,,},
xlabel=iterations (primal solves)]
\addplot [black, solid, thick, mark=*, mark repeat=10]  table[skip first n=5, x expr=\coordindex+1, y index=1] {foil2d/cfd/flap00/periodic/prim/fixedpt/testcase00/ascii/app1_mesh-how1-0-3_cfd-1p4-0p2-1000-0p72-0-1-_1-0__tdisc2-0p2-100-_1-1_-3.fixedpt.conv}; \label{line:fixedpt}
\addplot [cyan, solid, thick, mark=triangle*, mark repeat=10]  table[x expr=\coordindex+1, y index=1] {foil2d/cfd/flap00/periodic/prim/steepest/testcase00/neq5/ascii/app0_mesh-how1-0-3_cfd-1p4-0p2-1000-0p72-0-1-_1-0__tdisc2-0p2-100-_1-1_-3.steepest.conv}; \label{line:steepest}
\addplot [yellow!90!black, solid, thick, mark=square*, mark repeat=10]  table[x expr=\coordindex+1, y index=1] {foil2d/cfd/flap00/periodic/prim/lbfgs/testcase00/neq5/ascii/app0_mesh-how1-0-3_cfd-1p4-0p2-1000-0p72-0-1-_1-0__tdisc2-0p2-100-_1-1_-3.lbfgs.conv}; \label{line:lbfgs}
\addplot [red, solid, thick, mark=o]  table[x expr=\coordindex+1, y index=1] {foil2d/cfd/flap00/periodic/prim/newtgmres/testcase04/neq5/ascii/app1_mesh-how1-0-3_cfd-1p4-0p2-1000-0p72-0-1-_1-0__tdisc2-0p2-100-_1-1_-3.newtraph.conv}; \label{line:newtgmres-2}
\addplot [blue, solid, thick, mark=x]  table[x expr=\coordindex+1, y index=1] {foil2d/cfd/flap00/periodic/prim/newtgmres/testcase05/neq5/ascii/app1_mesh-how1-0-3_cfd-1p4-0p2-1000-0p72-0-1-_1-0__tdisc2-0p2-100-_1-1_-3.newtraph.conv}; \label{line:newtgmres-3}
\addplot [magenta, solid, thick, mark=+]  table[x expr=\coordindex+1, y index=1] {foil2d/cfd/flap00/periodic/prim/newtgmres/testcase06/neq5/ascii/app1_mesh-how1-0-3_cfd-1p4-0p2-1000-0p72-0-1-_1-0__tdisc2-0p2-100-_1-1_-3.newtraph.conv}; \label{line:newtgmres-4}
\end{axis}
\end{tikzpicture}
  \end{subfigure}
\caption{Convergence comparison for numerical solvers for fully discrete
         time-periodically constrained partial differential equations
         (\ref{eqn:dirk}), nonlinearly preconditioned with $m$ fixed point
         iterations. Left: $m = 0$, middle: $m = 1$, right: $m = 5$.
         Solvers: fixed point iteration (\ref{line:fixedpt}),
         steepest decent (\ref{line:steepest}),
         L-BFGS (\ref{line:lbfgs}),
         Newton-GMRES: $\epsilon = 10^{-2}$ (\ref{line:newtgmres-2}),
         $\epsilon = 10^{-3}$ (\ref{line:newtgmres-3}),
         $\epsilon = 10^{-4}$ (\ref{line:newtgmres-4}). The optimization
         algorithms (steepest decent and L-BFGS) were not included in the
         $m = 0$ study due to lack of convergence issues.}
\label{fig:solvers1}
\end{figure}
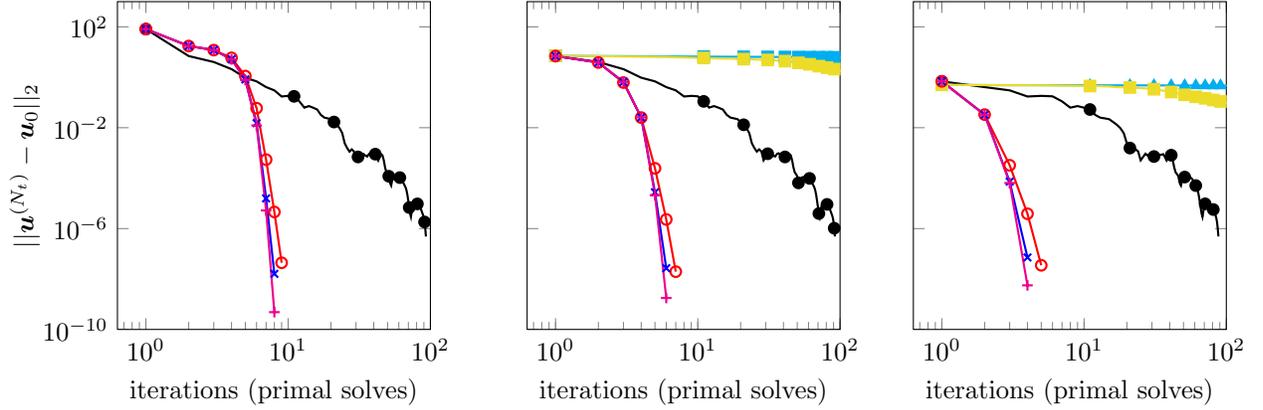

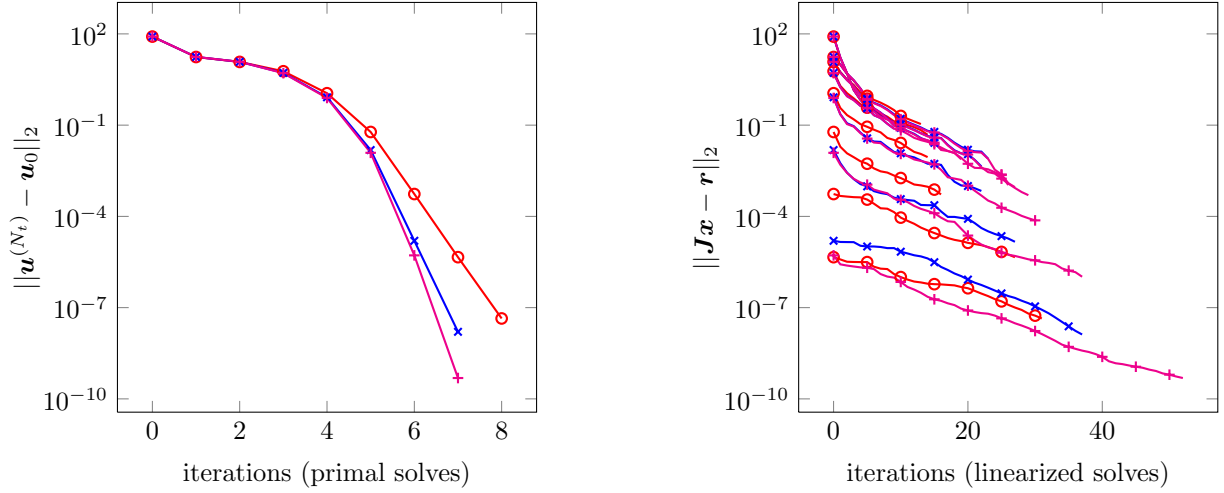
\begin{figure}
 \begin{subfigure}{0.45\textwidth}
  \begin{tikzpicture}

\begin{axis}[
name=plt1,
scale only axis,
width=0.75\textwidth,
height=0.25\textheight,
ymode=log,
xlabel=iterations (primal solves),
ylabel=$\norm{\ubm^{(N_t)} - \ubm_0}_2$]
\addplot [red, solid, thick, mark=o]  file {foil2d/cfd/flap00/periodic/prim/newtgmres/testcase04/ascii/app1_mesh-how1-0-3_cfd-1p4-0p2-1000-0p72-0-1-_1-0__tdisc2-0p2-100-_1-1_-3.newtraph.conv};
\addplot [blue, solid, thick, mark=x]  file {foil2d/cfd/flap00/periodic/prim/newtgmres/testcase05/ascii/app1_mesh-how1-0-3_cfd-1p4-0p2-1000-0p72-0-1-_1-0__tdisc2-0p2-100-_1-1_-3.newtraph.conv};
\addplot [magenta, solid, thick, mark=+]  file {foil2d/cfd/flap00/periodic/prim/newtgmres/testcase06/ascii/app1_mesh-how1-0-3_cfd-1p4-0p2-1000-0p72-0-1-_1-0__tdisc2-0p2-100-_1-1_-3.newtraph.conv};
\end{axis}
\end{tikzpicture}
 \end{subfigure}\hfill
 \begin{subfigure}{0.45\textwidth}
  \begin{tikzpicture}

\begin{axis}[
name=plt1,
scale only axis,
width=0.75\textwidth,
height=0.25\textheight,
ymode=log,
xlabel=iterations (linearized solves),
ylabel=$\norm{\Jbm\xbm - \rbm}_2$]

\addplot [red, solid, thick, mark=o, mark repeat=5, mark repeat=5]  file {foil2d/cfd/flap00/periodic/prim/newtgmres/testcase04/ascii/app1_mesh-how1-0-3_cfd-1p4-0p2-1000-0p72-0-1-_1-0__tdisc2-0p2-100-_1-1_-3.gmres0.conv};
\addplot [red, solid, thick, mark=o, mark repeat=5]  file {foil2d/cfd/flap00/periodic/prim/newtgmres/testcase04/ascii/app1_mesh-how1-0-3_cfd-1p4-0p2-1000-0p72-0-1-_1-0__tdisc2-0p2-100-_1-1_-3.gmres1.conv};
\addplot [red, solid, thick, mark=o, mark repeat=5]  file {foil2d/cfd/flap00/periodic/prim/newtgmres/testcase04/ascii/app1_mesh-how1-0-3_cfd-1p4-0p2-1000-0p72-0-1-_1-0__tdisc2-0p2-100-_1-1_-3.gmres2.conv};
\addplot [red, solid, thick, mark=o, mark repeat=5]  file {foil2d/cfd/flap00/periodic/prim/newtgmres/testcase04/ascii/app1_mesh-how1-0-3_cfd-1p4-0p2-1000-0p72-0-1-_1-0__tdisc2-0p2-100-_1-1_-3.gmres3.conv};
\addplot [red, solid, thick, mark=o, mark repeat=5]  file {foil2d/cfd/flap00/periodic/prim/newtgmres/testcase04/ascii/app1_mesh-how1-0-3_cfd-1p4-0p2-1000-0p72-0-1-_1-0__tdisc2-0p2-100-_1-1_-3.gmres4.conv};
\addplot [red, solid, thick, mark=o, mark repeat=5]  file {foil2d/cfd/flap00/periodic/prim/newtgmres/testcase04/ascii/app1_mesh-how1-0-3_cfd-1p4-0p2-1000-0p72-0-1-_1-0__tdisc2-0p2-100-_1-1_-3.gmres5.conv};
\addplot [red, solid, thick, mark=o, mark repeat=5]  file {foil2d/cfd/flap00/periodic/prim/newtgmres/testcase04/ascii/app1_mesh-how1-0-3_cfd-1p4-0p2-1000-0p72-0-1-_1-0__tdisc2-0p2-100-_1-1_-3.gmres6.conv};
\addplot [red, solid, thick, mark=o, mark repeat=5]  file {foil2d/cfd/flap00/periodic/prim/newtgmres/testcase04/ascii/app1_mesh-how1-0-3_cfd-1p4-0p2-1000-0p72-0-1-_1-0__tdisc2-0p2-100-_1-1_-3.gmres7.conv};

\addplot [blue, solid, thick, mark=x, mark repeat=5]  file {foil2d/cfd/flap00/periodic/prim/newtgmres/testcase05/ascii/app1_mesh-how1-0-3_cfd-1p4-0p2-1000-0p72-0-1-_1-0__tdisc2-0p2-100-_1-1_-3.gmres0.conv};
\addplot [blue, solid, thick, mark=x, mark repeat=5]  file {foil2d/cfd/flap00/periodic/prim/newtgmres/testcase05/ascii/app1_mesh-how1-0-3_cfd-1p4-0p2-1000-0p72-0-1-_1-0__tdisc2-0p2-100-_1-1_-3.gmres1.conv};
\addplot [blue, solid, thick, mark=x, mark repeat=5]  file {foil2d/cfd/flap00/periodic/prim/newtgmres/testcase05/ascii/app1_mesh-how1-0-3_cfd-1p4-0p2-1000-0p72-0-1-_1-0__tdisc2-0p2-100-_1-1_-3.gmres2.conv};
\addplot [blue, solid, thick, mark=x, mark repeat=5]  file {foil2d/cfd/flap00/periodic/prim/newtgmres/testcase05/ascii/app1_mesh-how1-0-3_cfd-1p4-0p2-1000-0p72-0-1-_1-0__tdisc2-0p2-100-_1-1_-3.gmres3.conv};
\addplot [blue, solid, thick, mark=x, mark repeat=5]  file {foil2d/cfd/flap00/periodic/prim/newtgmres/testcase05/ascii/app1_mesh-how1-0-3_cfd-1p4-0p2-1000-0p72-0-1-_1-0__tdisc2-0p2-100-_1-1_-3.gmres4.conv};
\addplot [blue, solid, thick, mark=x, mark repeat=5]  file {foil2d/cfd/flap00/periodic/prim/newtgmres/testcase05/ascii/app1_mesh-how1-0-3_cfd-1p4-0p2-1000-0p72-0-1-_1-0__tdisc2-0p2-100-_1-1_-3.gmres5.conv};
\addplot [blue, solid, thick, mark=x, mark repeat=5]  file {foil2d/cfd/flap00/periodic/prim/newtgmres/testcase05/ascii/app1_mesh-how1-0-3_cfd-1p4-0p2-1000-0p72-0-1-_1-0__tdisc2-0p2-100-_1-1_-3.gmres6.conv};

\addplot [magenta, solid, thick, mark=+, mark repeat=5]  file {foil2d/cfd/flap00/periodic/prim/newtgmres/testcase06/ascii/app1_mesh-how1-0-3_cfd-1p4-0p2-1000-0p72-0-1-_1-0__tdisc2-0p2-100-_1-1_-3.gmres0.conv};
\addplot [magenta, solid, thick, mark=+, mark repeat=5]  file {foil2d/cfd/flap00/periodic/prim/newtgmres/testcase06/ascii/app1_mesh-how1-0-3_cfd-1p4-0p2-1000-0p72-0-1-_1-0__tdisc2-0p2-100-_1-1_-3.gmres1.conv};
\addplot [magenta, solid, thick, mark=+, mark repeat=5]  file {foil2d/cfd/flap00/periodic/prim/newtgmres/testcase06/ascii/app1_mesh-how1-0-3_cfd-1p4-0p2-1000-0p72-0-1-_1-0__tdisc2-0p2-100-_1-1_-3.gmres2.conv};
\addplot [magenta, solid, thick, mark=+, mark repeat=5]  file {foil2d/cfd/flap00/periodic/prim/newtgmres/testcase06/ascii/app1_mesh-how1-0-3_cfd-1p4-0p2-1000-0p72-0-1-_1-0__tdisc2-0p2-100-_1-1_-3.gmres3.conv};
\addplot [magenta, solid, thick, mark=+, mark repeat=5]  file {foil2d/cfd/flap00/periodic/prim/newtgmres/testcase06/ascii/app1_mesh-how1-0-3_cfd-1p4-0p2-1000-0p72-0-1-_1-0__tdisc2-0p2-100-_1-1_-3.gmres4.conv};
\addplot [magenta, solid, thick, mark=+, mark repeat=5]  file {foil2d/cfd/flap00/periodic/prim/newtgmres/testcase06/ascii/app1_mesh-how1-0-3_cfd-1p4-0p2-1000-0p72-0-1-_1-0__tdisc2-0p2-100-_1-1_-3.gmres5.conv};
\addplot [magenta, solid, thick, mark=+, mark repeat=5]  file {foil2d/cfd/flap00/periodic/prim/newtgmres/testcase06/ascii/app1_mesh-how1-0-3_cfd-1p4-0p2-1000-0p72-0-1-_1-0__tdisc2-0p2-100-_1-1_-3.gmres6.conv};
\end{axis}
\end{tikzpicture}
 \end{subfigure}
\caption{Linear and nonlinear convergence of Newton-GMRES method for
         determining fully discrete time-periodic solutions with
         various linear system tolerances, $\epsilon$, i.e.
         $\norm{\Jbm\xbm - \Rbm} < \epsilon$, where $\rbm$ and
         $\Jbm$ are defined in (\ref{eqn:nlsys}) and (\ref{eqn:jac-prim}).
         Tolerances considered:
         $\epsilon = 10^{-2}$ (\ref{line:newtgmres-2}),
         $\epsilon = 10^{-3}$ (\ref{line:newtgmres-3}),
         $\epsilon = 10^{-4}$ (\ref{line:newtgmres-4}).}
\label{fig:solvers2}
\end{figure}

Figure~\ref{fig:solvers1} and Table~\ref{tab:summ} compare the solvers
under consideration for different levels of nonlinear preconditioning.
The iteration count is a good metric for comparison as the
  bulk of the time in any of these methods is timestepping the forward
  (primal) problem, the sensitivity equations, and/or the adjoint
  equations, which are the same from one method to the next. In our
  implementation, the linearized equations (sensitivity and adjoint) are
  about $2\times$ less expensive to solve than the nonlinear, primal equations
  -- see discussion in Section~\ref{subsec:tpde-soln}. The linear algebra
  involved inside e.g.~the L-BFGS and GMRES algorithms is negligible in
  comparison to these timestepping costs.
Regardless of nonlinear preconditioning, the Newton-GMRES solver converges
most rapidly for a range of linear system tolerances from $10^{-2}$ to
$10^{-4}$ and the optimization algorithms (steepest descent and L-BFGS)
converge most slowly. In fact, without any nonlinear preconditioning the
optimization algorithms fail to make progress toward the optimal solution and
were not included in the figure. Nonlinear preconditioning helps the
Newton-GMRES algorithm most substantially, particularly with $m = 5$, as this
appears to place the initial guess close enough to the solution that quadratic
convergence is obtained from the outset. This causes the number of Newton
iterations to be reduced from $8$ or $9$ to $3$ or $4$. From
Table~\ref{tab:summ}, this does not save many primal solvers -- since the
nonlinear preconditioning requires primal solves -- but requires far fewer
linear system iterations and therefore fewer sensitivity solutions.
Figure~\ref{fig:solvers2} isolates the Newton-GMRES solver (for $m = 0$,
i.e., the case without preconditioning) to highlight convergence rates for
different GMRES tolerances. It also shows the convergence of GMRES for each
nonlinear iteration and each tolerance considered. As expected, more GMRES
iterations are required near convergence as it becomes more difficult to
reduce the linear residual the prescribed orders of magnitude.

\begin{table}
\centering
\begin{tabular}{l c c c c} \toprule
$m = 0$ & $\norm{\ubm^{(N_t)} - \ubm_0}_2$ & Primal Solves (\ref{eqn:dirk2})
                 & Sensitivity Solves (\ref{eqn:dirk-sens})
                 & Adjoint Solves (\ref{eqn:uns-disc-adj-dirk}) \\\midrule
\rowcolor{black!00} Fixed Point Iteration     & 8.10e-07 & 90 & 0 & 0 \\ 
\rowcolor{black!10} Newton-Krylov ($10^{-2}$) & 4.41e-08 & 9 & 128 & 0 \\
\rowcolor{black!00} Newton-Krylov ($10^{-3}$) & 1.60e-08 & 8 & 156 & 0 \\
\rowcolor{black!10} Newton-Krylov ($10^{-4}$) & 4.85e-10 & 8 & 220 & 0 \\%
                                                                       \bottomrule
$m = 1$ & $\norm{\ubm^{(N_t)} - \ubm_0}_2$ & Primal Solves (\ref{eqn:dirk2})
                 & Sensitivity Solves (\ref{eqn:dirk-sens})
                 & Adjoint Solves (\ref{eqn:uns-disc-adj-dirk}) \\\midrule
\rowcolor{black!00} Fixed Point Iteration     & 8.10e-07 & 90 & 0 & 0 \\ 
\rowcolor{black!10} Steepest Decent           & 6.09e+00 & 121 & 0 & 121 \\ 
\rowcolor{black!00} L-BFGS                    & 1.36e+00 & 121 & 0 & 121 \\ 
\rowcolor{black!10} Newton-Krylov ($10^{-2}$) & 1.96e-08 & 8 & 104 & 0 \\
\rowcolor{black!00} Newton-Krylov ($10^{-3}$) & 2.69e-08 & 7 & 116 & 0 \\
\rowcolor{black!10} Newton-Krylov ($10^{-4}$) & 1.77e-09 & 7 & 149 & 0 \\%
                                                                       \bottomrule
$m = 5$ & $\norm{\ubm^{(N_t)} - \ubm_0}_2$ & Primal Solves (\ref{eqn:dirk2})
                 & Sensitivity Solves (\ref{eqn:dirk-sens})
                 & Adjoint Solves (\ref{eqn:uns-disc-adj-dirk}) \\\midrule
\rowcolor{black!00} Fixed Point Iteration     & 8.10e-07 & 90 & 0 & 0 \\ 
\rowcolor{black!10} Steepest Decent           & 4.65e-01 & 125 & 0 & 125 \\ 
\rowcolor{black!00} L-BFGS                    & 7.40e-02 & 125 & 0 & 125 \\ 
\rowcolor{black!10} Newton-Krylov ($10^{-2}$) & 3.50e-08 & 10 & 92 & 0 \\
\rowcolor{black!00} Newton-Krylov ($10^{-3}$) & 7.18e-08 & 9 & 88 & 0 \\
\rowcolor{black!10} Newton-Krylov ($10^{-4}$) & 5.61e-09 & 9 & 121 & 0 \\%
                                                                       \bottomrule
\end{tabular}

\caption{Table summarizing performance of numerical solvers for
         fully discrete time-periodic partial differential equations,
         considering nonlinear preconditioning via $m$ fixed point iterations.}
\label{tab:summ}
\end{table}

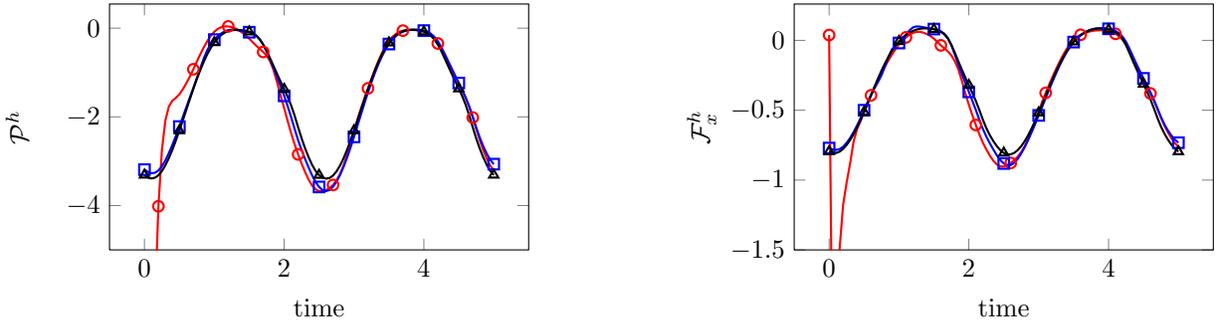
\begin{figure}
 \centering
 \begin{subfigure}{0.45\textwidth}
  \centering
  \begin{tikzpicture}

\begin{axis}[
  scale only axis,
  width=0.75\textwidth,
  height=0.15\textheight,
  ymin=-5,
  xlabel=time,
  ylabel={$\Pcal^h$}]
\addplot [red, solid, thick, mark=o, mark repeat=10]  table[x expr=0.05*\coordindex, y index=3] {foil2d/cfd/flap00/periodic/prim/newtgmres/testcase04/neq0/ascii/app1_mesh-how1-0-3_cfd-1p4-0p2-1000-0p72-0-1-_1-0__tdisc2-0p2-100-_1-1_-3_it000.allpost.instant}; \label{line:keq0}
\addplot [blue, solid, thick, mark=square, mark repeat=10]  table[x expr=0.05*\coordindex, y index=3] {foil2d/cfd/flap00/periodic/prim/newtgmres/testcase04/neq0/ascii/app1_mesh-how1-0-3_cfd-1p4-0p2-1000-0p72-0-1-_1-0__tdisc2-0p2-100-_1-1_-3_it001.allpost.instant};\label{line:keq1}
\addplot [black, solid, thick, mark=triangle, mark repeat=10]  table[x expr=0.05*\coordindex, y index=3] {foil2d/cfd/flap00/periodic/prim/newtgmres/testcase04/neq0/ascii/app1_mesh-how1-0-3_cfd-1p4-0p2-1000-0p72-0-1-_1-0__tdisc2-0p2-100-_1-1_-3_it008.allpost.instant};\label{line:keq8}
\end{axis}
\end{tikzpicture}
 \end{subfigure} \hfill
 \begin{subfigure}{0.45\textwidth}
  \centering
  \begin{tikzpicture}

\begin{axis}[
  scale only axis,
  width=0.75\textwidth,
  height=0.15\textheight,
  ymin=-1.5,
  xlabel=time,
  ylabel={$\Fcal_x^h$}]
\addplot [red, solid, thick, mark=o, mark repeat=10]  table[x expr=0.05*\coordindex, y index=0] {foil2d/cfd/flap00/periodic/prim/newtgmres/testcase04/neq0/ascii/app1_mesh-how1-0-3_cfd-1p4-0p2-1000-0p72-0-1-_1-0__tdisc2-0p2-100-_1-1_-3_it000.allpost.instant};
\addplot [blue, solid, thick, mark=square, mark repeat=10]  table[x expr=0.05*\coordindex, y index=0] {foil2d/cfd/flap00/periodic/prim/newtgmres/testcase04/neq0/ascii/app1_mesh-how1-0-3_cfd-1p4-0p2-1000-0p72-0-1-_1-0__tdisc2-0p2-100-_1-1_-3_it001.allpost.instant};
\addplot [black, solid, thick, mark=triangle, mark repeat=10]  table[x expr=0.05*\coordindex, y index=0] {foil2d/cfd/flap00/periodic/prim/newtgmres/testcase04/neq0/ascii/app1_mesh-how1-0-3_cfd-1p4-0p2-1000-0p72-0-1-_1-0__tdisc2-0p2-100-_1-1_-3_it008.allpost.instant};
\end{axis}
\end{tikzpicture}
 \end{subfigure}
\caption{Time history of power, $\Fcal_x^h(\ubm, \mubold, t)$, and $x$-directed
         force, $\Pcal^h(\ubm, \mubold, t)$, after $k$ Newton-GMRES iterations
         ($\epsilon = 10^{-2}$) starting from steady-state. Values of $k$: 
         $0$ (\ref{line:keq0}), $1$ (\ref{line:keq1}), and
         $8$ (\ref{line:keq8}).}
\label{fig:qoi-hist1}
\end{figure}

\begin{figure}
 \centering
 \begin{subfigure}{0.45\textwidth}
  \centering
  \begin{tikzpicture}

\begin{axis}[
  scale only axis,
  width=0.75\textwidth,
  height=0.20\textheight,
  ymode=log,
  xlabel=iteration,
  ylabel=$|W - W^*|$
]
\addplot [red, solid, thick, mark=o, mark repeat=1]  table[x index=0, y expr=abs(\thisrowno{4}+7.404719751811117412)] {foil2d/cfd/flap00/periodic/prim/newtgmres/testcase04/neq0/ascii/app1_mesh-how1-0-3_cfd-1p4-0p2-1000-0p72-0-1-_1-0__tdisc2-0p2-100-_1-1_-3.allpost.integrated};
\addplot [blue, solid, thick, mark=x, mark repeat=1]  table[x index=0, y expr=abs(\thisrowno{4}+7.404719751811117412)] {foil2d/cfd/flap00/periodic/prim/newtgmres/testcase05/neq0/ascii/app1_mesh-how1-0-3_cfd-1p4-0p2-1000-0p72-0-1-_1-0__tdisc2-0p2-100-_1-1_-3.allpost.integrated};
\addplot [magenta, solid, thick, mark=+, mark repeat=1]  table[x index=0, y expr=abs(\thisrowno{4}+7.404719751811117412)] {foil2d/cfd/flap00/periodic/prim/newtgmres/testcase06/neq0/ascii/app1_mesh-how1-0-3_cfd-1p4-0p2-1000-0p72-0-1-_1-0__tdisc2-0p2-100-_1-1_-3.allpost.integrated};
\end{axis}
\end{tikzpicture}
 \end{subfigure} \hfill
 \begin{subfigure}{0.45\textwidth}
  \centering
  \begin{tikzpicture}

\begin{axis}[
  scale only axis,
  width=0.75\textwidth,
  height=0.20\textheight,
  ymode=log,
  xlabel=iteration,
  ylabel=$|J_x - J_x^*|$
]
\addplot [red, solid, thick, mark=o, mark repeat=1]  table[x index=0, y expr=abs(\thisrowno{1}+1.555146353340866483)] {foil2d/cfd/flap00/periodic/prim/newtgmres/testcase04/neq0/ascii/app1_mesh-how1-0-3_cfd-1p4-0p2-1000-0p72-0-1-_1-0__tdisc2-0p2-100-_1-1_-3.allpost.integrated};
\addplot [blue, solid, thick, mark=x, mark repeat=1]  table[x index=0, y expr=abs(\thisrowno{1}+1.555146353340866483)] {foil2d/cfd/flap00/periodic/prim/newtgmres/testcase05/neq0/ascii/app1_mesh-how1-0-3_cfd-1p4-0p2-1000-0p72-0-1-_1-0__tdisc2-0p2-100-_1-1_-3.allpost.integrated};
\addplot [magenta, solid, thick, mark=+, mark repeat=1]  table[x index=0, y expr=abs(\thisrowno{1}+1.555146353340866483)] {foil2d/cfd/flap00/periodic/prim/newtgmres/testcase06/neq0/ascii/app1_mesh-how1-0-3_cfd-1p4-0p2-1000-0p72-0-1-_1-0__tdisc2-0p2-100-_1-1_-3.allpost.integrated};
\end{axis}
\end{tikzpicture}
 \end{subfigure}
\caption{Convergence of fully discrete quantities of interest to their
         values at the time-periodic solution, $W^*$ and $J_x^*$, for
         various solvers, without nonlinear preconditioning. Solvers:
         Newton-GMRES: $\epsilon = 10^{-2}$ (\ref{line:newtgmres-2}),
         $\epsilon = 10^{-3}$ (\ref{line:newtgmres-3}),
         $\epsilon = 10^{-4}$ (\ref{line:newtgmres-4}).}
\label{fig:qoi-conv1}
\end{figure}
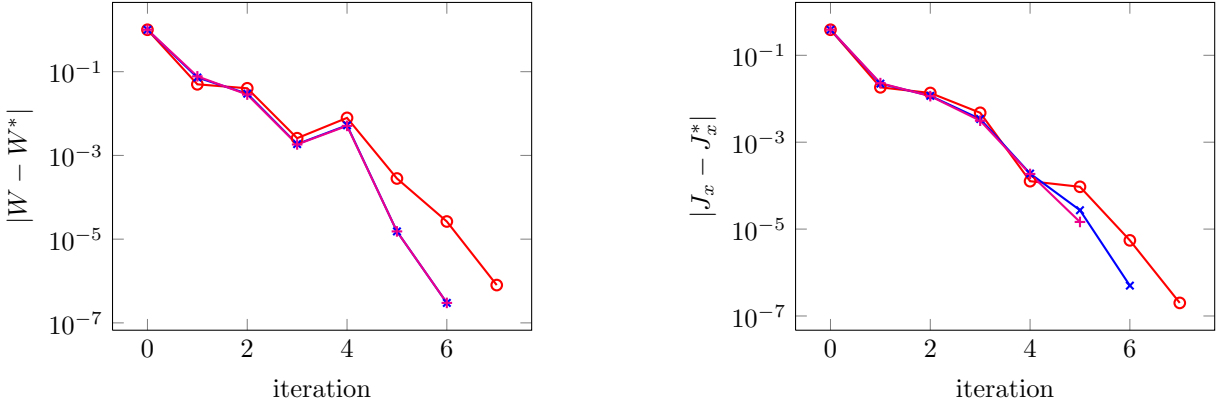

The time history of the instantaneous quantities of interest in
Figure~\ref{fig:qoi-hist1} illustrate the non-physical transients that
result from initializing the flow with the steady-state solution. While the
transients mostly vanish after a single Newton iteration, the trajectories of
these quantities of interest do not coincide with those of the true
time-periodic solution. The error between the integrated quantities
of interest -- $W$ and $J_x$ -- at the time-periodic flow versus intermediate
iterations is shown in Figure~\ref{fig:qoi-conv1}. Comparing
Figures~\ref{fig:solvers1}~and~\ref{fig:qoi-conv1}, it can be seen
that a tolerance of $10^{-8}$ on $\norm{\ubm^{(N_t)} - \ubm^{(0)}}_2$ leads
to an accuracy of $10^{-6}$ in the integrated quantities of the time-periodic
solution.

\begin{figure}[h]
  \centering
  \begin{subfigure}{0.45\textwidth}
    \centering
    \begin{tikzpicture}

\begin{axis}[
  scale only axis,
  width=0.9\textwidth,
  height=0.25\textheight,
  xmin=-1,
  xmax=1,
  ymin=-1,
  ymax=1,
  axis equal,
  xlabel=$\Re(\lambda)$,
  ylabel=$\Im(\lambda)$]
\addplot [black, solid, thick, mark=o, only marks]  table[select coords between index={0}{200}, x index=0, y index=1] {foil2d/cfd/flap00/pstab/testcase00/ascii/app0_mesh-how1-0-1_cfd-1p4-0p15-100-0p72-0-1-_1-0__tdisc2-0p2-50-_1-1_-3.stabmat.evals}; \label{line:evals}
\addplot [blue, thick] table[x index=0, y index=1] {dat/unit_circ.dat};
\end{axis}

\end{tikzpicture}
  \end{subfigure}
  \caption{First 200 eigenvalues (\ref{line:evals}) of
           $\pder{\ubm^{(N_t)}}{\ubm_0}$ -- evaluated at periodic solution --
           with largest magnitude. All eigenvalues lie in unit circle, thus
           the periodic orbit is stable.} \label{fig:stab}
\end{figure}
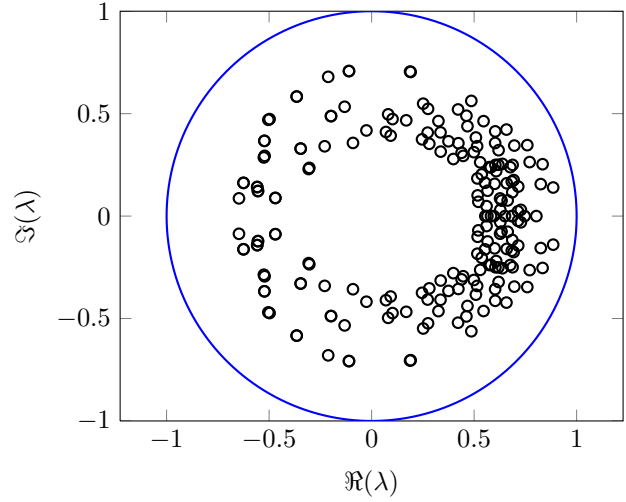

Next, the stability of the periodic orbit is verified by considering the
eigenvalues of $\displaystyle{\pder{\ubm^{(N_t)}}{\ubm_0}}$, evaluated at the
time-periodic solution. As discussed in Section~\ref{subsec:tpde-stab} and
many prior works \citep{coddington, kuchment2012floquet}, the periodic orbit
is stable if all eigenvalues of this matrix have modulus less than unity.
Figure~\ref{fig:stab} shows that the 200 eigenvalues of largest modulus lie
within the unit circle in the complex plane; thus, the periodic orbit is
stable for this problem.

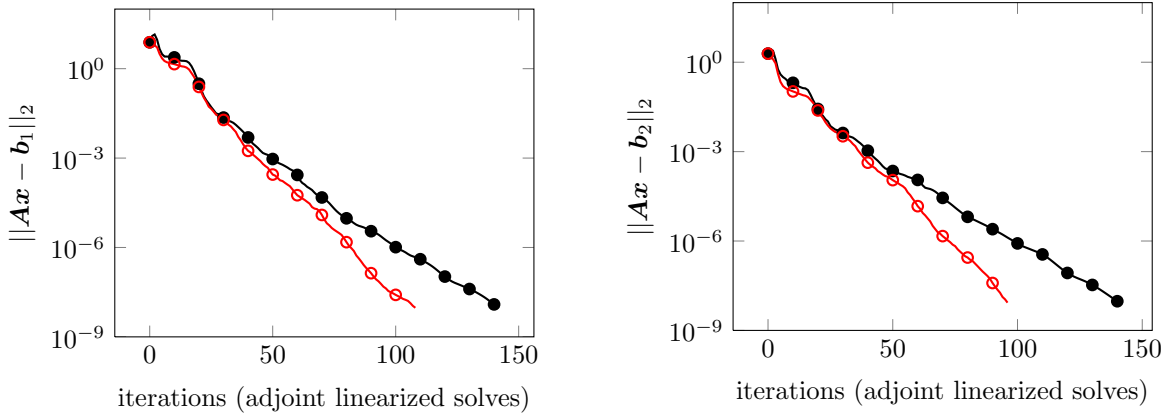
\begin{figure}
 \centering \hfill
 \begin{subfigure}{0.45\textwidth}
   \centering
   \begin{tikzpicture}

\begin{axis}[
name=plt1,
scale only axis,
width=0.75\textwidth,
height=0.20\textheight,
ymin = 1e-9,
ymax = 1e2,
ymode=log,
xlabel=iterations (adjoint linearized solves),
ylabel=$\norm{\Abm\xbm - \bbm_1}_2$]

\addplot [black, solid, thick, mark=*, mark repeat=10]  table[x index=0, y index=1] {foil2d/cfd/flap00/periodic/dual/totpow/fixedpt/testcase00/neq0/ascii/app0_mesh-how1-0-3_cfd-1p4-0p2-1000-0p72-0-1-_1-0__tdisc2-0p2-100-_1-1_-3.fixedpt.conv}; \label{line:dual-fixedpt}
\addplot [red, solid, thick, mark=o, mark repeat=10]  file {foil2d/cfd/flap00/periodic/dual/totpow/gmres/testcase00/neq0/ascii/app0_mesh-how1-0-3_cfd-1p4-0p2-1000-0p72-0-1-_1-0__tdisc2-0p2-100-_1-1_-3.gmres0.conv}; \label{line:dual-gmres}
\end{axis}
\end{tikzpicture}
 \end{subfigure}\hfill
 \begin{subfigure}{0.45\textwidth}
   \centering
   \vspace{2.5mm}
   \begin{tikzpicture}

\begin{axis}[
name=plt1,
scale only axis,
width=0.75\textwidth,
height=0.20\textheight,
ymin = 1e-9,
ymax = 1e2,
ymode=log,
xlabel=iterations (adjoint linearized solves),
ylabel=$\norm{\Abm\xbm - \bbm_2}_2$]

\addplot [black, solid, thick, mark=*, mark repeat=10]  table[x index=0, y index=1] {foil2d/cfd/flap00/periodic/dual/fx/fixedpt/testcase00/neq0/ascii/app0_mesh-how1-0-3_cfd-1p4-0p2-1000-0p72-0-1-_1-0__tdisc2-0p2-100-_1-1_-3.fixedpt.conv};
\addplot [red, solid, thick, mark=o, mark repeat=10]  file {foil2d/cfd/flap00/periodic/dual/fx/gmres/testcase00/neq0/ascii/app0_mesh-how1-0-3_cfd-1p4-0p2-1000-0p72-0-1-_1-0__tdisc2-0p2-100-_1-1_-3.gmres0.conv};
\end{axis}
\end{tikzpicture}
 \end{subfigure}
 \hfill
\caption{GMRES convergence for determining solution of adjoint equations
         corresponding to fully discrete time-periodic partial differential
         equation, i.e. a linear two-point boundary value problem. $\Abm$
         defined in (\ref{eqn:jac-dual}),
         $\displaystyle{\bbm_1 = \pder{W}{\ubm^{(N_t)}}}$, and
         $\displaystyle{\bbm_2 = \pder{J_x}{\ubm^{(N_t)}}}$ from
         (\ref{eqn:fixedpt-dual}), where $W$ is fully discrete approximation
         of the total work done by fluid on airfoil (Section~\ref{sec:app}) and
         $J_x$ is the x-directed impulse. Solvers:
         fixed point iteration (\ref{line:dual-fixedpt}) and
         GMRES (\ref{line:dual-gmres}).
         The linearization is performed about the
         time-periodic solution obtained with Newton-Krylov
         ($\epsilon = 10^{-4}$) method.}
\label{fig:solvers1-dual}
\end{figure}

This completes the discussion of the primal time-periodic problem and
attention is turned to the dual, or adjoint, problem. First, a brief
comparison of two potential solvers -- fixed point iteration and GMRES -- for
the periodic adjoint equation is provided. In contrast to the primal problem,
there is a less pronounced difference between the convergence of fixed point
iteration and the Krylov solver in the dual problem.
Figure~\ref{fig:solvers1-dual} shows the convergence history for two different
right-hand sides of $\Abm\xbm = \bbm$, each corresponding to the adjoint
method for a different quantity of interest. However, it should be noted that
the iterations for the GMRES solver are cheaper than those of the fixed point
solver as the terms $\displaystyle{\pder{F}{\ubm^{(n-1)}}}$ and
$\displaystyle{\pder{F}{\kbm_i^{(n)}}}$ are not computed.
Therefore, the
GMRES algorithm is superior to fixed point iterations as there are fewer
required iterations, each of which is cheaper.

\begin{figure}
 \centering
 \begin{subfigure}{0.45\textwidth}
  \centering
  \begin{tikzpicture}

\begin{axis}[
  scale only axis,
  width=0.75\textwidth,
  height=0.20\textheight,
  xmode=log,
  ymode=log,
  xlabel=$\tau$,
  ylabel=$\norm{\oder{W}{\mubold} - \frac{\Delta W}{\Delta \mubold}}_2 /
          \norm{\frac{\Delta W}{\Delta \mubold}}_2$
]
\addplot [black, solid, thick, mark=*, mark repeat=1]  table[x index=0, y index=1] {foil2d/cfd/flap00/pfindiff/testcase00/ascii/app0_mesh-how1-0-1_cfd-1p4-0p15-100-0p72-0-1-_1-0__tdisc2-0p2-50-_1-1_-3.rfindiff.cd2};
\end{axis}

\end{tikzpicture}
 \end{subfigure}\hfill
 \begin{subfigure}{0.45\textwidth}
  \centering
  \begin{tikzpicture}

\begin{axis}[
  scale only axis,
  width=0.75\textwidth,
  height=0.20\textheight,
  xmode=log,
  ymode=log,
  xlabel=$\tau$,
  ylabel=$\norm{\oder{J_x}{\mubold} - \frac{\Delta J_x}{\Delta \mubold}}_2 /
          \norm{\frac{\Delta J_x}{\Delta \mubold}}_2$
]
\addplot [black, solid, thick, mark=*, mark repeat=1]  table[x index=0, y index=2] {foil2d/cfd/flap00/pfindiff/testcase00/ascii/app0_mesh-how1-0-1_cfd-1p4-0p15-100-0p72-0-1-_1-0__tdisc2-0p2-50-_1-1_-3.rfindiff.cd2};
\end{axis}

\end{tikzpicture}
 \end{subfigure}
\caption{Verification of periodic adjoint-based gradient with second-order
         centered finite difference approximation, for a range of finite
         intervals, $\tau$. The computed gradient match the finite difference
         approximation to nearly $7$ digits before round-off errors degrade
         the accuracy.}
\label{fig:findiff}
\end{figure}
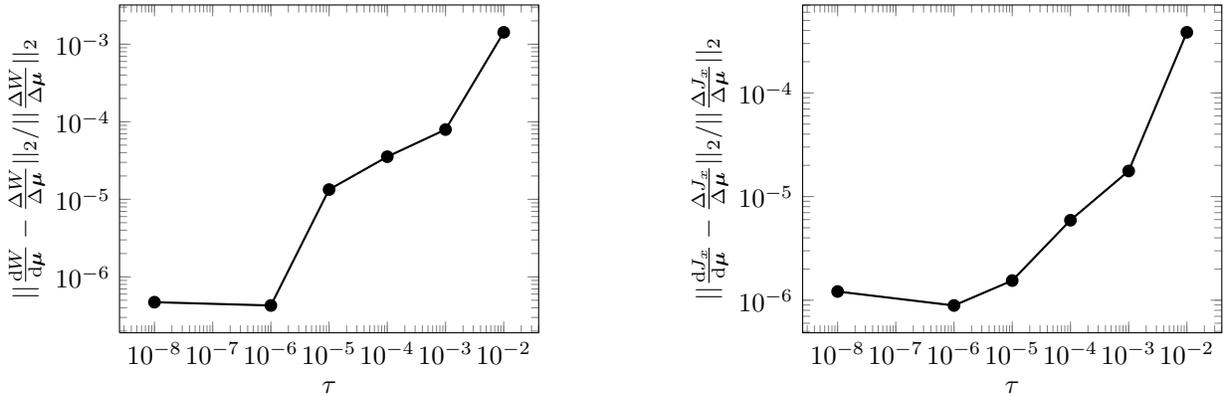

Finally, the adjoint method for computing gradients of quantities of
interest on the manifold of time-periodic solutions of the partial
differential equations is verified against a second-order finite difference
approximations. The finite difference approximation to gradients on the
aforementioned manifold requires finding the \emph{time-periodic} solution of
the governing equations \emph{at perturbations} about the nominal parameter
configuration in (\ref{eqn:nominal}). Figure~\ref{fig:findiff} and
Table~\ref{tab:findiff} show the
relative error between the gradients computed via the adjoint method in
Algorithm~\ref{alg:grad} to this finite difference approximation for a sweep
of finite difference intervals, $\tau$. To realize the sub-$10^{-6}$ finite
difference errors in the time-periodic gradient, tolerances of $10^{-12}$ were
used for the primal and dual time-periodic solutions. As expected, the error
starts to increase after $\tau$ drops too small due to the trade-off between
finite difference accuracy and round-off error.

\begin{table}
\centering
\begin{tabular}{l c c c c} \toprule
 & $\pder{W}{A_h}$ & $\pder{W}{\omega_h}$ & $\pder{W}{A_\theta}$
 & $\pder{W}{c_\theta}$
 \\\midrule 
 \rowcolor{black!00} Adjoint &
 -2.3091901647e+01 & -2.59357909093e+01 &
 -7.9956810720e+00 & 5.88159501747e-01 \\
 \rowcolor{black!10} Finite difference &
 -2.3091901345e+01 & -2.59357939542e+01 &
 -7.9956815124e+00 &  5.88159491776e-01 \\\toprule
 & $\pder{J_x}{A_h}$ & $\pder{J_x}{\omega_h}$ & $\pder{J_x}{A_\theta}$
 & $\pder{J_x}{c_\theta}$ \\\midrule
 \rowcolor{black!00} Adjoint &
 -1.8543679019e-01 & -1.02983075340e-01 &
  6.7297082256e+00 &  1.27010690706e-02 \\
 \rowcolor{black!10} Finite difference &
 -1.8543677450e-01 & -1.02983412630e-01 &
  6.7297089105e+00 &  1.27011295647e-02 \\\toprule
\end{tabular}

\caption{Comparison of non-zero derivatives of total energy, $W$, and
         $x$-impulse, $J_x$, computed with the adjoint method and a second-order
         finite difference approximation with step size $\tau = 10^{-6}$.}
\label{tab:findiff}
\end{table}

\subsection{Energetically Optimal Flapping with Thrust and Time-Periodicity
            Constraints} \label{subsec:app-opt}
In this section, the periodic adjoint method is used to solve an
optimal control problem with \emph{time-periodicity constraints} using
gradient-based optimization techniques. The optimization problem
is to determine the \emph{energetically optimal} flapping motion of the
NACA0012 airfoil in isentropic, viscous flow -- over a single
\emph{representative}, in-flight period -- such that the $x$-directed impulse
on the body is identically $0$. The continuous form of the optimal control
problem is given as
\begin{equation} \label{opt:cont}
  \begin{aligned}
    & \underset{\Ubm,~\mubold}{\text{minimize}}
    & & \Wcal(\Ubm, \mubold) \\
    & \text{subject to}
    & & \Jcal_x(\Ubm, \mubold) = 0 \\
    & & & \Ubm(\xbm,~0) = \Ubm(\xbm,~T) \\
    & & & \pder{\Ubm}{t} + \nabla \cdot \Fbm(\Ubm,~\nabla\Ubm) = 0 \quad
          \text{ in } \Omega(\mubold, t).
  \end{aligned}
\end{equation}
After spatial and temporal discretization via the high-order discontinuous
Galerkin and diagonally implicit Runge-Kutta schemes in Section~\ref{sec:app},
the continuous optimization problem in (\ref{opt:cont}) is replaced with its
fully discrete counterpart
\begin{equation} \label{opt:disc}
  \begin{aligned}
    & \underset{\substack{\ubm^{(0)},~\dots,~\ubm^{(N_t)} \in \Rbb^{N_\ubm},\\
                        \kbm_1^{(1)},~\dots,~\kbm_s^{(N_t)} \in \Rbb^{N_\ubm},\\
                        \mubold \in \Rbb^{N_\mubold}}}
                        {\text{minimize}}
    & & W(\ubm^{(0)},~\dots,~\ubm^{(N_t)},~\kbm_1^{(1)},~\dots,~\kbm_s^{(N_t)},~
          \mubold) \\
    & \text{subject to}
    & & J_x(\ubm^{(0)},~\dots,~\ubm^{(N_t)},~\kbm_1^{(1)},~\dots,
            ~\kbm_s^{(N_t)},~\mubold) = 0 \\
    & & &   \ubm^{(0)} = \ubm^{(N_t)} \\
    & & & \ubm^{(n)} = \ubm^{(n-1)} + \sum_{i = 1}^s b_i\kbm^{(n)}_i \\
    & & & \Mbb\kbm^{(n)}_i = \Delta t_n\rbm\left(\ubm_i^{(n)},~\mubold,~
                                                 t_{n-1} + c_i\Delta t_n\right).
  \end{aligned}
\end{equation}
The physical and numerical setup are identical to that in
Section~\ref{subsec:app-solvers} with the exception of the kinematic
parametrization. Instead of a single Fourier mode, the kinematic motion is
parametrized by cubic splines with $5$ equally spaced knots and boundary
conditions that enforce
\begin{equation} \label{eqn:mirror}
 \begin{aligned}
  h(\mubold,\,t) &= -h(\mubold,\,t+T/2) \\
  \theta(\mubold,\,t) &= -\theta(\mubold,\,t+T/2) 
 \end{aligned}
\end{equation} 
where $t$ is time and $T = 5$ is the fixed period of the flapping motion.
The vector of parameters, $\mubold$ -- used as optimization parameters --
are the \emph{knots} of the cubic splines. This leads to $N_\mubold = 8$
parameters; $4$ knots for the motion of $h(\mubold,\,t)$ and
$\theta(\mubold,\,t)$%
\footnote{There are only $4$ degrees of freedom since the mirror boundary
          condition in (\ref{eqn:mirror}) prescribes the value of one of the
          knots given the other four.}.
Notice that (\ref{eqn:mirror}) enforces the trajectories
of $h(\mubold,\,t)$ and $\theta(\mubold,\,t)$ in $[T/2,~T]$ to be the mirror
of those in $[0,~T/2]$, which implicitly enforces periodicity with period $T$.
The mapping $\Gcal$ from the reference to physical domain
              required for the DG-ALE formualtion is defined in
              (\ref{eqn:defmap}) with the new definition of $h(\mubold,\,t)$
              and $\theta(\mubold,\,t)$ with periodic cubic splines.

The optimization problem in (\ref{opt:disc}) is solved using the extension of
the nested framework for PDE-constrained optimization, or generalized
reduced-gradient method, introduced in Section~\ref{subsec:adj-grg}. The solvers
introduced in Section~\ref{subsec:tpde-soln} will be used to determine the
time-periodic flow around the airfoil. Given the results in the previous
section, the Newton-GMRES method with a tolerance of $\Delta = 10^{-3}$,
warm-started from $m = 5$ fixed-point iterations is employed. The flow is
deemed to be periodic if
\begin{equation} \label{eqn:periodic-tol}
 \norm{\ubm^{(0)} - \ubm^{(N_t)}}_2 \leq 10^{-10}.
\end{equation}
The periodic flow is used to compute quantities of
interest -- the total work and $x$-impulse. Then, the periodic adjoint method
will be used to compute gradients of the quantities of interest along the
manifold of \emph{time-periodic} solutions of the governing equation. GMRES
is used to solve the dual linear, periodic adjoint equations with a tolerance
of $\Delta = 10^{-4}$. Since there are two quantities of interest,
two periodic adjoint solves must be performed at each optimization iteration.
Finally, the quantities of interest and their gradients are passed to an
optimization solver -- SNOPT \cite{gill2002snopt} is used in this work -- and
progress is made toward a local minimum.

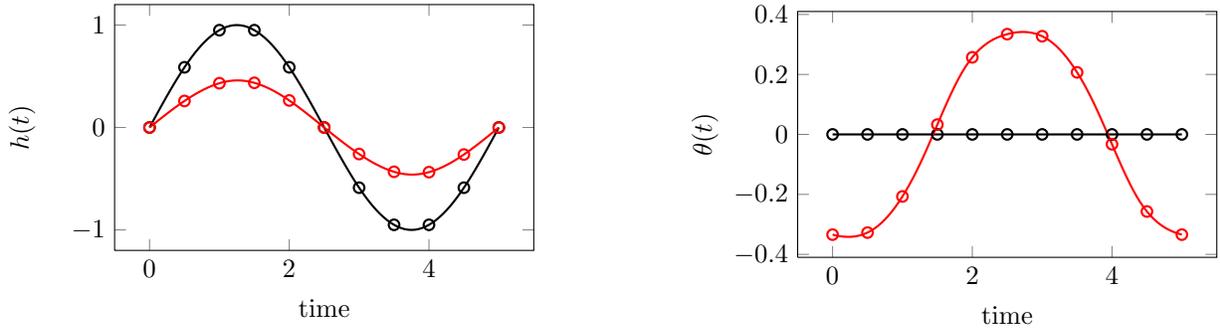
\begin{figure}
\begin{subfigure}{0.45\textwidth}
 \begin{tikzpicture}

\begin{axis}[
  scale only axis,
  width=0.75\textwidth,
  height=0.15\textheight,
  xlabel=time,
  ylabel=$h(t)$]
\addplot [black, solid, thick, mark=o, mark repeat=10]  table[x expr=0.05*\coordindex, y index=2] {foil2d/cfd/flap01/popt/0/testcase00/ascii/app0_mesh-how1-0-3_cfd-1p4-0p15-1000-0p72-0-1-_1-0__tdisc2-0p2-100-_1-1_-3_it000.xyth}; \label{line:traj-init}
\addplot [red, solid, thick, mark=o, mark repeat=10]  table[x expr=0.05*\coordindex, y index=2] {foil2d/cfd/flap01/popt/0/testcase00/ascii/app0_mesh-how1-0-3_cfd-1p4-0p15-1000-0p72-0-1-_1-0__tdisc2-0p2-100-_1-1_-3_it026.xyth}; \label{line:traj-opt}

\end{axis}

\end{tikzpicture}
\end{subfigure} \hfill
\begin{subfigure}{0.45\textwidth}
 \begin{tikzpicture}

\begin{axis}[
  scale only axis,
  width=0.75\textwidth,
  height=0.15\textheight,
  xlabel=time,
  ylabel=$\theta(t)$]
\addplot [black, solid, thick, mark=o, mark repeat=10]  table[x expr=0.05*\coordindex, y index=3] {foil2d/cfd/flap01/popt/0/testcase00/ascii/app0_mesh-how1-0-3_cfd-1p4-0p15-1000-0p72-0-1-_1-0__tdisc2-0p2-100-_1-1_-3_it000.xyth};
\addplot [red, solid, thick, mark=o, mark repeat=10]  table[x expr=0.05*\coordindex, y index=3] {foil2d/cfd/flap01/popt/0/testcase00/ascii/app0_mesh-how1-0-3_cfd-1p4-0p15-1000-0p72-0-1-_1-0__tdisc2-0p2-100-_1-1_-3_it026.xyth};
\end{axis}

\end{tikzpicture}
\end{subfigure}
\caption{Trajectories of $h(t)$ and $\theta(t)$ at initial guess
         (\ref{line:traj-init}) and optimal solution (\ref{line:traj-opt})
         for optimization problem in (\ref{opt:disc}).}
\label{fig:traj2}
\end{figure}

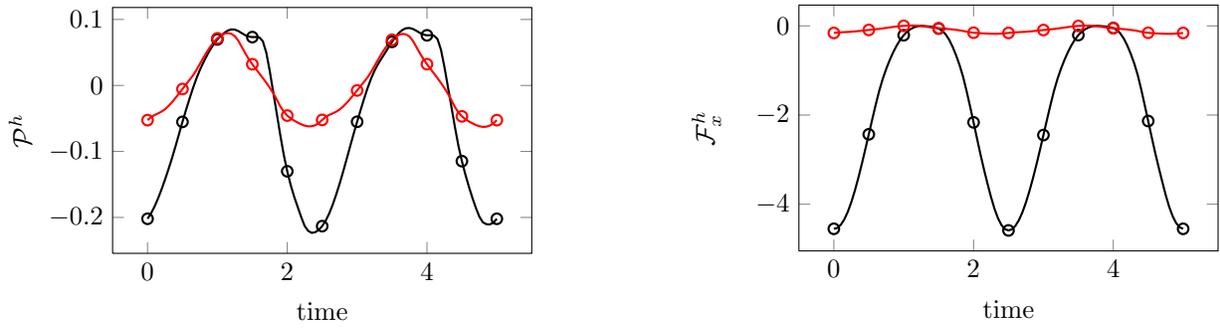
\begin{figure}
\begin{subfigure}{0.45\textwidth}
 \begin{tikzpicture}

\begin{axis}[
  scale only axis,
  width=0.75\textwidth,
  height=0.15\textheight,
  xlabel=time,
  ylabel={$\Pcal^h$}]
\addplot [black, solid, thick, mark=o, mark repeat=10]  table[x expr=0.05*\coordindex, y index=0] {foil2d/cfd/flap01/popt/0/testcase00/ascii/app0_mesh-how1-0-3_cfd-1p4-0p15-1000-0p72-0-1-_1-0__tdisc2-0p2-100-_1-1_-3_it000.allpost.instant}; \label{line:qoi-traj-init}
\addplot [red, solid, thick, mark=o, mark repeat=10]  table[x expr=0.05*\coordindex, y index=0] {foil2d/cfd/flap01/popt/0/testcase00/ascii/app0_mesh-how1-0-3_cfd-1p4-0p15-1000-0p72-0-1-_1-0__tdisc2-0p2-100-_1-1_-3_it026.allpost.instant}; \label{line:qoi-traj-opt}
\end{axis}

\end{tikzpicture}
\end{subfigure} \hfill
\begin{subfigure}{0.45\textwidth}
 \begin{tikzpicture}

\begin{axis}[
  scale only axis,
  width=0.75\textwidth,
  height=0.15\textheight,
  xlabel=time,
  ylabel={$\Fcal_x^h$}
]
\addplot [black, solid, thick, mark=o, mark repeat=10]  table[x expr=0.05*\coordindex, y index=3] {foil2d/cfd/flap01/popt/0/testcase00/ascii/app0_mesh-how1-0-3_cfd-1p4-0p15-1000-0p72-0-1-_1-0__tdisc2-0p2-100-_1-1_-3_it000.allpost.instant};
\addplot [red, solid, thick, mark=o, mark repeat=10]  table[x expr=0.05*\coordindex, y index=3] {foil2d/cfd/flap01/popt/0/testcase00/ascii/app0_mesh-how1-0-3_cfd-1p4-0p15-1000-0p72-0-1-_1-0__tdisc2-0p2-100-_1-1_-3_it026.allpost.instant};
\end{axis}
\end{tikzpicture}
\end{subfigure}
\caption{Time history of the power, $\Pcal^h(\ubm, \mubold, t)$, and $x$-directed
         force, $\Fcal_x^h(\ubm, \mubold, t)$, imparted onto foil by fluid at
         initial guess (\ref{line:qoi-traj-init}) and optimal solution
         (\ref{line:qoi-traj-opt}) for optimization problem in
         (\ref{opt:disc}).}
\label{fig:qoi-hist2}
\end{figure}

\begin{figure}
\begin{subfigure}{0.45\textwidth}
 \begin{tikzpicture}

\begin{axis}[
  scale only axis,
  width=0.75\textwidth,
  height=0.20\textheight,
  xlabel=optimization iteration,
  ylabel=$W$
]
\addplot [black, solid, thick, mark=o, mark repeat=1]  table[select coords between index={0}{17}, x index=0, y index=4] {foil2d/cfd/flap01/popt/0/testcase00/ascii/app0_mesh-how1-0-3_cfd-1p4-0p15-1000-0p72-0-1-_1-0__tdisc2-0p2-100-_1-1_-3.allpost.integrated.conv};
\end{axis}
\end{tikzpicture}
\end{subfigure} \hfill
\begin{subfigure}{0.45\textwidth}
 \begin{tikzpicture}

\begin{axis}[
  scale only axis,
  width=0.75\textwidth,
  height=0.20\textheight,
  xlabel=optimization iteration,
  ylabel=$J_x$
]
\addplot [black, solid, thick, mark=o, mark repeat=1]  table[select coords between index={0}{17}, x index=0, y index=1] {foil2d/cfd/flap01/popt/0/testcase00/ascii/app0_mesh-how1-0-3_cfd-1p4-0p15-1000-0p72-0-1-_1-0__tdisc2-0p2-100-_1-1_-3.allpost.integrated.conv};
\end{axis}
\end{tikzpicture}
\end{subfigure}
\caption{Convergence of quantities of interest, $W$ and $J_x$, with
         optimization iteration. Each optimization iteration requires the
         a periodic flow computation and its corresponding adjoint to
         evaluate the quantities of interest and their gradients.}
\label{fig:qoi-conv2}
\end{figure}
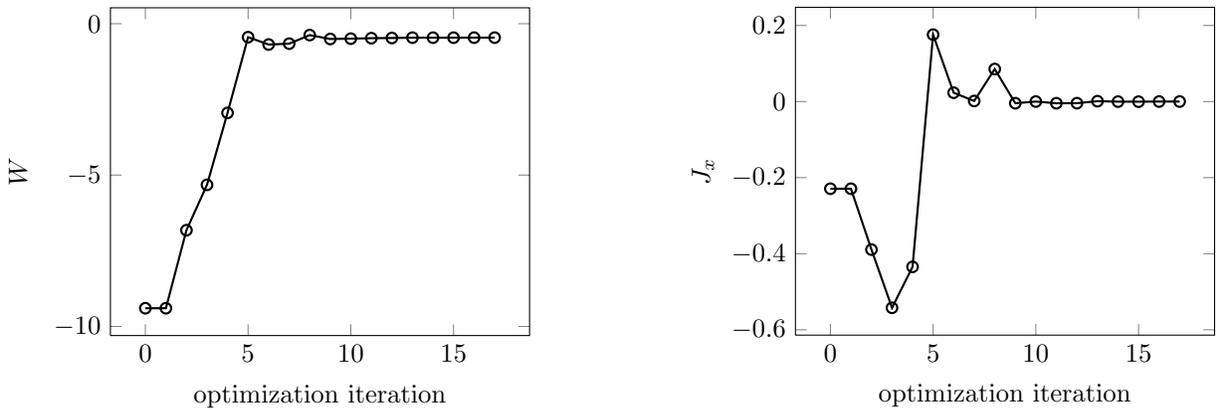

\begin{figure}
 \centering
 \begin{subfigure}{0.3\textwidth}
  \centering
  \includegraphics[width=\textwidth]{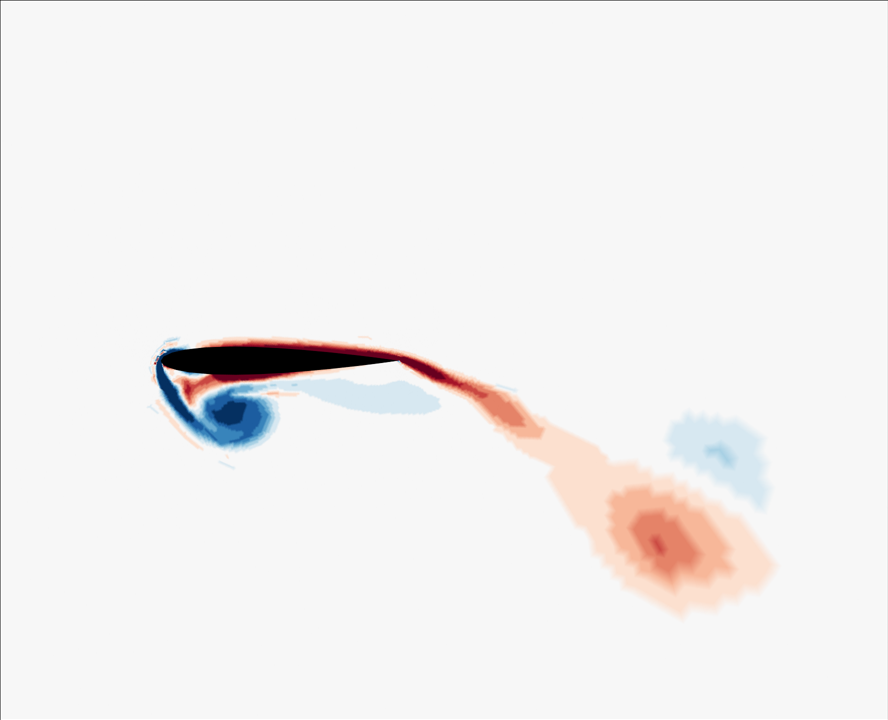}
 \end{subfigure} \hfill
 \begin{subfigure}{0.3\textwidth}
  \centering
  \includegraphics[width=\textwidth]{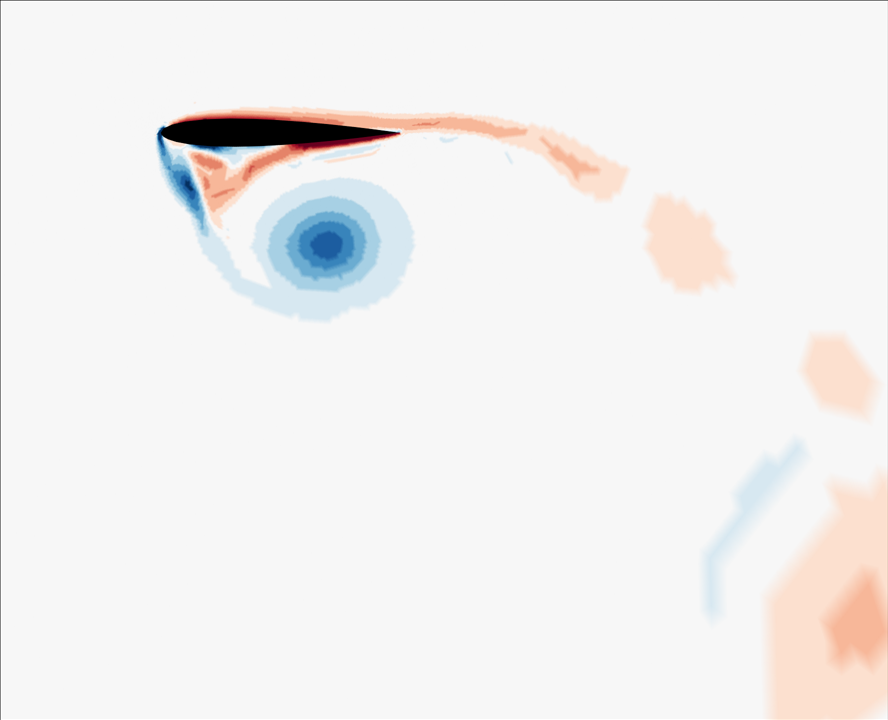}
 \end{subfigure} \hfill
 \begin{subfigure}{0.3\textwidth}
  \centering
  \includegraphics[width=\textwidth]{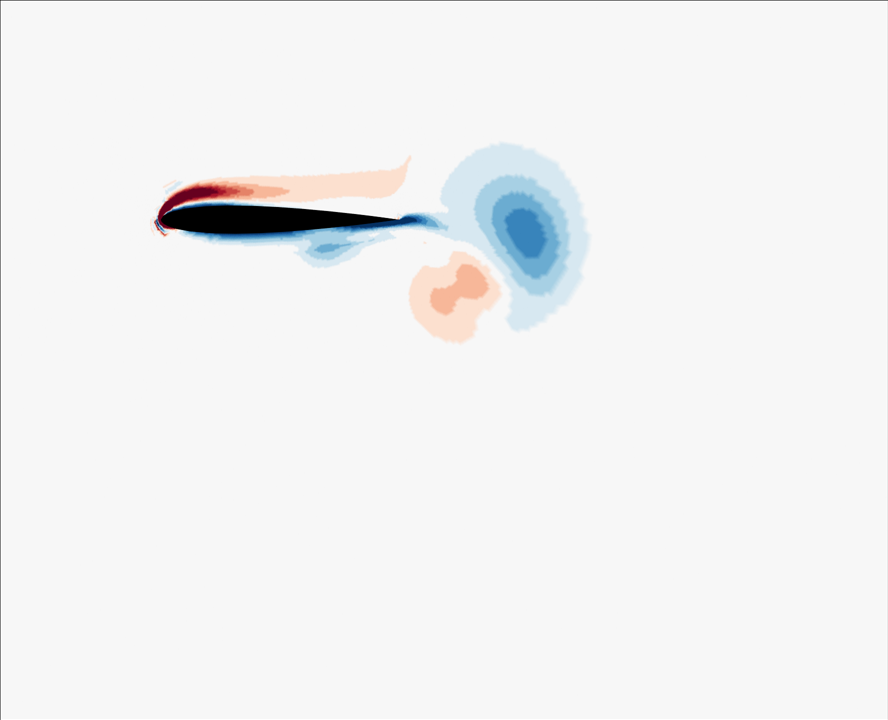}
 \end{subfigure} \\
 \begin{subfigure}{0.3\textwidth}
  \centering
  \includegraphics[width=\textwidth]{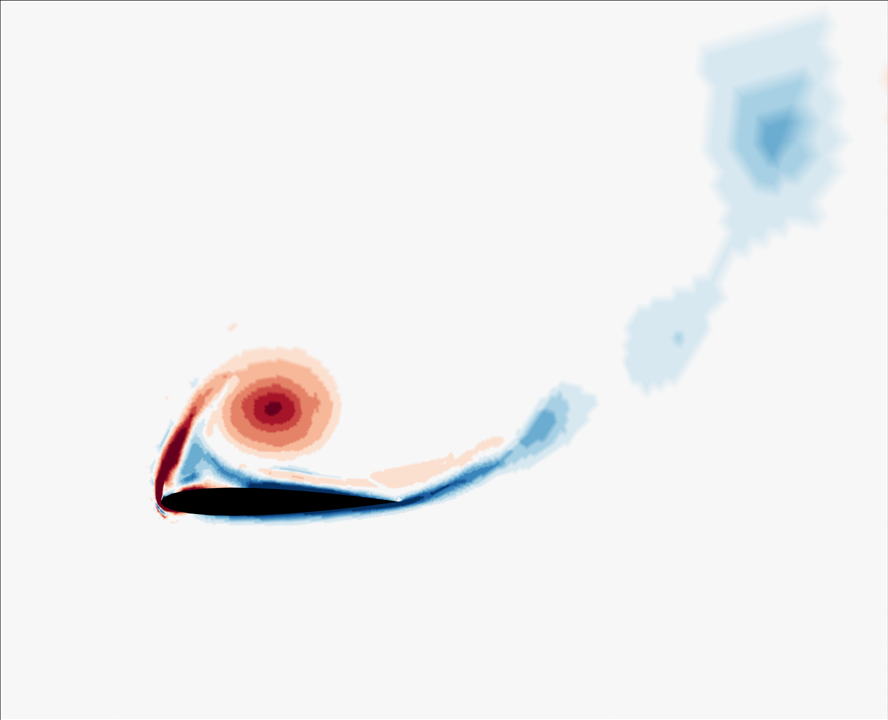}
 \end{subfigure} \hfill
 \begin{subfigure}{0.3\textwidth}
  \centering
  \includegraphics[width=\textwidth]{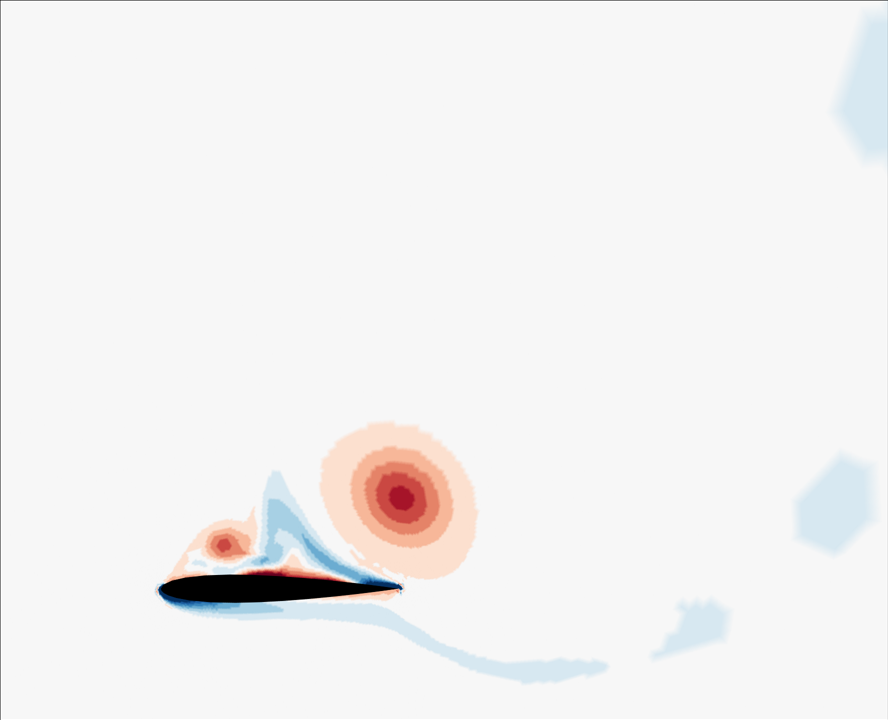}
 \end{subfigure} \hfill
 \begin{subfigure}{0.3\textwidth}
  \centering
  \includegraphics[width=\textwidth]{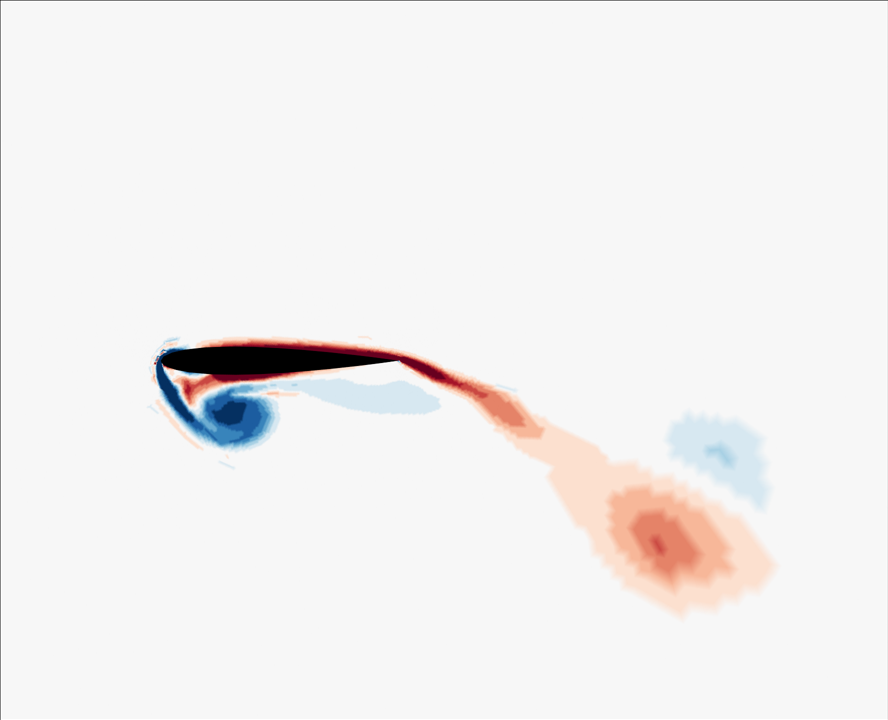}
 \end{subfigure}
 \caption{Trajectory of airfoil and flow vorticity at initial guess for
          optimization (pure heaving motion, see Figure~\ref{fig:traj2}).
          Snapshots taken at times $t = 0.0,~1.0,~2.0,~3.0,~4.0,~5.0$.}
 \label{fig:vort-init}
\end{figure}

\begin{figure}
 \centering
 \begin{subfigure}{0.3\textwidth}
  \centering
  \includegraphics[width=\textwidth]{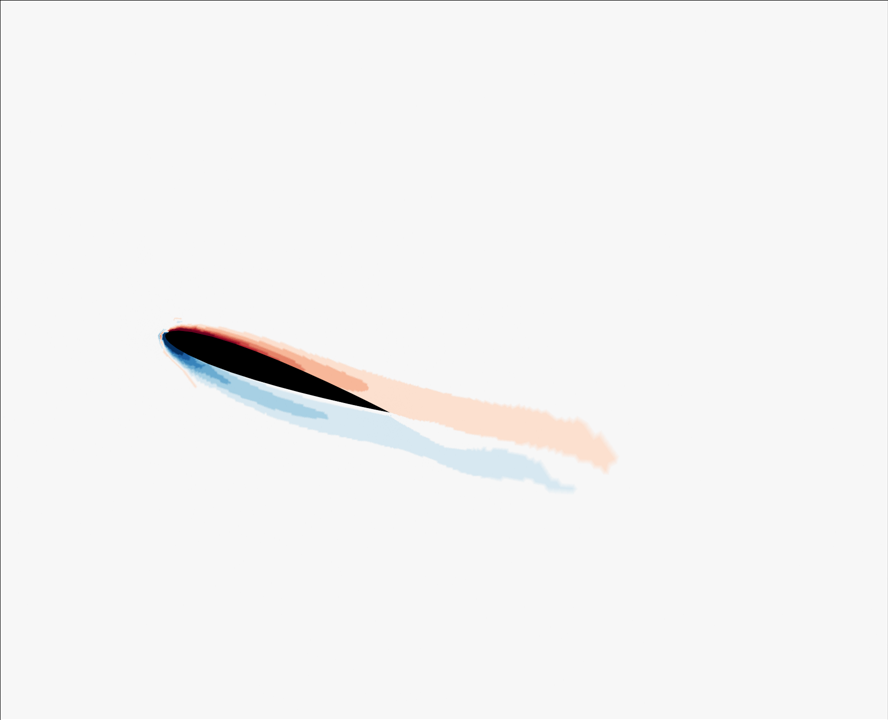}
 \end{subfigure} \hfill
 \begin{subfigure}{0.3\textwidth}
  \centering
  \includegraphics[width=\textwidth]{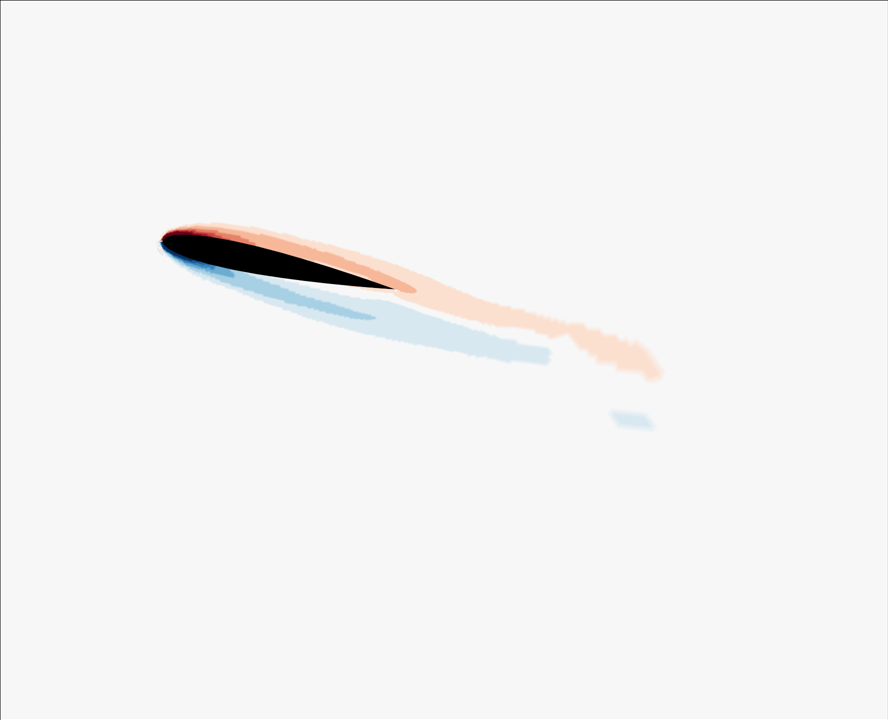}
 \end{subfigure} \hfill
 \begin{subfigure}{0.3\textwidth}
  \centering
  \includegraphics[width=\textwidth]{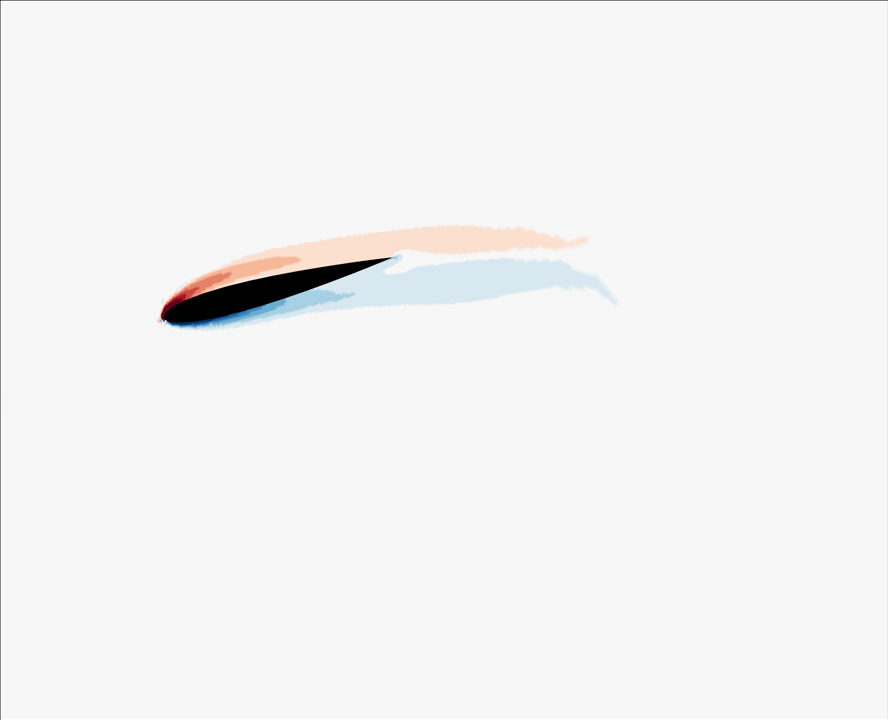}
 \end{subfigure} \\
 \begin{subfigure}{0.3\textwidth}
  \centering
  \includegraphics[width=\textwidth]{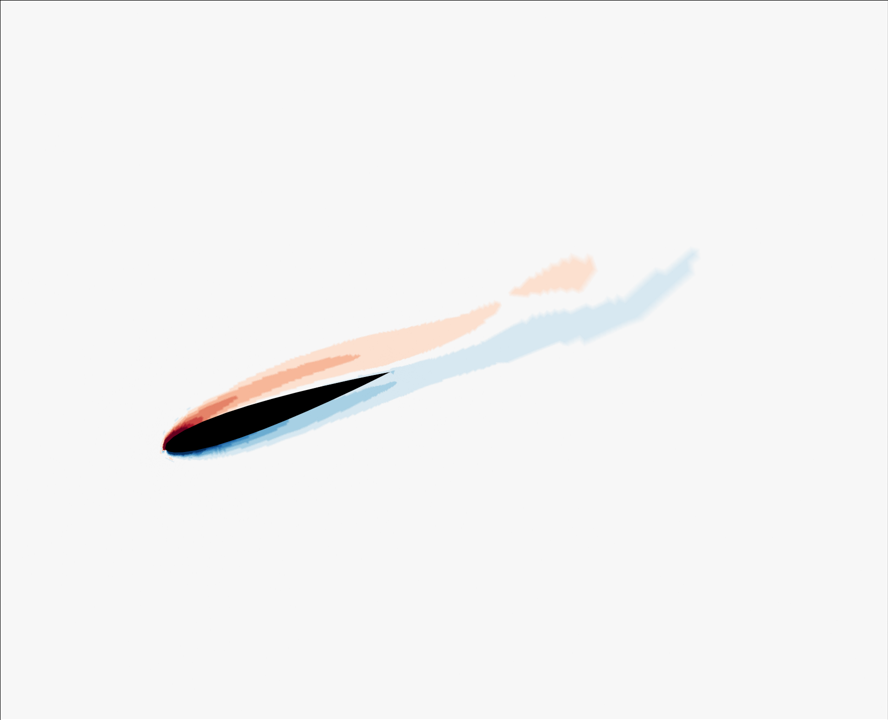}
 \end{subfigure} \hfill
 \begin{subfigure}{0.3\textwidth}
  \centering
  \includegraphics[width=\textwidth]{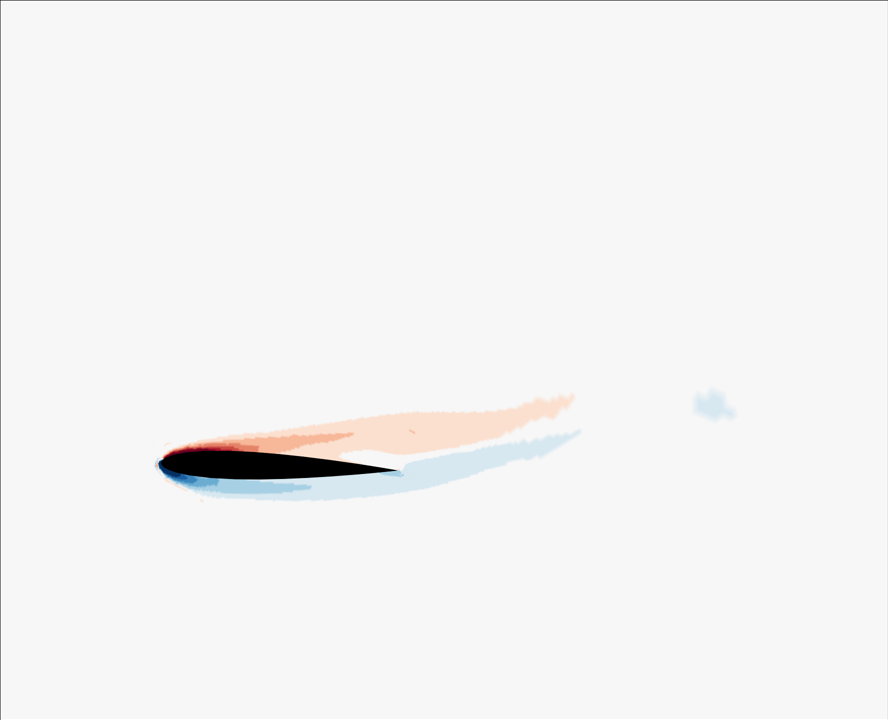}
 \end{subfigure} \hfill
 \begin{subfigure}{0.3\textwidth}
  \centering
  \includegraphics[width=\textwidth]{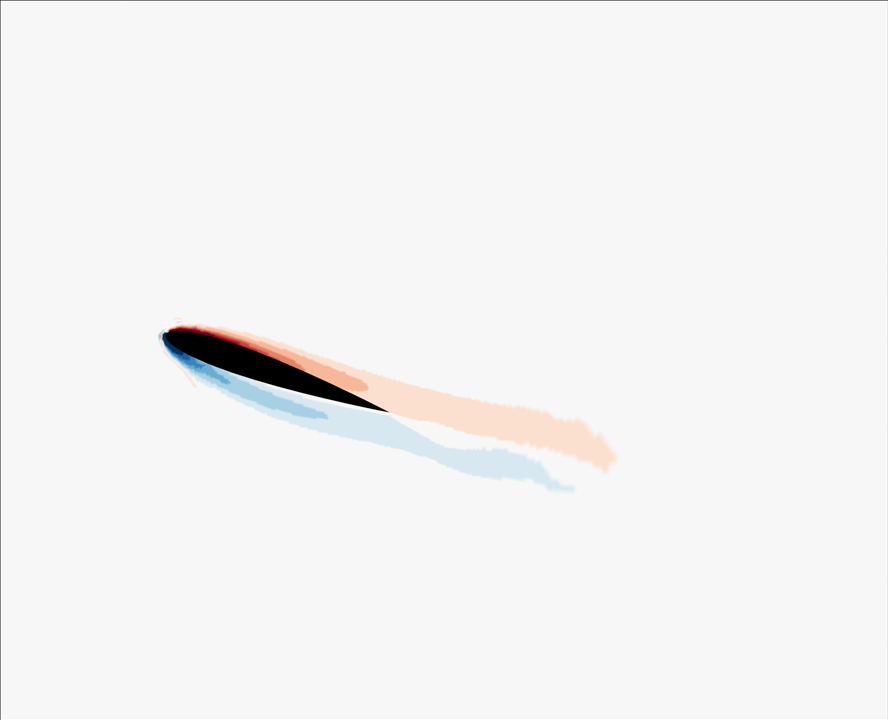}
 \end{subfigure}
 \caption{Trajectory of airfoil and flow vorticity at energetically optimal,
          zero-impulse flapping motion (see Figure~\ref{fig:traj2}).
          Snapshots taken at times $t = 0.0,~1.0,~2.0,~3.0,~4.0,~5.0$.}
 \label{fig:vort-opt}
\end{figure}

The initial condition for the optimization solver is shown in
Figure~\ref{fig:traj2}; the heaving motion is a sinusoid with amplitude $1$
and there is no pitch -- pure heaving motion. The vorticity snapshots in
Figure~\ref{fig:vort-init} show this motion induces a fairly violent flow with
shedding vortices. The corresponding time history of the power,
$\Pcal^h(\ubm, \mubold, t)$, and $x$-directed force,
$\Fcal_x^h(\ubm, \mubold, t)$, imparted onto the airfoil by the fluid are shown
in Figure~\ref{fig:qoi-hist2}.  After $16$ periodic optimization iterations, the
first-order optimality conditions have been reduced by two orders of
magnitude.  From Figure~\ref{fig:traj2}, the optimal airfoil motion is a
combination of heaving and pitching. From the initial guess, the amplitude of
the heaving motion has been reduced by more than a factor of two and the
pitching amplitude increased to $18.7^\circ$. The convergence history for the
optimization solver is given in Figure~\ref{fig:qoi-conv2}. At the optimal
solution, the total work required to perform the flapping motion is more than
an order of magnitude smaller than at the initial guess (pure heaving).
Figures~\ref{fig:vort-init}~and~\ref{fig:vort-opt} show snapshots of the flow
in time at the initial, purely heaving motion and the optimal flapping motion.
From these figures, it is clear that the flow corresponding to
the optimal motion is relatively benign with no shedding vortices, which
explains the reduction in required work. The efficiency of combined
pitching and heaving has been repeatedly observed experimentally
\cite{tuncer2005optimization, ramamurti2001simulation, platzer2008flapping} and
the phase angle of approximately $90^\circ$ between pitching and heaving,
as observed in Figure~\ref{fig:traj2}, has also been observed in experiments
\cite{tuncer2005optimization, ramamurti2001simulation, platzer2008flapping,
      oyama2009aerodynamic}.
The specific pitching and heaving amplitudes were determined by the
optimizer such that the thrust constraint is satisfied; if the thrust
requirement was increased, these magnitudes would increase and result in a
more violent flow field, eventually leading to vortex shedding
\cite{tuncer2005optimization, zahr2016dgopt}.

\section{Conclusion}
This document discussed a fully discrete framework for computing time-periodic
solutions of partial differential equations. The discussion included the
spatio-temporal discretization of the governing equations and a slew of
time-periodic shooting solvers, including optimization-based and Newton-Krylov
methods. These shooting methods consider the state at the final time to be a
nonlinear function of the initial condition and solve
$\ubm^{(N_t)}(\ubm_0) = \ubm_0$ using Newton-Raphson iterations or
optimization techniques to minimize its norm. The linear system of equations,
arising in the Newton-Raphson iterations, were solved using matrix-free GMRES
with matrix-vector products computed as the solution of the linearized,
sensitivity equations (with appropriate initial condition). The adjoint method
was used to compute the gradients in the gradient-based optimization solvers.
These periodic solvers were used to compute the time-periodic flow around a
flapping airfoil in isentropic, compressible, viscous flow, and their
performance compared. The Newton-Krylov solver exhibits superior convergence to
the optimization-based shooting methods, even when inexact tolerances were used
on the linear system solves, and fully leverages quality starting guesses.
An eigenvalue analysis is provided to show the periodic orbit of the flapping
problem is stable.

The main contribution of the document is the derivation of the adjoint
equations corresponding to the fully discrete time-periodically constraint
partial differential equations. As opposed to the backward-in-time evolution
equations, these equations constitute a linear, \emph{two-point boundary value
problem} that is provably solvable. The corresponding adjoint method was
introduced for computing \emph{exact} gradients of quantities of interest
along the manifold of time-periodic solutions of the discrete conservation law.
The gradients were verified against a second-order finite difference
approximation. These quantities of interest and their gradients were used in
the context of gradient-based optimization to solve an optimal control problem
with time-periodicity constraints, among others. In particular, the
energetically optimal flapping motion of a 2D airfoil in \emph{time-periodic},
isentropic, compressible, viscous flow that generates a prescribed
time-averaged thrust is sought. The proposed framework improves the nominal
flapping motion by reducing the flapping energy nearly an order of magnitude
and exactly satisfies the thrust constraint.

While this work is an initial step toward problems of engineering and scientific
relevance, additional development will be required to solve truly impactful
problems. One extension of this work is the development of robust solvers for
determining \emph{nearly} time-periodic solutions of problems where a
time-periodic solution does not exist, but exhibits quasi-cyclic behavior. An
example of such a problem is the 3D turbulent flow around periodically
driven bodies such as helicopter and windmill blades. Another extension will
be the development of faster numerical solvers to reduce the cost of computing
time-periodic solutions or solving optimization problems with time-periodicity
constraints. For example, economical, matrix-free preconditioners could
result in non-trivial speedups for the Newton-Krylov time-periodicity solver
and Krylov solver for the periodic adjoint equations. Model order reduction
techniques could dramatically reduce the cost of computing the solution of
the primal partial differential equations, and consequently the entire
time-periodic solver.


\appendix
\section{Existence and Uniqueness of Solutions of the Adjoint Equations
         of the Fully Discrete, Time-Periodically Constrained Partial
         Differential Equations}
\label{sec:exist}
This section proves existence and uniqueness of solutions of the adjoint
equations of the fully discrete, time-periodically constrained partial
differential equation. The strategy is to show the linear operator that
encapsulates them is the transpose of the linear operator that defines the
fully discrete, sensitivity equations, which is assumed non-singular at a
time-periodic solution.

Consider the initial-value problem (\ref{eqn:dirk2}), with the initial condition
parametrized by $\mubold$,
\begin{equation} \label{eqn:dirk2-app}
  \begin{aligned}
    \ubm^{(0)} &= \ubm_0(\mubold) \\
    \ubm^{(n)} &= \ubm^{(n-1)} + \sum_{i = 1}^s b_i\kbm^{(n)}_i \\
    \Mbb\kbm^{(n)}_i &= \Delta t_n\rbm\left(\ubm_i^{(n)},~\mubold,~
                                            t_{n-1} + c_i\Delta t_n\right).
  \end{aligned}
\end{equation}
The fully discrete adjoint equations corresponding to the primal equation in
(\ref{eqn:dirk2-app}) and the discrete quantity of
interest, $F(\ubm^{(0)}, \dots, \ubm^{(N_t)},
             \kbm_1^{(1)}, \dots, \kbm_s^{(N_t)}, \mubold)$
are
\begin{equation}
   \begin{aligned}
    \nubold^{(N_t)} &= \pder{F}{\ubm^{(N_t)}}^T \\
    \nubold^{(n-1)} &= \nubold^{(n)} +
                                 \pder{F}{\ubm^{(n-1)}}^T +
      \sum_{i=1}^s \Delta t_n\pder{\rbm}{\ubm}\left(\ubm_i^{(n)},~\mubold,~
                       t_{n-1}+c_i\Delta t_n\right)^T\taubold_i^{(n)} \\
    \Mbb^T\taubold_i^{(n)} &= \pder{F}{\kbm_i^{(n)}}^T +
                         b_i\nubold^{(n)} +
                         \sum_{j=i}^s a_{ji}\Delta t_n\pder{\rbm}{\ubm}
                         \left(\ubm_j^{(n)},
                         ~\mubold,~t_{n-1}+c_j\Delta t_n\right)^T
                         \taubold_j^{(n)},
  \end{aligned}
\end{equation}
and the gradient of the quantity of interest can be reconstructed as
\begin{equation}
   \oder{F}{\mubold} = \pder{F}{\mubold} +
                     {\nubold^{(0)}}^T\pder{\ubm_0}{\mubold} +
                     \sum_{n=1}^{N_t} \Delta t_n \sum_{i=1}^s
                 {\taubold_i^{(n)}}^T\pder{\rbm}{\mubold}(\ubm_i^{(n)},~
                                                \mubold,~t_{n-1}+c_i\Delta t_n),
\end{equation}
where $\nubold^{(n)}$ and $\taubold_i^{(n)}$ are the Lagrange multipliers.
These equations can be obtained using an identical derivation to that in
Section~\ref{subsec:adj-deriv}; see \cite{zahr2016dgopt}.
At this point, take $F = \vbm^T\ubm^{(N_t)}$ and $\mubold = \ubm_0$ for a fixed,
arbitrary vector $\vbm \in \Rbb^{N_\ubm}$. For this selection of $F$ and
$\mubold$, the above equations reduce to
\begin{equation} \label{eqn:adj-app}
   \begin{aligned}
    \nubold^{(N_t)} &= \vbm \\
    \nubold^{(n-1)} &= \nubold^{(n)} +
      \sum_{i=1}^s \Delta t_n\pder{\rbm}{\ubm}\left(\ubm_i^{(n)},~\mubold,~
                             t_{n-1}+c_i\Delta t_n\right)^T\taubold_i^{(n)} \\
    \Mbb^T\taubold_i^{(n)} &= b_i\nubold^{(n)} +
                         \sum_{j=i}^s a_{ji}\Delta t_n\pder{\rbm}{\ubm}
                         \left(\ubm_j^{(n)},
                         ~\mubold,~t_{n-1}+c_j\Delta t_n\right)^T
                         \taubold_j^{(n)}
  \end{aligned}
\end{equation}
and
\begin{equation}
   \oder{F}{\mubold}^T = \pder{\ubm^{(N_t)}}{\ubm_0}^T\vbm
                       = \nubold^{(0)}.
\end{equation}
The equations in (\ref{eqn:adj-app}) defining $\nubold^{(0)}$ are
\emph{identical} to those in (\ref{eqn:uns-disc-adj-dirk2}) defining
$\displaystyle{\pder{\lambdabold^{(0)}}{\lambdabold_{N_t}}}\vbm$, which leads
to the relation
\begin{equation}
 \pder{\ubm^{(N_t)}}{\ubm_0}^T\vbm =
 \pder{\lambdabold^{(0)}}{\lambdabold_{N_t}}\vbm
\end{equation}
for any $\vbm$. Thus, it can be concluded that
\begin{equation}
  \pder{\lambdabold^{(0)}}{\lambdabold_{N_t}} = \pder{\ubm^{(N_t)}}{\ubm_0}^T.
\end{equation}
Since the Jacobian of the time-periodic residual,
$\displaystyle{\pder{\ubm^{(N_t)}}{\ubm_0} - \Ibm}$, is non-singular at a
time-periodic solution, the matrix defining the linear, two-point boundary
value problem,
$\displaystyle{\pder{\lambdabold^{(0)}}{\lambdabold_{N_t}} - \Ibm}$
must also be non-singular. Thus, a solution of the linear, two-point boundary
value problem exists and is unique.


\section*{Acknowledgments}
This work was supported in part by the Department of Energy Computational
Science Graduate Fellowship Program of the Office of Science and National
Nuclear Security Administration in the Department of Energy under contract
DE-FG02-97ER25308 (MZ), and by the Director, Office of
Science, Computational and Technology Research, U.S. Department of
Energy under contract number DE-AC02-05CH11231 (PP and JW).
The content of this publication does not necessarily
reflect the position or policy of any of these supporters, and no official
endorsement should be inferred.

\bibliographystyle{ieeetr}
\bibliography{biblio}
\end{document}